\documentclass[10pt]{amsart}
\usepackage{amsbsy,amssymb,amscd,amsfonts,latexsym,amstext,delarray,
amsmath,amsthm, graphicx}

\input xypic
\usepackage{hyperref}

\newtheorem{The}{Theorem}[section]
\newtheorem{Cor}[The]{Corollary}
\newtheorem{Pro}[The]{Proposition}
\newtheorem{Exa}[The]{Example}
\newtheorem{Rem}[The]{Remark}
\newtheorem{Lem}[The]{Lemma}
\newtheorem{Def}[The]{Definition}
\newtheorem{Con}[The]{Problem}

\numberwithin{equation}{section}

\def\qqq{\,,\quad \forall\,}

\def\Ac{{\mathcal A}}

\def\Mc{{\mathcal M}}
\def\Sc{{\mathcal S}}
\def\wh{\widehat}
\def\wt{\widetilde}

\def\Oc{{\mathcal O}}

\def\proof{\vspace{2ex}\noindent{\bf Proof.} }
\def\endproof{\relax\ifmmode\expandafter\endproofmath\else
\unskip\nobreak\hfil\penalty50\hskip.75em\hbox{}\nobreak\hfil\bull
{\parfillskip= 0pt \finalhyphendemerits= 0 \bigbreak}\fi}
\def\endproofmath$${\eqno\bull$$\bigbreak}
\def\bull{\vbox{\hrule\hbox{\vrule\kern3pt\vbox{\kern6pt}\kern3pt\vrule}
\hrule}}

\newcommand{\ba}{\begin{eqnarray}}
\newcommand{\na}{\end{eqnarray}}

\newcommand{\scr}{\mathcal}

\newcommand{\A}{\mathbb{A}}
\newcommand{\C}{\mathbb{C}}
\newcommand{\E}{\mathbb{E}}
\newcommand{\bG}{\mathbb{G}}
\renewcommand{\H}{\mathbb{H}}

\newcommand{\F}{\mathbb{F}}
\newcommand{\K}{\mathbb{K}}
\renewcommand{\L}{\mathbb{L}}

\newcommand{\N}{\mathbb{N}}
\renewcommand{\P}{\mathbb{P}}
\newcommand{\Q}{\mathbb{Q}}
\newcommand{\R}{\mathbb{R}}

\newcommand{\Z}{\mathbb{Z}}

\newcommand{\cA}{\scr{A}}
\newcommand{\cB}{\scr{B}}

\newcommand{\cC}{\scr{C}}
\renewcommand{\cD}{\scr{D}}
\newcommand{\cE}{\scr{E}}
\newcommand{\cF}{\scr{F}}
\newcommand{\cG}{\scr{G}}
\renewcommand{\cH}{\scr{H}}

\newcommand{\cJ}{\scr{J}}
\newcommand{\cK}{\scr{K}}
\renewcommand{\cL}{\scr{L}}
\newcommand{\cM}{\scr{M}}
\newcommand{\cN}{\scr{N}}
\renewcommand{\O}{\scr{O}}
\newcommand{\cO}{\scr{O}}

\newcommand{\cS}{\scr{S}}

\newcommand{\cV}{\scr{V}}

\newcommand{\cX}{\scr{X}}

\newcommand{\cZ}{\scr{Z}}

\newcommand{\ie}{{\it i.e.\/}\ }
\newcommand{\eg}{{\it e.g.\/}\ }
\newcommand{\cf}{{\it cf.\/}\ }

\newcommand{\Spec}{{\rm Spec}}
\newcommand{\Sp}{{\rm Spec}}
\newcommand{\Ker}{{\rm Ker}}

\newcommand{\Gal}{{\rm Gal}}

\newcommand{\Tr}{{\rm Trace}}
\newcommand{\Trace}{{\rm Trace}}
\newcommand{\Aut}{{\rm Aut}}
\newcommand{\Inn}{{\rm Inn}}
\newcommand{\Irr}{{\rm Irr}}
\newcommand{\Out}{{\rm Out}}
\newcommand{\End}{{\rm End}}
\newcommand{\Ext}{{\rm Ext}}
\newcommand{\Tor}{{\rm Tor}}
\newcommand{\Hol}{{\rm Hol}}
\newcommand{\Hom}{{\rm Hom}}
\def\Prim{{\rm Prim}}

\newcommand{\Corr}{{\rm Corr}}

\title[NCG and motives]
{Noncommutative geometry and motives: \\ the thermodynamics of endomotives}

\author[Connes]{Alain Connes}
\author[Consani]{Caterina Consani}
\author[Marcolli]{Matilde Marcolli}
\address{A.~Connes: Coll\`ege de France \\
3, rue d'Ulm \\ Paris, F-75005 France
\\ I.H.E.S. and Vanderbilt
University} \email{alain\@@connes.org}
\address{C.~Consani: Mathematics Department \\ Johns Hopkins
University \\ Baltimore, MD 21218 USA} \email{kc\@@math.jhu.edu}
\address{M.~Marcolli: Max--Planck Institut f\"ur Mathematik  \\
Vivatsgasse 7 \\
Bonn, D-53111 Germany} \email{marcolli\@@mpim-bonn.mpg.de}

\begin{document}

\maketitle

\tableofcontents

\newpage

\section{Introduction}

A few unexpected encounters between noncommutative geometry and the
theory of motives have taken place recently. A first instance
occurred in \cite{AC}, where the Weil explicit formulae acquire the
geometric meaning of a Lefschetz trace formula over the
noncommutative space of adele classes. The reason why the
construction of \cite{AC} should be regarded as motivic is twofold.
On the one hand, as we discuss in this paper, the adele class space
is obtained from a noncommutative space, the Bost--Connes system,
which sits naturally in a category of noncommutative spaces
extending the category of Artin motives. Moreover, as we discuss
briefly in this paper and in full detail in our forthcoming work
\cite{CoCM}, it is possible to give a cohomological interpretation
of the spectral realization of the zeros of the Riemann zeta
function of \cite{AC} on the cyclic homology of a noncommutative
motive. The reason why this construction takes place in a category
of (noncommutative) motives is that the geometric space we need to
use is obtained as a cokernel of a morphism of algebras, which only
exists in a suitable abelian category extending the non-additive
category of algebras, exactly as in the context of motives and
algebraic varieties.

There are other instances of interactions between noncommutative
geometry and motives, some of which will in fact involve
noncommutative versions of more complicated categories of motives,
beyond Artin motives, involving pure and sometimes mixed motives.
One example is the problem of extending the results of \cite{AC} to
$L$-functions of algebraic varieties. The latter are, in fact,
naturally associated not to a variety but to a motive. Already in
\cite{ConsMar} it was shown that noncommutative geometry can be used
to model the geometry of the fibers at the archimedean places of
arithmetic varieties. This suggested the existence of a Lefschetz
trace formula for the local $L$-factors, and at least a semi-local
version for the $L$-function, over a noncommutative space obtained
as a construction over the adele class space. In this paper we give
such a Lefschetz interpretation for the archimedean local factors
introduced by Serre \cite{Se}. In order to obtain a similar
Lefschetz formula for the case of several places (archimedean and
non-archimedean), it seems necessary to develop a theory that
combines noncommutative spaces and motives.

A further example of an intriguing interaction between
noncommutative geometry and motives appears in the context of
perturbative renormalization in quantum field theory. The results of
\cite{cmln}, \cite{CM2} show the existence of a universal group of
symmetries that governs the structure of the divergences in
renormalizable quantum field theories and contains the
renormalization group as a 1-parameter subgroup. This universal
group is a motivic Galois group for a category of mixed Tate
motives. On the other hand, the treatment of divergent Feynman
integrals by dimensional regularization, used in this geometric
reformulation of renormalization, is best described in the framework
of noncommutative geometry. In fact, using the formalism of spectral
triples in noncommutative geometry, one can construct (\cf
\cite{CManom}) spaces with dimension a complex number $z$ and
describe geometrically the procedure of dimensional regularization
as a cup product of spectral triples. Once again it is desirable to
develop a common framework for noncommutative geometry and motives,
so as to obtain a more direct identification between the data of
perturbative quantum field theory and the objects of the category of
mixed motives governed by the same Galois group.

In a different but not unrelated perspective, Kontsevich recently
found another very interesting link between noncommutative geometry
and motives \cite{Maxim}, see also \cite{Kal}.

\medskip

The present paper is structured as follows. We begin by discussing
the problem of morphisms for noncommutative spaces. It is well
known, in fact, that morphisms of algebras are not sufficient and
one needs a larger class of morphisms that, for instance, can
account for the phenomenon of Morita equivalence.  We argue that
there are, in fact, two well developed constructions in
noncommutative geometry that make it possible to define morphisms as
correspondences, namely Kasparov's bivariant $KK$-theory on one hand
and modules over the cyclic category (and cyclic (co)homology) on
the other.

\smallskip

In $KK$-theory, one extends the notion of morphism to that of
(virtual) correspondences given by formal differences of
correspondences defined using bimodules. This yields an {\em
additive} category containing the category of $C^*$-algebras. The
product (composition of correspondences) plays a central role in the
understanding of $K$-theory of $C^*$-algebras. In the approach via
cyclic cohomology,  one constructs a functor $\cA\to \cA^\natural$
from the category of algebras and algebra homomorphisms to the {\em
abelian} category of $\Lambda$-modules introduced in \cite{CoExt}.
Cyclic cohomology  and homology appear then as derived functors,
namely as the $\Ext$ and $\Tor$ functors, much as in the case of the
absolute cohomology in the theory of motives. Working with the
resulting abelian category has the advantage of making all the
standard tools of  homological algebra available in noncommutative
geometry. It is well known that one obtains this way a good de Rham
(or crystalline) theory for noncommutative spaces, a general
bivariant cyclic theory, and  a framework for the Chern-Weil theory
of characteristic classes for general Hopf algebra actions
\cite{CMo}.

\smallskip

We show, first in the zero-dimensional case of Artin motives and
then in a more general setting in higher dimension, that the notion
of morphisms as (virtual) correspondences given by elements in
$KK$-theory is compatible with the notion of morphisms as used in
the theory of motives, where correspondences are given by algebraic
cycles in the product. In fact, we show  correspondences in the
sense of motives give rise to elements in $KK$-theory (using Baum's
geometric correspondences \cite{Baum} and \cite{CoSka}), compatibly
with the operation of composition.

\medskip

We focus on the zero-dimensional case of Artin motives. We extend
the category $\cC\cV^0_{\K,\E}$ of Artin motives over a field $\K$
with coefficients in a field $\E$ to a larger (pseudo)abelian
category $\cE\cV^0_{\K,\E}$ of {\em algebraic endomotives}, where
the objects are semigroup actions on projective systems of Artin
motives. The objects are described by semigroup crossed product
algebras $\cA_\K=A\rtimes S$ over $\K$ and the correspondences
extend the usual correspondences of Artin motives given by algebraic
cycles, compatibly with the semigroup actions.

We prove then that, in the case of a number field $\K$, taking
points over $\bar\K$ in the projective limits of Artin motives gives
a natural embedding of the category of algebraic endomotives into a
(pseudo)abelian category $C^*\cV^0_\K$ of {\em analytic
endomotives}. The objects are noncommutative spaces given by crossed
product $C^*$-algebras and the morphisms are given by geometric
correspondences, which define a suitable class of bimodules. The
Galois group $G$ acts as natural transformations of the functor
$\cF$ realizing the embedding. The objects in the resulting category
of analytic endomotives are typically noncommutative spaces obtained
as quotients of semigroup actions on totally disconnected compact
spaces.

We show that important examples like the Bost--Connes (BC) system of
\cite{BC} and its generalization for imaginary quadratic fields
considered in \cite{CMR} belong to this class. We describe a large
class of examples obtained from iterations of self-maps of algebraic
varieties.

\smallskip

We also show that the category of noncommutative Artin motives
considered here is especially well designed for the general program
of applications of noncommutative geometry to abelian class field
theory proposed in \cite{CM}, \cite{CMR2}. In fact, a first result
in the present paper shows that what we called ``fabulous states''
property in \cite{CM}, \cite{CMR2} holds for any noncommutative
Artin motive that is obtained from the action of a semigroup $S$ by
endomorphisms on a direct limit $A$ of finite dimensional algebras
that are products of abelian normal field extensions of a number
field $\K$. Namely, we show that there is a canonical action of the
absolute Galois group $G$ of $\K$ by automorphisms of a
$C^*$-completion $\bar\cA$ of the noncommutative space
$\cA_\K=A\rtimes S$ and, for pure states $\varphi$ on this
$C^*$-algebra which are induced from $A$, one has the intertwining
property
$$ \alpha \varphi(a)= \varphi(\alpha(a)) $$
for all $a\in A\rtimes S$ and all $\alpha\in G/[G,G]$.

\medskip

When passing from commutative to noncommutative algebras, new tools
of thermodynamical nature become available and make it possible to
construct an action of $\R^*_+$ on cyclic homology. We show that,
when applied to the noncommutative space (analytic endomotive)
$\cF(M)$ associated to an algebraic endomotive $M$, this
representation of $\R^*_+$ combines with a representation of the
Galois group $G$. In the particular case of the endomotive
associated to the BC system, the resulting representation of
$G\times \R^*_+$ gives the spectral realization of the zeros of the
Riemann zeta function and of the Artin $L$-functions for abelian
characters of $G$. One sees in this example that this construction
plays a role analogous to the action of the Weil group on the $\ell$-adic
cohomology and it can be thought of as a functor from the category
of endomotives to the category of representations of the group
$G\times \R^*_+$. In fact, here we think of the action of $\R^*_+$ as a
``Frobenius in characteristic zero'', hence of $G\times \R^*_+$ as
the corresponding Weil group.

Notice that it is the type III nature of the BC system that is at
the root of the thermodynamical behavior that makes the procedure
described above nontrivial. In particular, this procedure only
yields trivial results if applied to commutative and type II cases.

The construction of the appropriate ``motivic cohomology'' with the
``Frobenius'' action of $\R^*_+$ for endomotives is obtained through
a very general procedure, which we describe here in the generality
that suffices to the purpose of the present paper. It consists of
three basic steps, starting from the data of a noncommutative
algebra $\cA$ and a state $\varphi$, under specific assumptions that
will be discussed in detail in the paper. One considers the time
evolution $\sigma_t \in \Aut \cA$, $t\in \R$ naturally associated to
the state $\varphi$ (as in \cite{Co-th}).

\smallskip

The first step is what we refer to as {\em cooling}.  One considers
the space $\cE_\beta$ of extremal  KMS$_\beta$ states, for $\beta$
greater than critical. Assuming these states are of type I, one
obtains a morphism
$$
\pi : \cA \rtimes_{\sigma}\R \to \cS(\cE_\beta \times \R_+^*)\otimes
\cL^1,
$$
where $\cA$ is a dense subalgebra of a $C^*$-algebra $\bar\cA$, and
where $\cL^1$ denotes the ideal of trace class operators. In fact,
one considers this morphism restricted to the kernel of the
canonical trace $\tau$ on $\bar\cA \rtimes_{\sigma}\R$.

\smallskip

The second step is {\em distillation}, by which we mean the
following. One constructs the $\Lambda$-module $D(\cA,\varphi)$
given by the cokernel  of the cyclic morphism given by the
composition of $\pi$ with the trace $\Tr : \cL^1 \to \C$.

\smallskip

The third step is the {\em dual action}. Namely, one looks at the
spectrum of the canonical action of $\R_+^*$ on the cyclic homology
$$
HC_0(D(\cA,\varphi)).
$$

\smallskip

This procedure is quite general and applies to a large class of data
$(\cA,\varphi)$ and produces spectral realizations of zeros of
$L$-functions. In particular, we can apply it to the analytic
endomotive describing the BC-system. In this case, the
noncommutative space is that of commensurability classes of
$\Q$-lattices up to scaling. The latter can be obtained from the
action of the semigroup $\N$ on $\hat\Z$. However, the algebraic
endomotive corresponding to this action is not isomorphic to the
algebraic endomotive of the BC system, which comes from the algebra
$\Q[\Q/\Z]$. In fact, they are distinguished by the Galois action.
This important difference is a reflection of the distinction between
Artin and Hecke $L$-functions.

In this BC system, the corresponding quantum statistical system is
naturally endowed with an action of the group $G^{ab}=\,G/[G,G]$,
where $G$ is the Galois group of $\bar \Q$. Any Hecke character
$\chi$ gives rise to an idempotent $p_\chi$ in the ring of
endomorphisms of the $0$-dimensional noncommutative motive $X$ and
when applied to the range of $p_\chi$ the above procedure gives as a
spectrum the zeros of the Hecke $L$-function $L_\chi$.

\smallskip

We also show that the ``dualization" step, \ie the transition from
$\cA$ to $\cA \rtimes_{\sigma}\R$, is a very good analog in the case
of number fields of what happens for a function field $K$ in passing
to the extension $K\otimes_{\F_q} \bar \F_q$. In fact, in the case
of positive characteristic, the unramified extensions
$K\otimes_{\F_q}\F_{q^n}$ combined with the notion of places yield
the points $C(\bar \F_q)$ over $\bar \F_q $ of the underlying curve.
This has a good parallel in the theory of factors, as we are going
to discuss in this paper. In our forthcoming work \cite{CoCM} this
analogy will be developed further and will play an important role in
developing a setting in noncommutative geometry that parallels the
algebro-geometric framework that Weil used in his proof of RH for
function fields.

In the results discussed here and in \cite{CoCM}, instead of working
with the Hilbert space context of \cite{AC} the trace formula
involved in the spectral realization of \cite{AC} is given in
suitable nuclear function spaces associated to algebras of rapidly
decaying functions, as in the approach of \cite{Meyer}.

\smallskip

We also give an overview of some results that will be proved in full
detail in \cite{CoCM}, where we give a cohomological interpretation
of the spectral realization of the zeros of the zeta function as the
cyclic homology of a ``suitable" noncommutative space $M_\K$.
The sought for space $M_\K$ appears naturally as the cokernel of the
restriction morphism of functions on the adele class space $X_\K$ to
functions on the ``cooled down'' subspace $C_\K$ of idele classes.
This operation would not make sense in the usual category of
algebras (commutative or not) and algebra homomorphisms, since this
category is not additive and one cannot take cokernels, whereas this
can be done in a suitable category of motives.

\medskip

{\bf Acknowledgement.} This research was partially supported by the
third author's Sofya Kovalevskaya Award and by the second author's
NSERC grant 7024520. Part of this work was done during a visit of
the first and third authors to the Kavli Institute in Santa Barbara,
supported in part by the National Science Foundation under Grant No.
PHY99-07949, and during a visit of the first two authors to the Max
Planck Institute.

\bigskip

\bigskip
\section{Correspondences for
noncommutative spaces}\label{Scyclic}

It is often regarded as a problem in noncommutative geometry to
provide a good notion of {\em morphisms} for noncommutative spaces,
which accounts, for instance, for phenomena such as Morita
equivalence. Moreover, it is often desirable to apply to
noncommutative geometry the tools of homological algebra, which
require working with abelian categories. The category of algebras
with algebra homomorphisms is not satisfactory for both of the
reasons just mentioned. We argue here that the best approach to
dealing with a category of noncommutative spaces is an analog of
what happens in algebraic geometry, when one considers motives
instead of algebraic varieties and defines morphisms as {\em
correspondences}.

In the context of noncommutative geometry, there are well developed
analogs of correspondences and of motivic cohomology. For the
purposes of this paper we consider correspondences given by
Kasparov's $KK$-theory (\cf \cite{Kas}), although we will see later
that a more refined version is desirable, where for instance one
does not mod out by homotopy equivalence, and possibly based on a
relative version in the fibration relating algebraic and topological
$K$-theory (\cf \cite{CoKa}).

The analog of motivic and absolute cohomology is played by
modules over the cyclic category and cyclic cohomology
(\cf \cite{CoExt}, \cite{CoIHES}).

\subsection{$KK$-theory}

Similarly to what happens in the context of algebraic varieties,
$C^*$-algebras do not form an additive category. However, there is a
functor from the category of all separable $C^*$-algebras, with
$*$-homomorphisms, to an {\em additive} category $\cK\cK$, whose
objects are separable $C^*$-algebras and where, for $\bar\cA$ and
$\bar\cB$ in $Obj(\cK\cK)$, the morphisms are given by
$\Hom(\bar\cA,\bar\cB)=KK(\bar\cA,\bar\cB)$. Here
$KK(\bar\cA,\bar\cB)$ is Kasparov's bivariant $K$-theory
(\cite{Kas}, \cf also \S 8 and \S 9.22 of \cite{Bla}). Here
the notation $\bar\cA$ for $C^*$-algebras comes from the fact
that, later on in this paper, we will consider $C^*$-algebras that
are $C^*$-completions of certain natural dense subalgebras $\cA$,
though this assumption is not needed in this subsection.

In the $\cK\cK$ category morphisms are defined through
``correspondences'' in the following way. One considers the set
$M(\bar\cA,\bar\cB)$ of Kasparov modules, that is, of triples
$(E,\phi,F)$, where $E$ is a countably generated Hilbert module over
$\bar\cB$, $\phi$ is a $*$-homomorphism of $\bar\cA$ to bounded
linear operators on $E$ (\ie it gives $E$ the structure of an
$(\bar\cA,\bar\cB)$-bimodule) and $F$ is a bounded linear operator
on $E$ with the properties that the operators $[F,\phi(a)]$,
$(F^2-1)\phi(a)$, and $(F^*-F)\phi(a)$ are in $\cK(E)$ (compact
operators) for all $a\in \bar\cA$. Recall here that a Hilbert module
$E$ over $\bar\cB$ is a right $\bar\cB$-module with a positive
$\bar\cB$-valued inner product which satisfies $\langle x, yb
\rangle= \langle x, y\rangle b$ for all $x,y\in E$ and $b\in
\bar\cB$, and with respect to which $E$ is complete. Kasparov
modules are ``Morita type'' correspondences that generalize
$*$-homomorphisms of $C^*$-algebras (the latter are trivially seen
as Kasparov modules of the form $(\bar\cB,\phi,0)$).

One then considers on $M(\bar\cA,\bar\cB)$ the equivalence relation
of homotopy. Two elements are homotopy equivalent
$(E_0,\phi_0,F_0)\sim_h (E_1,\phi_1,F_1)$ if there is an element
$(E,\phi,F)$ of $M(\bar\cA,I\bar\cB)$, where $I\bar\cB=\{ f:[0,1]\to
\bar\cB\,|\, f \text{ continuous} \}$, such that
$(E\hat\otimes_{f_i}\bar\cB,f_i\circ \phi,f_i(F))$ is unitarily
equivalent to $(E_i,\phi_i,F_i)$, \ie there is a unitary in bounded
operators from $E\hat\otimes_{f_i}\bar\cB$ to $E_i$ intertwining the
morphisms $f_i\circ \phi$ and $\phi_i$ and the operators $f_i(F)$
and $F_i$. Here $f_i$ is the evaluation at the endpoints
$f_i:I\bar\cB\to \bar\cB$. By definition $KK(\bar\cA,\bar\cB)$ is
the set of homotopy equivalence classes of $M(\bar\cA,\bar\cB)$.
This is naturally an abelian group.

A slightly different formulation of $KK$-theory, which simplifies
the external tensor product $KK(\bar\cA,\bar\cB)\otimes
KK(\bar\cC,\bar\cD)\to KK(\bar\cA\otimes \bar\cC,\bar\cB\otimes
\bar\cD)$, is obtained by replacing in the data $(E,\phi,F)$ the
operator $F$ by an {\em unbounded} regular self-adjoint operator
$D$. The corresponding  $F$ is then given by $D(1+D^2)^{-1/2}$ (\cf
\cite{BaJu}).

The category $\cK\cK$ obtained this way is a universal enveloping
additive category for the category of $C^*$-algebras (\cf \cite{Bla}
\S 9.22.1). Though it is not an abelian category, it is shown in
\cite{MeyerNest} that $\cK\cK$ and the equivariant $\cK\cK_G$
admit the structure of a triangulated category.

In algebraic geometry the category of pure motives of algebraic
varieties is obtained by replacing the (non-additive) category of
smooth projective varieties by a category whose objects are triples
$(X,p,m)$ with $X$ a smooth projective algebraic variety, $p$ a
projector and $m$ an integer. The morphisms $$\Hom((X,p,m),(Y,q,n))=
q{\rm \Corr}^{n-m}(X,Y)p$$ are then given by correspondences
(combinations of algebraic cycles in the product $X\times Y$
with rational coefficients) modulo a suitable
equivalence relation. (Notice that the idempotent condition $p^2=p$
itself depends on the equivalence relation.) Depending on the choice
of the equivalence relation, one gets categories with rather
different properties. For instance, only the numerical equivalence
relation, which is the coarsest among the various relations
considered on cycles, is known to produce in general an abelian category (\cf
\cite{Jannsen}). The homological equivalence relation, for example,
is only known to produce a $\Q$-linear pseudo-abelian category, that is, an
additive category in which the ranges of projectors are included
among the objects. One of Grothendieck's standard conjectures would
imply the equivalence of numerical and homological equivalence.

In the case of (noncommutative) algebras, when we use $KK$-theory,
the analysis of the defect of surjectivity of the assembly map is
intimately related to a specific idempotent (the $\gamma$-element)
in a $KK$-theory group. Thus the analogy with motives suggests that
one should similarly enlarge the additive category $\cK\cK$ of
$C^*$-algebras with ranges of idempotents in $KK$. This amounts to
passing to the pseudo-abelian envelope of the additive $\cK\cK$
category.

\subsection{The abelian category of $\Lambda$-modules and
cyclic (co)homology}

We now interpret cyclic cohomology in a setting which is somewhat
analogous to the one relating algebraic varieties and motives. The
category of pure motives $\cM$ (with the numerical equivalence) is a
semi-simple abelian $\Q$-linear category with finite dimensional ${\rm Hom}$
groups. In the theory of motives, to a smooth projective algebraic
variety $X$ over a field $\K$ one wants to associate a {\em
motivic cohomology} $H^i_{mot}(X)$, defined as an object in
category $\cM_{\K}$,\ie
$H^i_{mot}(X)$ is a pure motive of weight $i$ (\cf \cite{Del3}).
In this way, one views pure motives as a universal
cohomology theory for algebraic varieties. Namely, the motivic
cohomology has the property that, for any given cohomology theory
$H$, with reasonable properties, there exists a realization functor
$R_H$ satisfying
\begin{equation}\label{realize}
H^n(X)=R_H\,H^n_{mot}(X).
\end{equation}
On a category of mixed motives, one would also have
{\em absolute cohomology} groups given by (\cf \cite{Del3})
\begin{equation}\label{ExtM}
 H^i_{abs}(M)={\rm Ext}^i(1,M).
\end{equation}
for $M$ a motive. Here the ${\rm Ext}^i$ are taken in a
suitable triangulated category $\cD(\K)$ whose heart is the
category of mixed motives.
For a variety $X$ the absolute cohomology would then be obtained
from the motivic cohomology via a spectral sequence
\begin{equation}\label{absmot}
 E^{pq}_2 = H^p_{abs}(H^q_{mot}(X)) \Rightarrow H^{p+q}_{abs}(X).
\end{equation}

In the context of noncommutative
geometry, we may argue that
the cyclic category and cyclic (co)homology can be thought of
as an analog of motivic and absolute cohomology.

In fact, the approach to cyclic cohomology via the cyclic
category (\cf \cite{CoExt}) provides a way to embed the nonadditive
category of algebras and algebra homomorphisms in an {\em abelian}
category of modules. The latter is the category of
$\K(\Lambda)$-modules, where $\Lambda$ is the cyclic category and
$\K(\Lambda)$ is the group ring of $\Lambda$ over a given field
$\K$. Cyclic cohomology is then obtained as an $\Ext$ functor (\cf
\cite{CoExt}).

The cyclic category $\Lambda$ has the same objects as the small
category $\Delta$ of totally ordered finite sets and increasing
maps, which plays a key role in simplicial topology. Namely,
$\Delta$ has one object $[n]$ for each integer $n$, and is generated
by faces $\delta_i: [n-1] \to [n]$ (the injection that misses $i$),
and degeneracies $\sigma_j: [n+1] \to [n] $ (the surjection which
identifies $j$ with $j+1$), with the relations
\begin{equation}\label{ad}
\delta_j  \delta_i = \delta_i  \delta_{j-1} \ \ \ \hbox{for} \ i < j
, \ \ \ \ \ \sigma_j  \sigma_i = \sigma_i  \sigma_{j+1}, \ \ \  i
\leq j
\end{equation}
\begin{equation}\label{ax}
\sigma_j  \delta_i = \left\{ \begin{array}{ll} \delta_i  \sigma_{j-1}  & i < j \\
1_n  &\text{if} \ i=j \ \text{or} \ i = j+1 \\
\delta_{i-1}  \sigma_j  & i > j+1  . \end{array} \right.
\end{equation}
To obtain the cyclic category $\Lambda$ one adds for each $n$ a new
morphism $\tau_n: [n] \to [n]$, such that
\begin{equation}\label{ae}
\begin{array}{lll} \tau_n  \delta_i = \delta_{i-1} \tau_{n-1} & 1 \leq i \leq
n , & \tau_n  \delta_0 = \delta_n \\[3mm] \tau_n  \sigma_i =
\sigma_{i-1} \tau_{n+1} &1 \leq i \leq n , &\tau_n  \sigma_0 =
\sigma_n
\tau_{n+1}^2 \\[3mm] \tau_n^{n+1} = 1_n. &&  \end{array}
\end{equation}
 Alternatively,
$\Lambda$ can be defined by means of its ``cyclic covering'',
through the category $\tilde\Lambda$. The latter has one object $(\Z
, n)$ for each $n \geq 0$ and the morphisms $f : (\Z , n) \to (\Z ,
m)$ are given by non decreasing maps $f : \Z \to \Z \ $, such that $
f(x+n) = f(x)+m \,,\quad \forall  x \in \Z$. One has $\Lambda =
\tilde\Lambda / \Z$, with respect to the obvious action of $\Z$ by
translation.

To any algebra $\cA$ over a field $\K$, one associates a module
$\cA^{\natural}$ over the category $\Lambda$ by assigning to each
$n$ the $(n+1)$ tensor power $\cA\otimes \cA\cdots \otimes \cA$. The
cyclic morphisms correspond to the cyclic permutations of the
tensors while the face and degeneracy maps correspond to the algebra
product of consecutive tensors and the insertion of the unit (\cf
\cite{CoExt}). A trace $\varphi: \cA \to \K$ gives rise to a
morphism $\varphi^\natural$ of $\cA^{\natural}$ to the
$\Lambda$-module $\K^\natural$ such that
\begin{equation}\label{taunat}
\varphi^\natural(a_0\otimes \cdots \otimes a_n)=\,\varphi(a_0 \cdots
a_n) .
\end{equation}

Using the natural isomorphism of  $\Lambda$ with the opposite small
category $\Lambda^o$ one obtains by dualization that any algebra
$\cA$ gives rise canonically to a module $\cA^\sharp$ over the small
category $\Lambda$, by assigning to each integer $n \geq 0$ the
vector space $C^n$ of $n+1$-linear forms $\varphi (x^0 , \ldots ,
x^n)$ on $\cA$, with the basic operations $\delta_i: C^{n-1} \to
C^n$ and $\sigma_i: C^{n+1} \to C^n$ given by
\begin{equation}\label{ag}
\begin{array}{rl} (\delta_i  \varphi) (x^0 , \ldots , x^n) =& \varphi (x^0 ,
\ldots , x^i x^{i+1} , \ldots , x^n), \,\, i=0,1,\ldots , n-1 \\[3mm]
(\delta_n  \varphi) (x^0 , \ldots , x^n) =& \varphi (x^n  x^0 , x^1
,
\ldots , x^{n-1})  \\[3mm] (\sigma_j  \varphi) (x^0 , \ldots , x^n)
=& \varphi (x^0 , \ldots , x^j , 1 , x^{j+1} , \ldots , x^n), \,\,
j=0,1,\ldots , n  \\[3mm] (\tau_n  \varphi) (x^0 , \ldots , x^n) =&
\varphi (x^n , x^0 , \ldots , x^{n-1}).  \end{array}
\end{equation}
These operations satisfy the relations (\ref{ad}) (\ref{ax}) and
(\ref{ae}). This shows that $\cA^\sharp$ is a $\Lambda$-module.

It is then possible (\cf \cite{CoExt}, \cite{Lo}), to interpret the
cyclic cohomology $HC^n$ as $\Ext^n$ functors, namely
\begin{equation}\label{ExtC}
HC^n(\cA) = \Ext^n (\cA^\natural, \K^\natural)= \Ext^n (\K^\sharp,
\cA^\sharp),
\end{equation}
for $\cA$ an algebra over  $\K$.
Both terms express the derived functor of the functor
which assigns to a $\Lambda$-module its space of traces.

The formula \eqref{ExtC} may be regarded as an analog
of \eqref{ExtM}, except for the fact that here the ${\rm Ext}^n$ are taken
in the abelian category of $\K(\Lambda)$-modules.

Similarly to \eqref{ExtC}, one also has an
interpretation of cyclic homology as $\Tor_n$
functors,
\begin{equation}\label{Tor}
HC_n(\cA) = \Tor_n (\K^\natural,\cA^\natural)\,.
\end{equation}

A canonical projective biresolution of $\Z^\natural$ is obtained by
considering a bicomplex of $\Lambda$-modules $(C^{n,m},d_1,d_2)$
(for $n,m>0$) where $C^{n,m}=C^m$ for all $n$ and the component
$C^k_j$ of $C^k$ is the free abelian group on the set
$\Hom(\Lambda_k,\Lambda_j)$ so that $\Hom_\Lambda (C^k,\cE)=\cE_k$.
See \cite{CoExt} for the differentials $d_1: C^{n+1,m}\to C^{n,m}$
and $d_2: C^{n,m+1} \to C^{n,m}$. This projective resolution of
$\Z^\natural$ determines a bicomplex computing the cyclic cohomology
$\Ext^n (\Z^\natural,\cE)$ of a $\Lambda$-module $\cE$ analogous to
the $(b,B)$-bicomplex $C^{n-m}(\cA)$.

All of the general properties of cyclic cohomology such as the long
exact sequence relating it to Hochschild cohomology are shared by
$\Ext$ of general $\Lambda$-modules and can be attributed to the
equality of the classifying space $B\Lambda$ of the small category
$\Lambda$ with the classifying space $BS^1$ of the compact
one-dimensional Lie group $S^1$. One has
\begin{equation}\label{af}
B\Lambda = BS^1 =\P^{\infty}(\C)
\end{equation}
This follows from the fact that $\Hom (\Lambda_n,\Lambda_m)$ is
nonempty and that ${\rm Hom}(\Lambda_0,\Lambda_0)=1$, which imply
that $B\Lambda$ is connected and simply connected, and from the
calculation (\cf \cite{CoExt}) of $\Ext^*_\Lambda
(\Z^\natural,\Z^\natural)=\Z[\sigma]$ as a polynomial ring in the
generator $\sigma$ of degree two. The cohomology ring
$H^*(B\Lambda,\Z)$ is given by $\Z[\sigma]$ where $\sigma$
corresponds to a map $f: B\Lambda \to \P^{\infty}(\C)$ which is a
homotopy equivalence.

\medskip

Finally we note that there is a natural way to associate a cyclic
morphism to a correspondence viewed as an $\cA$--$\cB$ bimodule
$\cE$, provided that the following finiteness holds

\begin{Lem} \label{cyreal} Let $\cE$ be an $\cA$--$\cB$ bimodule,
which is finite projective as a right $\cB$-module. Then there is an
 associated cyclic morphism
$$
\cE^\natural \in \Hom(\cA^\natural,\cB^\natural).
$$
\end{Lem}

\proof Since $\cE$ is finite projective it is a direct summand of
the right module $\cB^n$ for some n, so that $\End_{\cB}(\cE)\subset
M_n(\cB)=\cB\otimes M_n(\C)$. Thus one obtains $\tau^\natural$ from
the left action of $\cA$ which gives a homomorphism $\rho$ from
$\cA$ to $\End_{\cB}(\cE)$ and hence to $\cB\otimes M_n(\C)$. One
then takes  the composition of $\rho^\natural$ with the trace map
from $\cB\otimes M_n(\C)$ to $\cB$ which is also a cyclic morphism
as follows from the more general Proposition \ref{partial} below.
\endproof

\begin{Rem}
Note that the above cyclic morphism depends upon the choice of a
connection on the right $\cB$-module $\cE$, which is used in order
to write $\cE$ as a direct summand of the right module $\cB^n$ for
some n.
\end{Rem}

\medskip

All of the discussion above, regarding the cyclic module associated
to an algebra, extends to the context of locally convex topological
algebras. In that context one needs to work with topological tensor
products in the construction of $\cA^\natural$ and of continuous
multilinear forms in the construction of  $\cA^\sharp$.

Moreover, an important issue arises, since the ranges of continuous
linear maps are not necessarily closed subspaces. In order to
preserve the duality between cyclic homology and cyclic cohomology
we shall define the cokernel of a cyclic map $T: \cA^\natural \to
\cB^\natural$ as the quotient of $\cB^\natural$ by the closure of
the range of $T$. In a dual manner, the kernel of the transposed map
$T^t :  \cB^\sharp \to \cA^\sharp$ is automatically closed and is
the dual of the above.

\bigskip

\section{Artin Motives and noncommutative spaces}\label{Sartin}

In this section we show how certain categories of motives can be
embedded faithfully into categories of noncommutative spaces. We
first recall some general facts about motives.

Historically, the roots of the idea of motives introduced by
Grothendieck can be found in the $L$-functions $L(X,s)$ of algebraic
varieties, which can be written as an alternating product

\begin{equation}\label{Lprod}
L(X,s)=\prod_{i=0}^n L(H^i(X),s)^{(-1)^{i+1}}.
\end{equation}

Typically, such combined $L$-function as \eqref{Lprod} tends to have
infinitely many poles, corresponding to the zeros of the factors in
the denominator, whereas in isolating a single factor $L(H^i(X),s)$
one avoids this problem and also gains the possibility of having a
nice functional equation. It is convenient to think of such a factor
$L(H^i(X),s)$ as the $L$-function, not of a variety, but of a pure
motive of weight $i$, viewed as a direct summand of the variety $X$.

In algebraic geometry, the various constructions of categories of
pure motives of algebraic varieties are obtained by replacing the
(non additive) category $\cV_\K$ of algebraic varieties over a field
$\K$ by a (pseudo)abelian category $\cM_\K$ constructed in three
steps (\cf \cite{Man-mot}): (1) Enlarging the set of morphisms of
algebraic varieties to correspondences by algebraic cycles (modulo a
suitable equivalence relation); (2) Adding the ranges of projectors
to the collection of objects; (3) Adding the ``non-effective'' Tate
motives.

\subsection{Complex coefficients}\label{CcoeffS}

Since in the following we deal with reductions of motives by Hecke
characters which are complex valued, we need to say a few words
about coefficients. In fact, we are interested in considering
categories of motives $\cM_{\K,\E}$ over a field $\K$ with
coefficients in an extension $\E$ of $\Q$. This may sometime be a
finite extension, but we are also interested in treating the case
where $\E=\C$, \ie the category $\cM_{\K,\C}$ of motives over $\K$
with complex coefficients.

By this we mean the category obtained as follows. We let $\cM_\K$
denote, as above the (pseudo)abelian category of pure motives over
$\K$. We then define $\cM_{\K,\C}$ as the (pseudo)abelian envelope
of the additive category with $Obj(\cM_\K)\subset Obj(\cM_{\K,\C})$
and morphisms
\begin{equation}\label{Ccoeff}
\Hom_{\cM_{\K,\C}}(M,N)=\Hom_{\cM_{\K}}(M,N)\otimes \C.
\end{equation}
Notice that, having modified the morphisms to \eqref{Ccoeff} one in
general will have to add new objects to still have a (pseudo)abelian
category.

\smallskip

Notice that the notion of motives with coefficients that we use here
is compatible with the apparently different notion often adopted in
the literature, where motives with coefficients are defined for $\E$
a {\em finite} extension of $\K$, \cf \cite{Ta} and \cite{De} \S
2.1. In this setting one says that a motive $M$ over $\K$ has
coefficients in $\E$ if there is a homomorphism $\E \to \End(M)$.
This description is compatible with defining, as we did above,
$\cM_{\K,\E}$ as the smallest (pseudo)abelian category containing
the pure motives over $\K$ with morphisms
$\Hom_{\cM_{\K,\E}}(M,N)=\Hom_{\cM_{\K}}(M,N)\otimes \E$. This is
shown in \cite{De} \S 2.1: first notice that to a motive $M$ over
$\K$ one can associate an element of $\cM_{\K,\E}$ by taking
$M\otimes \E$, with the obvious $\E$-module structure. If $M$ is
endowed with a map $\E \to \End(M)$, then the corresponding object
in $\cM_{\K,\E}$ is the range of the projector onto the largest
direct summand on which the two structures of $\E$-module agree.
Thus, while for a finite extension $\E$, both notions of ``motives
with coefficients'' make sense and give rise to the same category
$\cM_{\K,\E}$, for an infinite extension like $\C$ only the
definition we gave above makes sense.

\subsection{Artin motives}

Since our main motivation comes from the ad\`ele class space of
\cite{AC}, we focus on zero-dimensional examples, like the
Bost--Connes system of \cite{BC} (whose dual system gives the
ad\`ele class space).  Thus, from the point of view of motives, we
consider mostly the simpler case of Artin motives. We return to
discuss some higher dimensional cases in Sections \ref{Shigher} and
\ref{SLfactors}.

Artin motives are motives of zero-dimensional smooth algebraic varieties
over a field $\K$
and they provide a geometric counterpart for finite dimensional
linear representations of the absolute Galois group of $\K$.

We recall some facts about the category of Artin motives (\cf
\cite{Andre}, \cite{Ta}). We consider in particular the case where
$\K$ is a number field, with an algebraic closure $\bar \K$ and with
the absolute Galois group $G=\Gal(\bar \K/\K)$. We fix an embedding
$\sigma:\K\hookrightarrow \C$ and write $X(\C)$ for $\sigma X (\C)$.

Let $\cV_\K^0$ denotes the category of zero-dimensional smooth
projective varieties over $\K$. For $X\in Obj(\cV_\K^0)$ we denote
by $X(\bar \K)$ the set of algebraic points of $X$. This is a finite
set on which $G$ acts continuously. Given $X$ and $Y$ in $\cV_\K^0$
one defines $M(X,Y)$ to be the finite dimensional $\Q$-vector space
\begin{equation}\label{morArt}
M(X,Y):=\Hom_G(\Q^{X(\bar \K)},\Q^{Y(\bar \K)})=(\Q^{X(\bar
\K)\times Y(\bar \K)})^G.
\end{equation}

Thus, a (virtual) correspondence $Z \in M(X,Y)$ is a formal linear
combination of connected components of $X\times Y$ with coefficients
in $\Q$. These are identified with $\Q$-valued $G$-invariant
functions on $X(\bar \K)\times Y(\bar \K)$ by $Z=\sum a_i Z_i
\mapsto \sum a_i \chi_{Z_i(\bar \K)}$, where $\chi$ is the
characteristic function, and the $Z_i$ are connected components of
$X\times Y$.

We denote by $\cC\cV_\K^0$ the (pseudo)abelian envelope of the
additive category with objects $Obj(\cV_\K^0)$ and morphisms as in
\eqref{morArt}. Thus, $\cC\cV_\K^0$ is the smallest (pseudo)abelian
category with $Obj(\cV_\K^0)\subset Obj(\cC\cV_\K^0)$ and morphisms
$\Hom_{\cC\cV_\K^0}(X,Y)=M(X,Y)$, for $X,Y\in Obj(\cV_\K^0)$. In this
case it is in fact an abelian category.

Moreover, there is a fiber functor
\begin{equation}\label{ArtinFiber}
\omega: X \mapsto H^0(X(\C),\Q)=\Q^{X(\bar \K)},
\end{equation}
which gives an identification of the category $\cC\cV_\K^0$ with the
category of finite dimensional $\Q$-linear representations of
$G=\Gal(\bar \K/\K)$.

\smallskip

Notice that, in the case of Artin motives, the issue of different
properties of the category of motives depending on the different
possible choices of the equivalence relation on cycles (numerical,
homological, rational) does not arise. In fact, in dimension zero
there is no equivalence relation needed on the cycles in the
product, \cf \eg \cite{Andre}. Nonetheless, one sees that passing to
the category $\cC\cV_\K^0$ requires adding new
objects in order to have a (pseudo)abelian category. One can see
this in a very simple example, for $K=\Q$. Consider the field
$\L=\Q(\sqrt 2)$. Then consider the one dimensional non-trivial
representation of the Galois group $G$ that factors through the
character of order two of $\Gal(\L/\Q)$. This representation does
not correspond to an object in $Obj(\cV_\Q^0)$ but can be obtained
as the range of the projector $p=(1-\sigma)/2$, where $\sigma$ is
the generator of $\Gal(\L/\Q)$. Namely it is a new object in
$Obj(\cC\cV_\Q^0)$.

\subsection{Endomotives}

We now show that there is a very simple way to give a faithful
embedding of the category $\cC\cV_{\K}^0$ of Artin motives in an
additive category of noncommutative spaces. The construction extends
to arbitrary coefficients $\cC\cV_{\K,\E}^0$ (in the sense of \S
\ref{CcoeffS}) and in particular to $\cC\cV_{\K,\C}^0$. The objects
of the new category are obtained from projective systems of Artin
motives with actions of semigroups of endomorphisms. For this reason
we will call them {\em endomotives}. We construct a category
$\cE\cV^0_{\K,\E}$ of algebraic endomotives, which makes sense over
any field $\K$. For $\K$ a number field, this embeds in a category
$C^*\cV^0_\K$ of analytic endomotives, which is defined in the
context of $C^*$-algebras.

\smallskip

One can describe the category $\cV_\K^0$ as the
category of reduced finite dimensional commutative algebras over
$\K$. The correspondences for Artin motives given by \eqref{morArt}
can also be described in terms of bimodules (\cf \eqref{intpro}
and \eqref{FUkk} below). We consider a specific class of
noncommutative algebras to which all these notions extend
naturally in the zero-dimensional case.

\smallskip

The noncommutative algebras we deal with are of the form
\begin{equation}\label{algdef}
\cA_\K=\,A\rtimes S ,
\end{equation}
where $A$ is an inductive limit of reduced finite dimensional
commutative algebras over the field $\K$ and $S$ is a unital abelian
semigroup of algebra endomorphisms $\rho:  A \to A$. The algebra $A$
is unital and, for $\rho \in S$, the image $e=\rho(1)\in A$ is an
idempotent. We assume that $\rho$ is an isomorphism of $A$ with the
compressed algebra $eAe$.

\smallskip

The crossed product algebra $\cA_\K=\,A\rtimes S $ is obtained by
adjoining to $A$ new generators $U_\rho$ and $U^*_\rho$,  for $\rho
\in S$, with algebraic rules given by
\begin{equation}\label{algrulesS}
\begin{array}{lll}
U^*_\rho U_\rho =1, & U_\rho U^*_\rho=\rho(1), & \forall
\rho \in S \\[2mm]
U_{\rho_1\,\rho_2}=\,U_{\rho_1}\,U_{\rho_2}, &
U^*_{\rho_2\,\rho_1}=\,U^*_{\rho_1}\,U^*_{\rho_2}, & \forall
\rho_j \in S \\[2mm]
U_\rho\,x=\,\rho(x)\,U_\rho, &  x\,U^*_\rho=\,U^*_\rho\,\rho(x), &
\forall \rho \in S, \, \forall x\in A .
\end{array}
\end{equation}

\medskip
Since $S$ is abelian these rules suffice to show that $A\rtimes S $
is the linear span of the monomials $U^*_{\rho_1}\,a\,U_{\rho_2}$
for $a\in A$ and $\rho_j\in S$. In fact one gets
\begin{equation}\label{monom}
U^*_{\rho_1}\,a\,U_{\rho_2}\,U^*_{\rho_3}\,b\,U_{\rho_4}=\,U^*_{\rho_1
\rho_3}\,\rho_3(a)\,\rho_2\rho_3(1)\,\rho_2(b)\,U_{\rho_2
\rho_4}\qqq \rho_j\in S\,,a,b\in A\,.
\end{equation}

\medskip

By hypothesis, $A$ is an inductive limit of reduced finite
dimensional commutative algebras $A_\alpha$ over $\K$. Thus, the
construction of the algebraic semigroup crossed product $\cA_\K$
given above is determined by assigning the following data: a
projective system $\{ X_\alpha \}_{\alpha\in I}$ of varieties in
$\cV^0_\K$, with morphisms $\xi_{\alpha,\beta}: X_\beta\to
X_\alpha$, and a countable indexing set $I$.  The graphs
$\Gamma(\xi_{\alpha,\beta})$ of these morphisms are $G$-invariant
subsets of $X_\beta(\bar \K)\times X_\alpha(\bar \K)$. We denote by
$X$ the projective limit, that is, the pro-variety
\begin{equation}\label{projlimX}
X =\varprojlim_\alpha X_\alpha ,
\end{equation}
with maps $\xi_\alpha: X \to X_\alpha$.

\smallskip

The endomorphisms $\rho$ give isomorphisms $$ \tilde \rho:  X\to
X^e$$ of $X$ with the subvariety $X^e$ associated to the idempotent
$e=\rho(1)$, \ie corresponding to the compressed algebra $eAe$.

\smallskip

The noncommutative space described by the algebra $\cA_\K$ of
\eqref{algdef} is the quotient of $X(\bar \K)$ by the  action of
$S$, \ie of the $\tilde \rho$.

\smallskip

By construction, the Galois group $G$ acts on $X(\bar \K)$ by
composition. Thus, if we view the elements of $X(\bar \K)$ as
characters, \ie as $\K$-algebra homomorphisms $\chi: A\to \bar \K$,
we can write the Galois action of $G$ as
\begin{equation}\label{charaction}
\alpha(\chi)=\,\alpha \circ \chi : A\to \bar \K, \ \ \forall \alpha
\in G=\,\Aut_\K(\bar \K),\ \ \forall \chi \in X(\bar \K)\,.
\end{equation}

This action commutes with the maps $ \tilde \rho$. In fact, given a
character  $\chi \in X(\bar \K)$ of $A$, one has
$\chi(\rho(1))\in\{0,1\}$ (with $\chi(\rho(1))=1 \Leftrightarrow
\chi \in X^e$) and $(\alpha \circ \chi)\circ \rho=\,\alpha \circ
(\chi\circ \rho)$.

\smallskip

Thus,  the whole construction is $G$-equivariant. Notice that this
does not mean that $G$ acts by automorphisms on $\cA_\K$, but it
does act on $X(\bar \K)$, and on the noncommutative quotient $X(\bar
\K)/S$.

\smallskip

The purely algebraic construction of the crossed product algebra
$\cA_\K=A\rtimes S$ and of the action of $G$ and $S$ on $X(\bar \K)$
given above makes sense also when $\K$ has positive characteristic.
One can extend the basic notions of Artin motives in the above
generality. For the resulting crossed product algebras
$\cA_\K=A\rtimes S$ one can define the set of correspondences
$M(\cA_\K,\cB_\K)$ using $\cA_\K$--$\cB_\K$-bimodules (which are
finite and projective as right modules). One can then use Lemma
\ref{cyreal} to obtain a realization in the abelian category of
$\K(\Lambda)$-modules. This may be useful in order to consider
systems like the generalization of the BC system of \cite{BC} to
rank one Drinfeld modules recently studied by Benoit Jacob
\cite{Jacob}.

\medskip

In the following, we concentrate mostly on the case of
characteristic zero, where $\K$ is a number field. We fix as above
an embedding $\sigma:\K\hookrightarrow \C$ and we take  $\bar \K$ to
be the algebraic closure of $\sigma(\K)\subset \C$ in $\C$. We then
let
\begin{equation}\label{complexify}
A_\C= A \otimes_\K \C= \varinjlim_\alpha A_\alpha\otimes_\K \C ,
 \    \    \   \cA_\C= \cA_\K\otimes_\K \C= A_\C \rtimes S .
\end{equation}
We define an embedding of algebras $A_\C\subset C(X)$ by
\begin{equation}\label{embed}
a \in A \to \hat a, \ \ \ \hat a(\chi)=\,\chi(a)\ \ \ \forall \chi
\in X.
\end{equation}
Notice that $A$ is by construction a subalgebra $A\subset C(X)$ but
it is not in general an involutive subalgebra for the canonical
involution of $C(X)$. It is true, however, that the $\C$-linear span
$A_\C\subset C(X)$ is an involutive subalgebra of $C(X)$ since it is
the algebra of locally constant functions on $X$. The
$C^*$-completion $R=C(X)$ of $A_\C$ is an abelian AF $C^*$-algebra.
We then let $\bar\cA=C(X)\rtimes S$ be the $C^*$-crossed product of
the $C^*$-algebra $C(X)$ by the semi-group action of $S$. It is the
$C^*$-completion of the algebraic crossed product $A_\C\rtimes S$
and it is defined by the algebraic relations \eqref{algrulesS} with
the involution which is apparent in the formulae (\cf  \cite{lr1}
\cite{l1}  and the references there for the general theory of
crossed products by semi-groups).

\smallskip

It is important, in order to work with cyclic (co)homology as we do
in Section \ref{Frobenius} below, to be able to restrict from
$C^*$-algebras $C(X)\rtimes S$ to a canonical dense subalgebra
\begin{equation}\label{densesub}
\cA= C^\infty(X)\rtimes_{alg} S \subset \bar\cA=C(X)\rtimes S
\end{equation}
where $\cA=C^\infty(X)\subset C(X)$ is the algebra of locally
constant functions, and the crossed product $C^\infty(X)\rtimes S$
is the algebraic one. It is to this category of (smooth) algebras,
rather than to that of $C^*$-algebras, that cyclic homology applies.

\medskip

\begin{Pro}\label{fab} 1) The action \eqref{charaction} of $G$ on
$X(\bar \K)$ defines a canonical action of $G$ by automorphisms of
the $C^*$-algebra $\bar\cA=C(X)\rtimes S$, preserving globally
$C(X)$ and such that, for any pure state $\varphi$ of $C(X)$,
\begin{equation}\label{fabulous0}
\alpha \,\varphi(a)=\, \varphi (\alpha^{-1}(a)) \qqq a\in A \,,\quad
\alpha \in G.
\end{equation}

2) When the Artin motives $A_\alpha$ are abelian and normal, the
subalgebras $A\subset C(X)$ and $\cA_\K\subset \bar\cA=C(X)\rtimes
S$ are globally invariant under the action of $G$ and the states
$\varphi$ of $R\rtimes S$ induced by pure states of $R$ fulfill
\begin{equation}\label{fabulous}
\alpha \,\varphi(a)=\, \varphi (\theta(\alpha)(a)) \qqq a\in \cA_\K
\,,\quad \theta(\alpha)=\alpha^{-1} \qqq \alpha \in G^{ab}=\,G/[G,G]
\end{equation}
\end{Pro}

\proof 1) The action of the Galois group $G$ on $C(X)$ is given,
using the notations of \eqref{charaction}, by
\begin{equation}\label{funaction}
\alpha(h)(\chi)=\,h(\alpha^{-1}(\chi)) \qqq \chi \in X(\bar
\K)\,,\quad h\in C(X)\,.
\end{equation}
For $a\in A \subset C(X)$ and $\varphi$ given by the evaluation at
$\chi \in X(\bar \K)$ one gets
$$
\alpha \,\varphi(a)=\, (\alpha \circ \chi)(a)= \,\alpha(\chi)(a)=\,
\hat a(\alpha(\chi))=\, \varphi (\alpha^{-1}(\hat a))
$$
using \eqref{funaction}. Since the action of $G$ on $X$ commutes
with the action of $S$, it extends to the crossed product
$C(X)\rtimes S$.

2) Let us show that $A\subset C(X)$ is globally invariant under the
action of $\alpha \in G$.  By construction $A=\cup\,A_\beta$ where
the finite dimensional $\K$-algebras $A_\beta$ are of the form
$$ A_\beta=\,\prod \L_i$$ where each $\L_i$ is a normal abelian
extension of $\K$. In particular the automorphism of $\L_i$
 defined by $$ \tilde \alpha=\, \chi^{-1}\circ \alpha \circ
\chi $$ is independent of the  choice of the embedding $\chi\,:
\,L_i\to \bar\K$. These assemble into an automorphism $\tilde
\alpha$ of $A_\beta$ such that for any $\chi\in X$ one has
$\chi\circ \tilde \alpha=\alpha \circ \chi $. One has by
construction for $a\in A$ with $\hat a $ given by \eqref{embed}, and
$b=\, \tilde \alpha(a)$, the equality
$$
\hat b(\chi)=\, \chi(b)=\,\chi \tilde \alpha(a)=\, \alpha \chi(a)=\,
\alpha^{-1}(\hat a)(\chi)
$$
so that $b=\,\alpha^{-1}(\hat a)$ and one gets the required global
invariance for $A$. The invariance of $\cA_\K$ follows. Finally note
that $\theta(\alpha)=\alpha^{-1}$ is a group homomorphism for the
abelianized groups.
\endproof

\smallskip
The intertwining property \eqref{fabulous} is (part of) the {\em
fabulous state} property of \cite{CM}, \cite{CMR}, \cite{CMR2}. The
subalgebra $\cA_\K$ then qualifies as a rational subalgebra of
$\bar\cA=C(X)\rtimes S$ in the normal abelian case. In the normal
non-abelian case the global invariance of $A$ no longer holds. One
can define an action of $G$ on $A$ but only in a  non-canonical
manner which depends on the choices of base points. Moreover in the
non-abelian case one should not confuse the contravariant
intertwining  stated in \eqref{fabulous0} with the covariant
intertwining of \cite{CM}.

\medskip

\begin{Exa}\label{BCex}{\rm
The prototype example of the data described above is the
Bost--Connes system. In this case we work over $\K=\Q$, and we
consider the projective system of the $X_n=\Sp(A_n)$, where
$A_n=\Q[\Z/n\Z]$ is the group ring of the abelian group $\Z/n\Z$.
The inductive limit is the group ring $A=\Q[\Q/\Z]$ of $\Q/\Z$. The
endomorphism $\rho_n$ associated to an element $n\in S$ of the
(multiplicative) semigroup $S=\N=\Z_{>0}$ is given on the canonical
basis $e_r \in \Q[\Q/\Z]$, $r\in \Q/\Z$, by
\begin{equation} \label{rhoBC}
\rho_n(e_r)=\,\frac{1}{n}\,\sum_{ns=r}\,e_s
\end{equation}
The direct limit $C^*$-algebra is $R=C^*(\Q/\Z)$ and the crossed
product $\bar\cA=R\rtimes S$ is the $C^*$-algebra of the BC-system.
One is in the normal abelian situation, so that proposition
\ref{fab} applies. The action of the Galois group $G$ of $\bar \Q$
over $\Q$ factorizes through the cyclotomic action and coincides
with the symmetry group of the BC-system.  In fact, the action on
$X_n=\Sp(A_n)$ is obtained by composing a character $\chi: A_n \to
\bar\Q$ with an element $g$ of the Galois group $G$. Since $\chi$ is
determined by the $n$-th root of unity $\chi(e_{1/n})$, we just
obtain the cyclotomic action. The subalgebra $\cA_\Q\subset \bar\cA$
coincides with the rational subalgebra of \cite{BC}. Finally the
fabulous state property of extreme KMS$_\infty$ states follows from
Proposition \ref{fab} and the fact (\cf \cite{BC}) that such states
are induced from pure states of $R$}.
\end{Exa}

\begin{Rem}\label{OJex}{\rm Fourier transform gives an algebra
isomorphism $$ C^*(\Q/\Z)\sim C(\hat \Z)$$ which is compatible with
the projective system $\hat \Z=\varprojlim_n \Z/n\Z $, while the
endomorphisms $\rho_n$ are given by division by $n$. It is thus
natural to consider the example given by $Y_n=\Sp(B_n)$ where $B_n$
is simply the algebra of functions $C(\Z/n\Z,\Q)$, with  the
corresponding action of $S=\N$. It is crucial to note that the
corresponding system is {\em not isomorphic} to the one described
above in Example \ref{BCex}. Indeed what happens is that the
corresponding action of the Galois group $G$ is now trivial. More
generally, one could construct examples of that kind by taking $S$
to be the ring of integers $\O$ of a number field and consider its
action on the quotients $\O/J$, by ideals $J$. The action is
transitive at each finite level, though the action of $\O$ on the
profinite completion $\hat\O$ is not. The resulting $C^*$-algebra is
$C(\hat\O)\rtimes \O$. The basic defect of these examples is that
they do not yield the correct Galois action.}
\end{Rem}

\begin{Exa}\label{CMR}
{\rm As a corollary of the results of \cite{CMR} one can show that
given an imaginary quadratic extension $\K$ of $\Q$ there is a
canonically associated noncommutative space of the form
\eqref{algdef}, with the correct non-trivial Galois action of the
absolute Galois group of $\K$. The normal abelian hypothesis of
Proposition \ref{fab} is fulfilled. The rational subalgebra obtained
from Proposition \ref{fab} is the same as the rational subalgebra of
\cite{CMR} and the structure of noncommutative Artin motive is
uniquely dictated by the structure of inductive limit of finite
dimensional $\K$-algebras of the rational subalgebra of \cite{CMR}.
The construction of the system (\cf \cite{CMR}) is much more
elaborate than in Remark \ref{OJex}. }
\end{Exa}

\medskip

We now describe a large class of examples of noncommutative spaces
\eqref{algdef}, obtained as semigroup actions on projective systems
of Artin motives, that arise from self-maps of algebraic varieties.

We consider a pointed algebraic variety $(Y,y_0)$ over $\K$ and a
countable unital abelian semigroup $S$ of finite endomorphisms of
$(Y,y_0)$, unramified over $y_0\in Y$. Finite endomorphisms have
a well defined degree. For $s\in S$ we let $\deg s$
denote the degree of $s$. We construct a projective
system of Artin motives $X_s$ over $\K$ from these data as follows.
For $s\in S$, we set
\begin{equation}\label{Xs}
X_s=\{ y\in Y:\, s(y)=y_0 \}.
\end{equation}
For a pair $s,s'\in S$, with $s'=sr$, the map $\xi_{s,s'}: X_{sr}\to
X_s$ is given by
\begin{equation}\label{Xsr}
X_{sr} \ni y \mapsto r(y)\in X_s.
\end{equation}
This defines a projective system indexed by the semigroup $S$ itself
with partial order given by divisibility. We let $X=\varprojlim_s
X_s$.

Since $s(y_0)=y_0$, the base point $y_0$ defines a component $Z_s$
of $X_s$ for all $s\in S$. Let $\xi_{s,s'}^{-1}(Z_s)$ be the inverse
image of $Z_s$ in $X_{s'}$. It is a union of components of $X_{s'}$.
This defines a projection $e_s$ onto an open and closed subset
$X^{e_s}$ of the projective limit $X$.

\begin{Pro}\label{propXs}
The semigroup $S$ acts on the projective limit $X$ by partial
isomorphisms $\beta_s: X \to X^{e_s}$ defined by the property that
\begin{equation}\label{endoS}
\xi_{su}(\beta_s(x))=\xi_u(x), \, \forall u\in S, \forall x\in X.
\end{equation}
\end{Pro}

\proof The map $\beta_s$ is well defined as a map $X\to X$ by
\eqref{endoS}, since the set $\{ su:\, u\in S \}$ is cofinal and
$\xi_u(x)\in X_{su}$. In fact, we have $su \xi_u(x)=s(y_0)=y_0$. The
image of $\beta_s$ is in $X^{e_s}$, since by \eqref{endoS} we have
$\xi_s(\beta_s(x))=\xi_1(x)=y_0$. We show that $\beta_s$ is an
isomorphism of $X$ with $X^{e_s}$. The inverse map is given by
\begin{equation}\label{invrhos}
\xi_u(\beta_s^{-1}(x))=\xi_{su}(x), \, \forall x\in X^{e_s}, \forall
u\in S.
\end{equation}
In fact, for $x\in X^{e_s}$, we have $\xi_{su}(x)\in X_u$.
\endproof

The corresponding algebra morphisms $\rho_s$ are then given by
\begin{equation}\label{endoS1}
\rho_s(f)(x)=f(\beta_s^{-1}(x)),\, \forall x\in X^{e_s}, \ \
\rho_s(f)(x)=0 ,\, \forall x\notin X^{e_s}\,.
\end{equation}

In this class of examples one has an ``equidistribution'' property,
by which the uniform normalized counting measures $\mu_s$ on $X_s$
are compatible with the projective system and define a probability
measure on the limit $X$. Namely, one has
\begin{equation}\label{measlim}
\xi_{s',s} \mu_s = \mu_{s'},\, \forall s,s'\in S.
\end{equation}
This follows from the fact that the number of preimages of a point
under $s\in S$ is equal to $\deg s$.

\medskip

\begin{Exa}\label{BCrevised} {\rm Let $Y$ be the affine group scheme
$\bG_m$ (the multiplicative group). Let $S$ be the semigroup of
non-zero endomorphisms of $\bG_m$. These correspond to maps of the
form $u\mapsto u^n$ for some non-zero $n\in \Z$. In fact one can
restrict to $\N\subset \Z\smallsetminus \{ 0 \}$.}
\end{Exa}

\medskip

\begin{Pro}
The construction of Example \ref{BCrevised}, applied to the pointed
algebraic variety $(\bG_m(\Q),1)$, gives the BC system.
\end{Pro}

\proof We restrict to the semi-group $S=\N$ and determine $X_n$ for
$n\in \N$. One has by construction $X_n =\Spec (A_n)$ where $A_n$ is
the quotient of the algebra $\Q[u(n),u(n)^{-1}]$ of Laurent
polynomials by the relation $ u(n)^n=\,1 $. For $n|m$ the map
$\xi_{m,n}$ from $X_m$ to $X_n$ is given by the algebra homomorphism
that sends the generators $u(n)^{\pm 1}\in A_n$ to $u(m)^{\pm a}\in
A_m$ with $a=\,m/n$. Thus, one obtains an isomorphism of the
inductive limit algebra $A$ with the group ring $\Q[\Q/\Z]$ of the
form
\begin{equation}\label{BClim}
\iota\;:\;\varinjlim \,A_n \to \,\Q[\Q/\Z]\,,\quad
\iota(u(n))=\,e_{\frac{1}{n}} .
\end{equation}
Let us check that the partial isomorphisms $\rho_n$ \eqref{rhoBC} of
the group ring correspond to those given by \eqref{endoS} under the
isomorphism $\iota$. We identify the projective limit $X$ with the
space of characters of $\Q[\Q/\Z]$ with values in $\bar \Q$. The
projection $\xi_m(x)$ is simply given by the restriction of $x\in X$
to the subalgebra $A_m$. Let us compute the projections of the
composition of $x\in X$ with the endomorphism $\rho_n$
\eqref{rhoBC}. One has
$$
x(\rho_n(e_r))=\,\frac{1}{n}\,\sum_{ns=r}\,x(e_s)\,.
$$
 This is non-zero iff the restriction $x|_{A_n}$ is
the trivial character, \ie iff $\xi_n(x)=1$. Moreover, in that case
one has
$$
x(\rho_n(e_r))=\,x(e_s)\qqq s \,,\quad ns = r\,,
$$
and in particular \begin{equation}\label{actrho}
x(\rho_n(e_{\frac{1}{k}}))=\,x(e_{\frac{1}{nk}}).
\end{equation}

For $k|m$ the inclusion $X_k\subset X_m$ is given at the algebra
level  by the homomorphism
$$
j_{k,m}\;:\; A_m\to A_k\,,\quad j_{k,m}(u(m))=\, u(k) .
$$
Thus, one can rewrite \eqref{actrho} as
\begin{equation}\label{actrho1}
x\circ\rho_n\circ j_{k,nk}=\,x|_{A_{nk}},
\end{equation}
that is, one has
$$
 \xi_{nk}(x)=\, \xi_k(x\circ\rho_n),
$$
which, using \eqref{invrhos}, gives the desired equality of the
$\rho$'s of \eqref{rhoBC} and \eqref{endoS1}.
\endproof

\medskip

The construction of Example \ref{BCrevised} continues to make sense
for $\bG_m(\K)$ for other fields, \eg for a field of positive
characteristic. In this case one obtains new systems, different from
BC.

\begin{Exa}\label{ellcurveex} {\rm Let $Y$ be an elliptic curve
defined over $\K$. Let $S$ be the semigroup of non-zero
endomorphisms of $Y$. This gives rise to an example in the general
class described above. When the elliptic curve has complex
multiplication, this gives rise to a system which, in the case of a
maximal order, agrees with the one constructed in \cite{CMR}. In the
case without complex multiplication, this provides an example of a
system where the Galois action does not factor through an abelian
quotient. }
\end{Exa}

\medskip

In general, we say that a system $(X_\alpha,S)$ with $X_\alpha$ a
projective system of Artin motives and $S$ a semigroup of
endomorphisms of the limit $X$ is {\em uniform} if the normalized
counting measures $\mu_\alpha$ on $X_\alpha$ satisfy
\begin{equation}\label{measX}
\xi_{\beta,\alpha}\mu_\alpha=\mu_\beta
\end{equation}
as above.

\medskip

In order to define morphisms in the category of algebraic
endomotives, it is best to encode the datum $(X_\alpha,S)$ by the
algebraic groupoid obtained as the crossed product $\cG=X\rtimes S$.
We let $\tilde S$ be the (Grothendieck) group of the abelian
semigroup $S$ and may assume, using the injectivity of the partial
action of $S$,  that $S$ embedds in $\tilde S$. Then the action of
$S$ on $X$ extends to a partial action of $\tilde S$. For
$s=\rho_1/\rho_2\in \tilde S$ the projections
\begin{equation}\label{EsFs}
E(s)=\rho_1^{-1}(\rho_2(1)\rho_1(1)) \ \ \ \text{ and } \ \ \
F(s)=\rho_2^{-1}(\rho_2(1)\rho_1(1))
\end{equation}
only depend on $s$ and the map $(\rho_2)^{-1}\rho_1$ is an
isomorphism  of the reduced algebras
\begin{equation}\label{redalgebrasiso}
 s=\rho_2^{-1}\rho_1\, :\; A_{E(s)}\to A_{F(s)}\,.
\end{equation}
One checks that
\begin{equation}\label{redalgebrasiso0}
E(s^{-1})=F(s)=s(E(s))\,, \ \ F(s\,s')\geq \,F(s)s(F(s')) \,,\  \
\forall s,s'\in \tilde S\,.
\end{equation}
The algebraic groupoid $\cG$ is the disjoint union
\begin{equation}\label{groupoidstack}
\cG=\sqcup_{s\in \tilde S}\,X^{F(s)}
\end{equation}
which corresponds to the algebra $\oplus_{s\in \tilde S} A_{F(s)}$,
which is the commutative algebra direct sum of the reduced algebras
$A_{F(s)}$. The range and source  maps are given by the natural
projection from $\cG$ to $X$ and by its composition with the
antipode ${\bf S}$ which is given, at the algebra level, by
\begin{equation}\label{groupoidstack1}
{\bf S}(a)_s= s(a_{s^{-1}}) \,,\  \   \forall s\in \tilde S \,.
\end{equation}
The composition in the groupoid corresponds to the product of
monomials
\begin{equation}\label{groupoidstack2}
a \,U_s\,b\,U_t=\,a\,s(b)\,U_{st}\,.
\end{equation}

Given algebraic endomotives $(X_\alpha,S)$ and $(X'_{\alpha'},S')$,
with associated groupoids $\cG=\cG(X_\alpha,S)$,
$\cG'=\cG(X'_\alpha,S')$ a geometric correspondence is given by a
$\cG-\cG'$-space $Z$ where the right action of $\cG'$ fulfills a
suitable {\em etale} condition which we now describe. Given a space
such as $\cG$, \ie a disjoint union $Z=\Spec\, C$ of zero
dimensional pro-varieties over $\K$, a  {\em right action} of $\cG$
on $Z$ is specified by a map $g: Z\to X$ and a collection of partial
isomorphisms
\begin{equation}\label{groupoidact11}
z \in g^{-1}(F(s)) \mapsto z\cdot s\in  g^{-1}(E(s))
\end{equation}
fulfilling the obvious rules for a partial action of the abelian
group $\tilde S$. More precisely one requires that
\begin{equation}\label{groupoidact12}
g(z \cdot s) = g(z)\cdot s \,, \  \ z\cdot (s s')=(z\cdot s)\cdot s'
\ \text{on} \ g^{-1}(F(s)\cap s(F(s')))
\end{equation}
where we denote   by $x \mapsto x\cdot s $ the given partial action
of $\tilde S$ on $X$.

One checks that such an action gives to the $\K$-linear space $C$
(for $Z=\Spec\, C$) a structure of right module over $\cA_\K$.

\begin{Def}\label{algendocorr0}
We say that the action of $\cG$ on $Z$ is {\em etale} if the
corresponding module $C$ is finite and projective over $\cA_\K$.
\end{Def}
 Given algebraic endomotives
$(X_\alpha,S)$ and $(X'_{\alpha'},S')$, an etale correspondence is a
$\cG(X_\alpha,S)-\cG(X'_{\alpha'},S')$ space $Z$ such that the right
action of $\cG(X'_{\alpha'},S')$ is etale. We  define the
$\Q$-linear space of (virtual) correspondences
$\Corr((X_\alpha,S),(X'_{\alpha'},S'))$ as formal linear
combinations $U=\sum_i a_i Z_i$ of etale correspondences $Z_i$
modulo isomorphism
 and the equivalence  $Z\sqcup Z'\sim Z + Z'$ where $Z\sqcup Z'$ is the
disjoint union. The composition of correspondences is given by the
fiber product over a groupoid. Namely, for algebraic endomotives
$(X_\alpha,S)$,
  $(X'_{\alpha'},S')$, $(X''_{\alpha''},S'')$ and correspondences $Z$ and
  $W$, their composition is given by
\begin{equation}\label{corrfiberprodalg}
Z \circ W = Z \times_{\cG'} W,
\end{equation}
which is the fiber product over the groupoid
$\cG'=\cG(X'_{\alpha'},S')$. We can then form a category of
algebraic endomotives in the following way.
\begin{Def}\label{ExtArtinAlg}
The category $\cE\cV^0_{\K,\E}$ of {\em algebraic endomotives} with
coefficients in $\E$ is the (pseudo)abelian category generated by
the following objects and morphisms. The objects are uniform systems
$(X_\alpha,S)$ as above, with $X_\alpha$ Artin motives over $\K$.
For $M=(X_\alpha,S)$ and $M'=(X'_{\alpha'},S')$ we set
\begin{equation}\label{homM}
\Hom_{\cE\cV^0_{\K,\E}}(M,M')=
\Corr((X_\alpha,S),(X'_{\alpha'},S'))\otimes \E.
\end{equation}
\end{Def}

\medskip

When taking points over $\bar\K$, the algebraic endomotives
described above yield $0$-dimensional singular quotient spaces
$X(\bar\K)/S$, which can be described through locally compact etale
groupoids $\cG(\bar\K)$ (and the associated crossed product
$C^*$-algebras $C(X(\bar\K))\rtimes S$). This gives rise to a
corresponding category of analytic endomotives given by locally
compact etale groupoids  and geometric correspondences. This
category provides the data for the construction in Section
\ref{Frobenius} below.

\medskip

Let $X$ be a totally disconnected compact space, $S$ an abelian
semigroup of homeomorphisms $X^s\to X$, $x\mapsto x\cdot s$, with
$X^s$ a closed and open subset of $X$, and $\mu$ a probability
measure on $X$ with the property that the Radon--Nikodym derivatives
\begin{equation}\label{RNder}
\frac{ds^*\mu}{d\mu}
\end{equation}
are locally constant functions on $X$. We let $\cG=X \rtimes S$ be
the corresponding etale locally compact groupoid, constructed in the
same way as in the above algebraic case. The crossed product
$C^*$-algebra $C(X(\bar\K))\rtimes S$ coincides with the
$C^*$-algebra $C^*(\cG)$ of the groupoid $G$.

\smallskip

The notion of right (or left) action of $\cG$ on a  totally
disconnected locally compact spaces $\cZ$ is defined as in the
algebraic case by \eqref{groupoidact11} and \eqref{groupoidact12}. A
right action of $\cG$ on $\cZ$ gives on the space $C_c(\cZ)$ of
continuous functions with compact support on $\cZ$ a structure of
right module over $C_c(\cG)$ as in the algebraic case.
 When the fibers of the map $g$ are discrete
(countable) subsets of $\cZ$ one can define on $C_c(\cZ)$ an inner
product with values in $C_c(\cG)$ by
\begin{equation}\label{groupoidact3loc}
\langle \xi,\eta\rangle(x,s)=\,\sum_{z\in g^{-1}(x)}\,\bar
\xi(z)\,\eta(z\circ s)
\end{equation}
We define the notion of {\em etale} action by
\begin{Def}\label{analyetaleright}
A right action of $\cG$ on $\cZ$ is {\em etale} if and only if the
fibers of the map $g$ are discrete and the  identity is a compact
operator in the right $C^*$-module $\cE_\cZ$ over $C^*(\cG)$ given
by \eqref{groupoidact3loc}.
\end{Def}
An etale correspondence is a $\cG(X_\alpha,S)-\cG(X'_{\alpha'},S')$
space $\cZ$ such that the right action of $\cG(X'_{\alpha'},S')$ is
etale. We consider the $\Q$-vector space
$$ \Corr((X,S,\mu),(X',S',\mu')) $$
of linear combinations of etale   correspondences $\cZ$ modulo the
equivalence relation $\cZ\cup \cZ'=\cZ+\cZ'$ for disjoint unions.

For $M=(X,S,\mu)$, $M'=(X',S',\mu')$, and $M''=(X'',S'',\mu'')$, the
composition of correspondences
$$ \Corr(M,M')\times \Corr(M', M'') \to \Corr(M,M'') $$ is given as above
by the fiber product over $\cG'$. In fact, a correspondence, in the
sense above, gives rise to a bimodule $\cM_\cZ$ over the algebras
$C(X)\rtimes S$ and $C(X')\rtimes S'$ and the composition of
correspondences translates into the tensor product of bimodules.

\begin{Def}\label{ExtArtCalg}
The category $C^*\cV^0_\E$ of {\em analytic endomotives} is the
(pseudo)abelian category generated by objects of the form
$(X,S,\mu)$ with the properties listed above and morphisms given as
follows. For $M=(X,S,\mu)$ and $M'=(X',S',\mu')$ we set
\begin{equation}\label{CorrMCalg}
\Hom_{C^*\cV^0_\E}(M,M')= \Corr(M,M')\otimes \E.
\end{equation}
\end{Def}

\medskip
It would also be possible to define a category where morphisms are
given by $KK$-classes $KK(C^*(\cG), C^*(\cG'))\otimes \E$. The
definition we gave above for the category $C^*\cV^0_\E$ is more
refined. In particular, we do not divide by homotopy equivalence.
The associated $KK$-class $k(\cZ)$ to a correspondence $\cZ$  is
simply given by taking the equivalence class of the triple
$(E,\phi,F)$ with $(E,\phi)$ given by the bimodule $\cM_\cZ$ with
the trivial grading $\gamma=1$ and the zero endomorphism $F=0$.
Since the correspondence is etale any endomorphism of the right
$C^*$-module $\cM_\cZ$ on $C(\cX')\rtimes S'$ is  compact. Thus in
particular $F^2-1=-1$  is compact. This is exactly the reason why
one had to require that the correspondence is etale. There is a
functor $k:C^*\cV^0_\E\to \cK\cK\otimes \E$ for the category
$\cK\cK$ described in Section \ref{Scyclic} above, which extends by
$\E$-linearity the map $\cZ\mapsto k(\cZ)$. The version of
Definition \ref{ExtArtCalg} is preferable to the $KK$-formulation,
since we want morphisms that will act on the cohomology theory that
we introduce in Section \ref{Frobenius} below.

\medskip

\begin{The}\label{NCArtinThm}
The categories introduced above are related as follows.
\begin{enumerate}
\item The map $\cG\mapsto \cG(\bar\K)$ determines a
tensor functor
\begin{equation}\label{Ffunctoralgan}
\cF :\cE\cV^0_{\K,\E} \to C^*\cV^0_\E
\end{equation}
from algebraic to analytic endomotives.
\item The Galois group $G=\Gal(\overline{\sigma(\K)}/\sigma(\K))$
acts by natural transformations of $\cF$.
\item The category $\cC\cV^0_{\K,\E}$ of Artin motives embeds as a
full subcategory of $\cE\cV^0_{\K,\E}$.
\item The composite functor
\begin{equation}\label{kFfunctor}
k\circ \cF: \cE\cV^0_{\K,\E} \to \cK\cK\otimes \E
\end{equation}
maps the full subcategory $\cC\cV^0_{\K,\E}$ of Artin motives
faithfully to the category $\cK\cK_{G,\E}$ of $G$-equivariant
$KK$-theory with coefficients in $\E$.
\end{enumerate}
\end{The}

\proof (1) The functor $\cF$ defined by taking points over $\bar\K$
maps
\begin{equation}\label{FfunctObj}
\cF: (X_\alpha,S)\mapsto (\cX,S)
\end{equation}
where $\cX=X(\bar\K)=\varprojlim_\alpha X_\alpha(\bar\K)$. If
$(X_\alpha,S)$ is a uniform algebraic endomotive then
$\cF(X_\alpha,S)=(X,S,\mu)$ is a measured analytic endomotive with
$\mu$ the projective limit of the normalized counting measures. On
correspondences the functor gives
\begin{equation}\label{FfunctCorr}
\cF: Z  \mapsto \cZ= Z(\bar\K).
\end{equation}
The obtained correspondence is etale since the corresponding right
$C^*$-module is finite and projective as an induced module of the
finite projective module of the algebraic situation. This is
compatible with composition of correspondences since in both cases
the composition is given by the fiber product over the middle
groupoid (\cf \eqref{corrfiberprodalg}).

\smallskip

(2) We know from \eqref{charaction} that the group $G$ acts by
automorphisms on $\cF(X_\alpha,S)$. This action is compatible with
the morphisms. This shows that $G$ acts on the functor $\cF$ by
natural transformations.

(3) The category of Artin motives is embedded in the category
$\cE\cV^0_{\K,\E}$ by the functor $\cJ$ that maps an Artin motive
$M$ to the system with $X_\alpha=M$ for all $\alpha$ and $S$
trivial. Morphisms are then given by the same geometric objects.

(4) Let $X=\Spec \; A$, $X'=\Spec \; B$ be $0$-dimensional
varieties. Given a component $Z=\Spec \; C$ of $X\times X'$ the two
projections turn $C$ into an $A-B$-bimodule $k(Z)$. Moreover, the
composition of correspondences translates into the tensor product of
bimodules \ie one has
\begin{equation}\label{intpro}
k(Z)\otimes_{B} k(L)\simeq k(Z\circ L)\,.
\end{equation}
We can then compose $k$ with the natural functor $A \to
A_\C=A\otimes_K \C$ which associates to an $A-B$-bimodule $\cE$ the
$A_\C-B_\C$-bimodule $\cE_\C=\,\cE\otimes_\K \C$. Thus, we view the
result as an element in $KK$.

Recall that we have, for any coefficients $\E$, the functor $\cJ$
which to a variety $X\in Obj(\cV_{\K,\E}^0)$ associates the system
$(X_\alpha,S)$ where the system $X_\alpha$ consists of a single $X$
and the semigroup $S$ is trivial. Let  $U\in
\Hom_{\cC\cV_{\K,\E}^0}(X,X')$  be given by a correspondence $U=\sum
a_i \chi_{Z_i}$, with coefficients in $\E$. We consider the element
in $KK$ (with coefficients in $\E$) given by the sum of bimodules
\begin{equation}\label{FUkk}
k(U)=\sum a_i k(Z_i).
\end{equation}

Using the isomorphism \eqref{embed} we thus obtain a $G$-equivariant
bimodule for the corresponding $C^*$-algebras, and hence owing to
finite dimensionality a corresponding class in $KK_G$. Thus and for
any coefficients $\E$ we get a faithful functor
\begin{equation}\label{HomHom}
k\circ \cF \left(\Hom_{\cC\cV_{\K,\E}^0}(X,X')\right)\subset
KK_G(k\cF(X),k\cF(X'))\otimes\E
\end{equation}
compatibly with the composition of correspondences.

Since a correspondence $U=\sum a_i \chi_{Z_i}$, with coefficients in
$\E$, is uniquely determined by the corresponding map
$$
K_0(A_\C)\otimes \E\;\to\;K_0(B_\C)\otimes \E
$$
we get the required faithfulness of the restriction of $k\circ \cF$
to Artin motives.
\endproof

\medskip

The fact of considering only a zero dimensional setting in Theorem
\ref{NCArtinThm} made it especially simple to compare the
composition of correspondences between Artin motives and
noncommutative spaces. This is a special case of a more general
result that holds in higher dimension as well and which we discuss
in Section \ref{Shigher} below. The comparison between
correspondences given by algebraic cycles and correspondences in
$KK$-theory is based in the more general case on the description
given in \cite{CoSka} of correspondences in $KK$-theory.

\medskip

The setting described in this section is useful in order to develop
some basic tools that can be applied to a class of ``zero
dimensional'' noncommutative spaces generalizing Artin motives,
which inherit natural analogs of arithmetic notions for motives. It
is the minimal one that makes it possible to understand the
intrinsic role of the absolute Galois group in dynamical systems
such as the BC-system. This set-up has no pretention at defining the
complete category of noncommutative Artin motives. This would
require in particular to spell out exactly the finiteness conditions
(such as the finite dimensionality of the center) that should be
required on the noncommutative algebra   $\cA$.

\medskip

We now turn to a general procedure which should play in
characteristic zero a role similar to the action of the Frobenius on
the $\ell$-adic cohomology. Its main virtue is its generality. It
makes essential use of {\em positivity} \ie of the key feature of
$C^*$-algebras. It is precisely for this reason that it was
important to construct the above functor $\cF$ that effects a bridge
between the world of Artin motives and that of noncommutative
geometry.

\section{Scaling as Frobenius in characteristic zero}\label{Frobenius}

In this section we describe a general cohomological procedure which
starting from a noncommutative space given by a pair $(\cA,
\varphi)$ of a unital  involutive algebra $\cA$ over $\C$ and a
state $\varphi$ on $\cA$ yields a representation of the
multiplicative group $\R_+^*$. We follow all the steps in a
particular example: the BC-system described above in Example
\ref{BCex}. We show that, in this case, the spectrum of the
representation is the set of non-trivial zeros of Hecke
$L$-functions with Gr\"ossencharakter.

What is striking is the  generality of the procedure and the analogy
with the role of the Frobenius in positive characteristic. The
action of the multiplicative group $\R_+^*$ on the cyclic homology
$HC_0$ of the distilled dual system plays the role of the action of
the Frobenius on the $\ell$-adic cohomology $H^1(C(\bar
\F_q),\Q_\ell)$. We explain in details that passing to the dual
system is the analog in characteristic zero of the transition from
$\F_q$ to its algebraic closure $\bar \F_q$. We also discuss the
analog of the intermediate step $\F_q \to \F_{q^n}$. What we call
the ``distillation" procedure gives a $\Lambda$-module
$D(\cA,\varphi)$ which plays a role similar to a motivic $H^1_{mot}$
(\cf \eqref{absmot} above). Finally, the dual action of $\R_+^*$ on
cyclic homology $HC_0(D(\cA,\varphi))$ plays a role similar to the
action of the Frobenius on the $\ell$-adic cohomology.

While here we only illustrate the procedure in the case of the
BC-system, in the next section we discuss it in the more general
context of global fields and get in that way a cohomological
interpretation of the spectral realization of zeros of Hecke
$L$-functions with Gr\"ossencharakter of \cite{AC}.

There is a deep relation outlined in Section \ref{IJ} between this
procedure and the classification of type III factors \cite{Co-th}.
In particular, one expects interesting phenomena when the
noncommutative space $(\cA, \varphi)$ is of type III and examples
abound where the general cohomological procedure should be applied,
with the variations required in more involved cases. Other type
III$_1$ examples include the following.

\begin{itemize}

\item The spaces of leaves of Anosov foliations and their codings.

\item The noncommutative space
$(\cA, \varphi)$ given by the restriction of the vacuum state on the
local algebra in the Unruh model of black holes.

\end{itemize}

One should expect that the analysis of the thermodynamics will get
more involved in general, where more that one critical temperature
and phase transitions occur, but that a similar distillation
procedure will apply at the various critical temperatures.

The existence of a ``canonically given" time evolution for
noncommutative spaces (see \cite{Co-th} and subsection \ref{IJ}
below) plays a central role in noncommutative geometry and the
procedure outlined below provides a very general substitute for the
Frobenius when working in characteristic zero.

\medskip
\subsection{Preliminaries}

Let $\cA$ be a unital involutive algebra over $\C$. We assume that
$\cA$ has a countable basis (as a vector space over $\C$) and that
the expression
\begin{equation}\label{normdef}
 ||x||=\,\sup \,||\pi(x)||
 \end{equation}
defines a {\em finite} norm on $\cA$, where $\pi$ in \eqref{normdef}
ranges through all unitary Hilbert space representations of the
unital involutive algebra $\cA$.

We let $\bar \cA$ denote the completion of $\cA$ in this norm. It is
a $C^*$-algebra.

This applies in particular to the unital involutive algebras
$\cA=\,C^\infty(X)\rtimes_{alg} S $ of \eqref{densesub} and the
corresponding $C^*$-algebra is $\bar \cA=\,C(X)\rtimes S$.

A state $\varphi$ on $\cA$ is a linear form $\varphi : \bar\cA
\rightarrow \C$ such that
$$
\varphi (x^* x) \geq 0, \ \ \forall \, x \in \bar \cA \ \ \text{ and
} \ \ \ \ \varphi (1) = 1 \,.
$$
One lets $\cH_\varphi$ be the Hilbert space of the Gelfand--Naimark--Segal
(GNS) construction,
\ie the completion of $\cA$ associated to the sesquilinear form
$$
\langle x,y \rangle = \varphi (y^* x) \, , \ x,y \in \cA \, .
$$

\begin{Def}\label{defregular}
A state $\varphi$ on a $C^*$-algebra $\bar\cA$ is {\em regular} iff
the vector $1\in \cH_\varphi$ is cyclic for the commutant of the
action of $\bar\cA$ by left multiplication on $\cH_\varphi$.
\end{Def}

Let $M$ be the von-Neumann algebra weak closure of the action of
$\bar\cA$ in $\cH_\varphi$ by left multiplication.

The main result of Tomita's theory \cite{tt} is that, if we let
$\Delta_{\varphi}$ be the {\it modular operator}
\begin{equation}\label{modop}
\Delta_{\varphi} = S_{\varphi}^* \, S_{\varphi} \, ,
\end{equation}
which is the {\it module} of the involution $ S_{\varphi}$, $x
\rightarrow x^*$, then we have
$$
\Delta_{\varphi}^{it}\,M\,\Delta_{\varphi}^{-it}=\,M \qqq t\in \R
\,.
$$
The corresponding one parameter group $ \sigma_t^{\varphi}\in
\Aut(M) $ is called the modular automorphism group of $\varphi$ and
it satisfies the KMS$_1$ condition relative to $\varphi$ (\cf
\cite{tt}). In general, given a unital involutive algebra $\cA$, a
one parameter group of automorphisms $\sigma_t\in \Aut(\cA) $ and a
state $\varphi$ on $\cA$, the KMS$_\beta$-condition is defined as
follows.

\begin{Def}\label{KMSbetadef}
A triple $(\cA, \sigma_t,\varphi)$ satisfies the
Kubo-Martin-Schwinger (KMS) condition at inverse temperature $0\leq
\beta < \infty$, if the following holds. For all $x,y\in \cA$, there
exists a holomorphic function $F_{x,y}(z)$ on the strip $0< {\rm
Im}(z)< \beta$, which extends as a continuous function on the
boundary of the strip, with the property that
\begin{equation}\label{KMSdef}
 F_{x,y}(t)=\varphi(x\sigma_t(y)) \ \ \ \text{ and } \ \ \
F_{x,y}(t+i\beta)=\varphi(\sigma_t(y)x), \ \ \ \ \forall t\in \R.
\end{equation}
\end{Def}

We also say that $\varphi$ is a KMS$_\beta$ state for \((\cA,
\sigma_t)\). The set $\Sigma_\beta$ of KMS$_\beta$ states is a
compact convex Choquet simplex \cite[II \S5]{BR} whose set of
extreme points $\cE_\beta$ consists of the factor states. One can
express any KMS$_\beta$ state uniquely in terms of extremal states,
because of the uniqueness of the barycentric decomposition of a
Choquet simplex.

\bigskip

\subsection{The cooling morphism}

Let $(\cA, \varphi)$ be a pair of an algebra and a state as above.
We assume that the modular automorphism group $\sigma_t^\varphi$
associated to the state $\varphi$ leaves the algebra $\cA$ globally
invariant and we denote by $\sigma_t\in \Aut(\cA)$ its restriction
to $\cA$.  By construction the state $\varphi$ is a KMS$_1$ state
for the quantum statistical system $(\bar \cA, \sigma)$. Here
$\bar\cA$ is the $C^*$-completion in \eqref{normdef} as above.

The process of cooling down the system consists of investigating
KMS$_\beta$-states for $\beta >1$ (\ie at lower temperatures).

The general theory of quantum statistical mechanics shows that the
space $\Sigma_\beta$ of KMS$_\beta$-states is a Choquet simplex. We
let $\cE_\beta$ denote the set of its extremal points. Each
$\varepsilon \in \cE_\beta$ is a factor state on $\cA$.

We require that, for sufficiently large $\beta$, all the $\epsilon
\in \cE_\beta$ are type I$_\infty$ factor states.

In general, the original state $\varphi$ will be of type III (it is
of type III$_1$ for the BC-system), so that the type I$_\infty$
property we required means that, under cooling, the system
simplifies and becomes ``commutative". In fact, under this
assumption, we obtain a ``cooled down dual algebra''
$C(\tilde\Omega_\beta,\cL^1)$, using the set $\Omega_\beta$ of
regular extremal KMS$_\beta$ states (\cf \eqref{cool} below), that
is Morita equivalent to a commutative algebra.

\begin{Lem}\label{typeIcond}
Let $\varphi$ be a regular extremal KMS$_\beta$ state on the
$C^*$-algebra $\bar \cA$ with the time evolution $\sigma_t$. Assume
that the corresponding factor $M_\varphi$ obtained through the GNS
representation $\cH_\varphi$ is of type I$_\infty$. Then there
exists an irreducible representation $\pi_\varphi$ of $\bar\cA$ in a
Hilbert space $\cH(\varphi)$ and an unbounded self-adjoint operator
$H$ acting on $\cH(\varphi)$, such that $e^{-\beta H}$ is of trace
class. Moreover, one has
\begin{equation}\label{state}
\varphi(x)=\, \Trace(\pi_{\varphi}(x)\,e^{-\beta
H})/\Trace(e^{-\beta H}) \,,\quad \forall x\in \bar\cA,
\end{equation}
and
\begin{equation}\label{imp}
e^{itH}\,\pi_\varphi(x)\,e^{-itH}=\,\pi_\varphi(\sigma_t(x)).
\end{equation}
\end{Lem}

\proof Let $1_\varphi\in \cH_\varphi$ be the cyclic and separating
vector associated to the regular state $\varphi$ in the GNS
representation. The factor $M_\varphi$ is of type I$_\infty$ hence
one has a factorization
\begin{equation}\label{factorH}
 \cH_\varphi = \cH(\varphi)\otimes \cH',  \ \ \ M_\varphi= \{
T\otimes 1 | T\in \cL(\cH(\varphi)) \}.
\end{equation}
The restriction of the vector state $1_\varphi$ to $M_\varphi$ is
faithful and normal, hence it can be uniquely written in the form
\begin{equation}\label{vectstate}
\langle (T\otimes 1) 1_\varphi, 1_\varphi \rangle = \Tr (T \rho),
\end{equation}
for a density matrix $\rho$, that is, a positive trace class
operator with $\Tr(\rho)=1$. This can be written in the form
$\rho=e^{-\beta H}$. It then follows that the vector state of
\eqref{vectstate} satisfies the KMS$_\beta$ condition on
$M_\varphi$, relative to the 1-parameter group implemented by
$e^{itH}$.

The factorization \eqref{factorH} shows that the left action of
$\bar\cA$ in the GNS representation $\cH_\varphi$ is of the form
$a\mapsto \pi_\varphi(a)\otimes 1$. The representation $\pi_\varphi$
on $\cH(\varphi)$ is irreducible.

Let $\tilde H$ be the generator of the 1-parameter group of unitary
operators $\sigma_t$ acting on $\cH_\varphi$. Since $M_\varphi$ is
% correction bar
the weak closure of $\bar\cA$ in $\cL(\cH_\varphi)$, it follows that
the extension of the state $\varphi$ to $M_\varphi$ fulfills the
KMS$_\beta$ condition with respect to the 1-parameter group
implemented by $e^{it\tilde H}$. This extension of $\varphi$ is
given by the vector state \eqref{vectstate}. This shows that the
1-parameter groups implemented by $e^{it\tilde H}$ and $e^{itH}$
agree on $M_\varphi$. This proves \eqref{imp}. Equation
\eqref{state} follows from \eqref{vectstate}, which implies in
particular that the operator $H$ has discrete spectrum bounded
below.
\endproof

In general, we let $\Omega_\beta\subset \cE_\beta$ be the set of
regular (in the sense of Definition \ref{defregular}) extremal
KMS$_\beta$ states of type I$_\infty$.

The equation \eqref{imp} does not determine $H$ uniquely, but only
up to an additive constant. For any choice of $H$, the corresponding
representation of the crossed product algebra $\hat \cA=\,\cA
\rtimes_{\sigma}\R$ is given by
\begin{equation}
\pi_{\varepsilon,H}(\int \,x(t)\,U_t\,dt)=\,\int
\,\pi_\varepsilon(x(t))\,e^{itH}\,dt .
\end{equation}
Let us define carefully the crossed product algebra $\hat \cA=\,\cA
\rtimes_{\sigma}\R$ by specifying precisely which maps $t\to x(t)\in
\cA$ qualify as elements of $\hat \cA$. One defines the Schwartz
space $\cS(\R,\cA)$ as the direct limit of the Schwartz space
$\cS(\R,V)$ over {\em finite dimensional} subspaces $V\subset \cA$ :
\begin{equation}\label{schwar}
\cS(\R,\cA)=\cup_V \,\cS(\R,V).
\end{equation}
Notice then that, since the vector space $\cA$ has countable basis
and the action of $\sigma_t$ is unitary, one can find a linear basis
of $\cA$ whose elements are eigenvectors for $\sigma_t$ and such
that the corresponding characters are unitary. It follows that the
Schwartz space $\cS(\R,\cA)$ is an algebra under the product
\begin{equation}\label{productalg}
(x \star y)(s)=\,\int\,x(t)\,\sigma_t(y(s-t))\,dt\,.
\end{equation}
We denote the algebra obtained in this way by
\begin{equation}\label{productalg1}
\hat \cA=\,\cA \rtimes_{\sigma_t}\R=\,\cS(\R,\cA).
\end{equation}

The dual action $\theta_\lambda \in \Aut(\hat \cA)$ of $\R^*_+$ is
given by
\begin{equation}\label{dualaction}
\theta_\lambda(\int \,x(t)\,U_t\,dt)=\,\int
\,\lambda^{it}\,x(t)\,U_t\,dt.
\end{equation}

By construction, the set $\tilde \Omega_\beta$ of pairs
$(\varepsilon,H)$ forms the total space of a principal $\R$-bundle
over $\Omega_\beta$, where the action of $\R$ is simply given by
translation of $H$ \ie in multiplicative terms, one has a principal
$\R_+^*$-bundle with
\begin{equation}
\lambda \,(\varepsilon,H)=\, (\varepsilon,H + \log \lambda) \,,\quad
\forall \lambda \in \R_+^*,
\end{equation}
which gives the fibration
\begin{equation}\label{fibration}
 \R_+^* \to \tilde \Omega_\beta \to \Omega_\beta\,.
\end{equation}
Notice that \eqref{fibration} admits a natural section given by the
condition
\begin{equation}\label{section}
 \Tr(e^{-\beta\,H})=\,1\,,
 \end{equation}
and a natural splitting
\begin{equation}\label{fibration1}
  \tilde \Omega_\beta \sim \Omega_\beta\times\, \R_+^* \,.
\end{equation}
The states $\varepsilon \in \Omega_\beta$ are in fact uniquely
determined by the corresponding representation $\pi_\varepsilon$,
thanks to the formula \eqref{state}, which does not depend on the
additional choice of $H$ fulfilling \eqref{imp}. The property
\eqref{state} also ensures that, for any function $f\in \cS(\R)$
with sufficiently fast decay at $\infty$, the operator $f(H)$ is of
trace class.

\medskip

We let $\hat \cA_\beta$ be the linear subspace of $\hat \cA$ spanned
by elements of the form $\int \,x(t)\,U_t\,dt$ where $x\in
\cS(I_\beta,\cA)$. Namely, $x\in \cS(\R,\cA)$ admits an analytic
continuation to the strip $I_\beta=\,\{z\,:\; \Im z\in [0,\beta]\}$,
which restricts to Schwartz functions on the boundary of $I_\beta$.

One checks that $\hat \cA_\beta$ is a subalgebra of $\hat \cA$.
Indeed the product \eqref{productalg} belongs to $\cS(I_\beta,\cA)$
provided that $x\in \cS(\R,\cA)$ and $y\in \hat \cA_\beta$. One then
gets the following.

\begin{Pro} \label{morphism}
\begin{enumerate}
\item For $(\varepsilon,H)\in
\tilde \Omega_\beta$, the representations $\pi_{\varepsilon,H}$ are
pairwise inequivalent irreducible representations of $\hat \cA=\,\cA
\rtimes_{\sigma}\R$.
\item Let $\theta_\lambda \in \Aut(\hat \cA)$ be the dual action
\eqref{dualaction}. This satisfies
\begin{equation}\label{equi}
\pi_{\varepsilon,H}\circ
\theta_\lambda=\,\pi_{\lambda(\varepsilon,H)} , \ \ \forall
(\varepsilon,H)\in \tilde \Omega_\beta .
\end{equation}
\item For $x\in \hat \cA_\beta$, one has
\begin{equation}
\pi_{\varepsilon,H}(\int
\,x(t)\,U_t\,dt)\in\,\cL^1(\cH(\varepsilon)) .
\end{equation}
\end{enumerate}
\end{Pro}

\begin{proof} (1) By construction the representations $\pi_\varepsilon$ of $\cA$ are
already irreducible. Let us show that the unitary equivalence class
of the representation $\pi_{\varepsilon,H}$ determines
$(\varepsilon,H)$. First the state $\varepsilon$ is uniquely
determined from \eqref{state}. Next the irreducibility of the
restriction to $\cA$ shows that only scalar operators can intertwine
$\pi_{\varepsilon,H}$ with $\pi_{\varepsilon,H + c}$, which implies
that $c=0$.

(2) One has $\pi_{\varepsilon,H}(U_t)=\,e^{itH}$, hence
$\pi_{\varepsilon,H}(\theta_\lambda(U_t))
=\,\lambda^{it}\,e^{itH}=\,\,e^{it(H+\log \lambda)}$, as required.

(3) Let us first handle the case of scalar valued $x(t)\in \C$ for
simplicity. One then has
$$
\pi_{\varepsilon,H}(\int \,x(t)\,U_t\,dt)=\,g(H) =\sum g(u) e(u),
$$
where the $e(u)$ are the spectral projections of $H$ and $g$ is the
Fourier transform in the form
$$
g(s)=\,\int \,x(t)\,e^{its}\,dt\,.
$$
Since $x\in \cS(\R)$ admits analytic continuation to the strip
$I_\beta=\,\{z\,:\; \Im z\in [0,\beta]\}$ and since we are assuming
that the restriction to the boundary of the strip $I_\beta$ is also
in Schwarz space, one gets that the function $g(s)\,e^{\beta s}$
belongs to $\cS(\R)$ and the trace class property follows from
$\Trace(e^{-\beta H})<\infty$. In the general case with  $x\in
\cS(I_\beta,\cA)$, with the notation above, one gets
$$
\pi_{\varepsilon,H}(\int \,x(t)\,U_t\,dt)=\, \int
\,\pi_{\varepsilon}(x(t)) \sum e^{its} e(s)\,dt =\,\sum
\,\pi_{\varepsilon}(g(s))\,e(s)
$$
where, as above, $e(s)$ is the spectral projection of $H$ associated
to the eigenvalue $s\in \R$. Again $||g(s)\,e^{\beta s}||<C$
independently of $s$ and $\sum \,e^{-\beta s}\,\Trace(e(s))<\infty$,
which gives the required trace class property.
\end{proof}

\medskip

In the case of the BC-system the above construction is actually
independent of $\beta>1$. Indeed, what happens is that the set of
irreducible representations of $\hat \cA$ thus constructed does not
depend on $\beta>1$ and the corresponding map at the level of the
extremal KMS$_\beta$ states is simply obtained by changing the value
of $\beta$ in formula \eqref{state}. In general the size of
$\Omega_\beta$ is a non-decreasing function of $\beta\;$ :
\medskip

\begin{Cor}\label{piCor} For any $\beta'>\beta$ the formula \eqref{state}
defines a canonical injection
\begin{equation}
c_{\beta',\beta}:\; \Omega_\beta \to \Omega_{\beta'}
\end{equation}
which extends to an $\R_+^*$-equivariant map of the bundles $\tilde
\Omega_\beta$.
\end{Cor}

\proof First $\Trace(e^{-\beta H})<\infty$ implies
$\Trace(e^{-\beta' H})<\infty$ for all $\beta'>\beta$ so that the
map $c_{\beta',\beta}$ is well defined. Notice that the obtained
state given by \eqref{state} for $\beta'$ is KMS$_{\beta'}$ extremal
(since it is factorial), regular and of type I. The map
$c_{\beta',\beta}$ is injective since one recovers the original
KMS$_\beta$ state by the reverse procedure. The latter in general
only makes sense for those values of $\beta$ for which
$\Trace(e^{-\beta H})<\infty$. This set of values is a half line
with lower bound depending on the corresponding element of
$\Omega_\infty = \cup_\beta \Omega_\beta$.
\endproof

\smallskip

We fix a separable Hilbert space $\cH$ and let, as in
\cite{dixmier}, $\Irr \hat \cA$ be the space of irreducible
representations of the $C^*$-algebra $\hat \cA$ in $\cH$ endowed
with the topology of pointwise weak convergence. Thus
\begin{equation}\label{repconv}
\pi_\alpha \to \pi \;\Longleftrightarrow \langle
\xi,\pi_\alpha(x)\eta\rangle \to \langle \xi,\pi(x)\eta\rangle
\qqq\, \xi,\eta \in \cH ,\, x \in \hat\cA\,.
\end{equation}
This is equivalent to pointwise strong convergence by \cite{dixmier}
3.5.2. By \cite{dixmier} Theorem 3.5.8 the map
\begin{equation}\label{irrtoprim}
\pi \in \Irr \hat \cA \mapsto \ker \pi \in \Prim \,\hat \cA
\end{equation}
is continuous and open. Indeed by Definition 3.1.5 of \cite{dixmier}
the topology on the spectrum, \ie the space of equivalence classes
of irreducible representations, is {\em defined} as the pull back of
the Jacobson topology of the primitive ideal space. By Proposition
\ref{morphism} each $\pi_{\varepsilon,H}$ defines an irreducible
representation of $\hat\cA$ whose range is the elementary
$C^*$-algebra of compact operators in $\cH(\varepsilon)$. Thus, by
Corollaries 4.10 and 4.11 of \cite{dixmier}, this representation is
characterized by its kernel which is  an element of $\Prim \hat\cA$.
The map $(\varepsilon,H) \mapsto \ker \pi_{\varepsilon,H}$ gives
 an $\R^*_+$-equivariant embedding
\begin{equation}\label{embedombeta}
\tilde\Omega_\beta\subset \Prim\,\hat\cA\,.
\end{equation}

We  assume that the representations $\pi_{\varepsilon,H}$, for
$(\varepsilon,H)\in \tilde \Omega_\beta$  can be continuously
realized as representations $\tilde\pi_{\varepsilon,H}$ in $\cH$ \ie
more precisely that the injective  map \eqref{embedombeta} is lifted
to a continuous map from $\tilde \Omega_\beta$ to $\Irr \hat \cA $
such that
\begin{equation}\label{liftingmap}
(\varepsilon,H)\mapsto \tilde\pi_{\varepsilon,H}(\hat f(H))
\end{equation}
 is a
continuous map  to $\cL^1=\cL^1(\cH)$ for any $f\in \cS(I_\beta)$.

\begin{Cor}\label{pi} Under the above conditions,
the representations $\pi_{\varepsilon,H}$, for $(\varepsilon,H)\in
\tilde \Omega_\beta$, assemble to  give an $\R_+^*$-equivariant
morphism
\begin{equation}\label{cool}
\pi \;:\; \hat \cA_\beta \to C(\tilde \Omega_\beta, \cL^1)
\end{equation}
defined by
\begin{equation}
\pi(x)(\varepsilon,H)=\,\pi_{\varepsilon,H}(x)\,,\quad \forall
(\varepsilon,H)\in \tilde \Omega_\beta\,.
\end{equation}
\end{Cor}

\medskip

This continuity can be checked directly in the case of the BC-system
where $\Omega_\beta=\,\cE_\beta$ is compact for the weak topology
and where all the representations $\pi_{\varepsilon,H}$ take place
in the same Hilbert space.

\subsection{The distilled $\Lambda$-module $D(\cA,\varphi)$}

The above cooling morphism of \eqref{cool} composes with the trace
$$\Trace : \cL^1\to \C$$ to yield a cyclic morphism from the
$\Lambda$-module associated to $\hat \cA_\beta$ to one associated to
a commutative algebra of functions on $\tilde \Omega_\beta$. Our aim
in this section is to define appropriately the cokernel of this
cyclic morphism, which gives what we call the distilled
$\Lambda$-module $D(\cA,\varphi)$.

\medskip

A partial trace gives rise to a cyclic morphism as follows.

\begin{Pro} \label{partial}
Let $\cB$ be a unital algebra. The equality
$$
\Tr((x_0 \otimes t_0)\otimes (x_1 \otimes t_1) \otimes \ldots
\otimes (x_n \otimes t_n))=\, x_0 \otimes x_1 \otimes \ldots \otimes
x_n\, \Trace(t_0\,t_1\,\cdots \,t_n)
$$
defines a map of $\Lambda$-modules from $(\cB\otimes
\cL^1)^\natural$ to $\cB^\natural$.
\end{Pro}

\proof  We perform the construction with the $\Lambda$-module of the
inclusion
$$
\cB\otimes \cL^1 \subset \cB\otimes \tilde \cL^1
$$
and obtain by restriction to $(\cB\otimes \cL^1)^\natural$ the
required cyclic morphism. The basic cyclic operations are given by
$$
\delta_i  (x^0 \otimes \ldots \otimes x^n) =\, x^0 \otimes \ldots
\otimes x^i x^{i+1} \otimes \ldots \otimes x^n\, , \quad 0 \leq i
\leq n-1 ,
$$
$$
\delta_n   (x^0 \otimes\ldots \otimes x^n) =\, x^n  x^0 \otimes x^1
\otimes \ldots \otimes x^{n-1},
$$
$$
\sigma_j   (x^0 \otimes \ldots \otimes x^n) =\, x^0 \otimes\ldots
\otimes x^j\otimes 1 \otimes x^{j+1} \otimes \ldots \otimes x^n,
\quad 0 \leq j \leq n ,
$$
$$
\tau_n   (x^0 \otimes \ldots \otimes x^n) =\, x^n \otimes x^0
\otimes\ldots \otimes x^{n-1}\,.
$$

\medskip

One has, for $0 \leq i \leq n-1$,
$$
\Tr\;\delta_i  ((x_0 \otimes t_0) \otimes \ldots \otimes (x_n
\otimes t_n)) = $$
$$  x_0 \otimes \ldots \otimes x_i x_{i+1} \otimes
\ldots \otimes x_n \, \Trace(t_0\,t_1\,\cdots \,t_n) =
$$
$$
\delta_i\;\Tr((x_0 \otimes t_0)\otimes (x_1 \otimes t_1) \otimes
\ldots \otimes (x_n \otimes t_n)) .
$$
Similarly,
$$
\Tr\;\delta_n ((x_0 \otimes t_0) \otimes \ldots \otimes (x_n \otimes
t_n)) = $$
$$  x_n\,x_0 \otimes \ldots  \otimes x_{n-1}\,
\Trace(t_n\,t_0\,t_1\,\cdots \,t_{n-1}) =
$$
$$
\delta_n( x_0  \otimes \ldots \otimes x_n)\,
\Trace(t_0\,t_1\,\cdots \,t_n)  = $$
$$  \delta_n\;\Tr((x_0 \otimes t_0)
\otimes \ldots \otimes (x_n \otimes t_n))  . $$

One has, for $0 \leq j \leq n$,
$$
\Tr\;\sigma_j  ((x_0 \otimes t_0) \otimes \ldots \otimes (x_n
\otimes t_n)) = $$
$$ x_0  \otimes \ldots\otimes x_j\otimes 1 \otimes
x_{j+1} \otimes \ldots \otimes x_n \,\, \Trace(t_0\,t_1\,\cdots
\,t_n) =
$$
$$
\sigma_j( x_0  \otimes \ldots \otimes x_n)\,
\Trace(t_0\,t_1\,\cdots \,t_n)= $$
$$  \sigma_j\;\Tr((x_0 \otimes t_0)
\otimes \ldots \otimes (x_n \otimes t_n)) . $$

Finally, one has
$$
\Tr\;\tau_n((x_0 \otimes t_0) \otimes \ldots \otimes (x_n \otimes
t_n)) = $$
$$ x_n  \otimes x_0  \otimes \ldots \otimes
x_{n-1}\,\Trace(t_0\,t_1\,\cdots \,t_n) =
$$
$$
\tau_n\;\Tr((x_0 \otimes t_0) \otimes \ldots \otimes (x_n
\otimes t_n)).
$$
\endproof

When dealing with non-unital algebras $\cA$, the maps $x \mapsto
x\otimes 1$ in the associated $\Lambda$-module $\cA^\natural$ are
not defined. Thus one uses the algebra $\tilde \cA$ obtained by
adjoining a unit, and defines $\cA^\natural$ as the submodule of
${\tilde \cA}^\natural$ which is the kernel of the augmentation
morphism,
$$
\varepsilon^\natural: {\tilde \cA}^\natural \to \C
$$
associated to the algebra homomorphism $\varepsilon:{\tilde \cA}\to
\C$.

We take this as a definition of $\cA^\natural$ for non unital
algebras. One obtains a $\Lambda$-module which coincides with $\cA$
in degree $0$ and whose elements in degree $n$ are tensors
$$
\sum \,a_0\otimes \cdots \otimes a_n\,,\quad a_j\in \tilde \cA
\,,\quad \sum \, \varepsilon(a_0) \cdots \varepsilon(a_n)=\,0
$$
While the algebra $\tilde \cA$ corresponds to the ``one point
compactification'' of the noncommutative space $\cA$, it is
convenient to also consider more general compactifications. Thus if
$\cA\subset \cA^{comp}$ is an inclusion of $\cA$ as an essential
ideal in a unital algebra $\cA^{comp}$ one gets a $\Lambda$-module
$(\cA,\cA^{comp})^\natural$ which in degree $n$ is given by tensors
$$
\sum \,a_0\otimes \cdots \otimes a_n\,,\quad a_j\in \cA^{comp} \,,
$$
where for each simple tensor $\,a_0\otimes \cdots \otimes a_n$ in
the sum, at least one of the $a_j$ belongs to $\cA$. One checks that
there is a canonical cyclic morphism
$$
\cA^\natural \to (\cA,\, \cA^{comp})^\natural
$$
obtained from the algebra homomorphism $\tilde \cA\to \cA^{comp}$.

\smallskip

Thus, we can apply the construction of Proposition \ref{partial} to
the case  $\cB=\,C(\tilde \Omega_\beta)$ under the hypothesis of
Corollary \ref{pi}. We shall briefly discuss in Remark
\ref{dixmierdouady} below how to adapt the construction in the
general case.

\medskip

The next step is to analyze the decay of the obtained functions on
$\tilde \Omega_\beta$. Let us look at the behavior on a fiber of the
fibration \eqref{fibration}. We first look at elements  $\int
\,x(t)\,U_t\,dt \in  \hat \cA_\beta$ with scalar valued $x(t)\in \C$
for simplicity. One then has, with $g(s)=\,\int
\,x(t)\,e^{its}\,dt$,
$$
\Trace( \pi_{\varepsilon,H + c}(\int \,x(t)\,U_t\,dt))=\,\Trace(g(H
+c))
$$
and we need to understand the behavior of this function when $c\to
\pm \infty$. Recall from the above discussion that both $g$ and the
function $g(s)\,e^{\beta s}$ belong to $\cS(\R)$.

\medskip

\begin{Lem} \label{plusdecay} (i)  For all $N>0$, one has $|\Trace(g(H
+c))|=\,O(e^{-\beta\,c}\,c^{-N})$, for $c\to +\infty$.

(ii) For any element $x \in \hat \cA_\beta$, one has
$$
|\Trace( \pi_{\varepsilon,H + c}(x))|=\,O(e^{-\beta\,c}\,c^{-N})
\,,\quad c\to +\infty\,.
$$
\end{Lem}

\begin{proof} (i) Let $g(s)=\,e^{-\beta\,s}\,h(s)$ where
$h\in \cS(\R)$. One has $|h(s+c)|\leq C \,c^{-N}$ for $s\in \Sp\,
H$, hence $|g(s+c)|\leq C \,e^{-\beta\,c}\,c^{-N}\,e^{-\beta\,s}$,
for $s\in \Sp \, H$, so that $|\Trace(g(H +c))|\leq C
\,e^{-\beta\,c}\,c^{-N}\,\Trace(e^{-\beta\,H})$.

(ii) In the general case with  $x\in \cS(I_\beta,\cA)$ and the
notation as in the proof of Proposition \ref{morphism}, one gets
$$
\Trace(\pi_{\varepsilon,H+c}(x))=\,\sum
\,\Trace(\pi_{\varepsilon}(g(s+c))\,e(s)),
$$
whose size is controlled by $\sum \,||g(s+c)||\,\Trace(e(s))$. As
above one has $||g(s+c)||\leq C
\,e^{-\beta\,c}\,c^{-N}\,e^{-\beta\,s}$, for $s\in \Sp \, H$, so
that using $\sum \,\Trace(e(s))\,e^{-\beta\,s}<\infty$ one gets the
required estimate.
\end{proof}

\medskip
Let us now investigate the behavior of this function when $c\to -
\infty$. To get some feeling, we first take the case of the
BC-system. Then the spectrum of $H$ is the set $\{\log(n)\,:\,n\in
\N=\Z_{>0} \}$, so that $\Trace(g(H +c))=\sum_{n\in
\N^*}\,g(\log(n)+c)$. Let $f(x)=\,g(\log(x))$ be extended to $\R$ as
a continuous even function. Then one has $f(0)=0$ and, for
$\lambda=\,e^c$,
$$
\Trace(g(H +c))=\frac{1}{2}\,\sum_{n\in \Z}\,f(\lambda \,n)\,.
$$
The Poisson summation formula therefore gives, with the appropriate
normalization of the Fourier transform $\hat f$ of $f$,
$$
\lambda \,\sum_{n\in \Z}\,f(\lambda \,n)=\sum_{n\in \Z}\,\hat
f(\lambda^{-1} \,n).
$$
This shows that, for $c\to -\infty$, \ie for $\lambda \to 0$, the
leading behavior is governed by
$$
\sum_{n\in \Z}\,f(\lambda \,n)\sim \lambda^{-1}\int\,f(x)\,dx .
$$
This was clear from the Riemann sum approximation to the integral.
The remainder is controlled by the decay of $\hat f$ at $\infty$,
\ie by the smoothness of $f$. Notice that the latter is not
guaranteed by that of $g$ because of a possible singularity at $0$.
One gets in fact the following result.
\medskip

\begin{Lem} \label{bcdecay}
In the BC case, for any $\beta >1$ and any element $x\in
\hat\cA_\beta$,  one has
\begin{equation}\label{bcbc}
\Trace( \pi_{\varepsilon,H + c}(x))=\,e^{-c}\,\tau(x) +\,O(|c|^{-N})
\,,\quad c\to -\infty\,,\quad \forall N>0 ,
\end{equation}
where $\tau$ is the canonical dual trace on $\hat \cA$.
\end{Lem}

We give the proof after Lemma \ref{Emap}. Combining Lemmata
\ref{plusdecay} and \ref{bcdecay} one gets, in the BC case and for
$x\in \hat\cA_\beta$, $\tau(x)=0$, that the function
\begin{equation}\label{feps} f_\varepsilon
(c)=\,\Trace( \pi_{\varepsilon,H + c}(x))
\end{equation}
has the correct decay for $c\to \pm \infty$ in order to have a
Fourier transform (in the variable $c$) that is holomorphic in the
strip $I_\beta$. Indeed one has, for all $N>0$,
\begin{equation}\label{decaycond}
\begin{array}{ll}
|f_\varepsilon (c)|=\,O(e^{-\beta\,c}\,|c|^{-N}) , & \text{ for } \
\ c\to +\infty \\[2mm]
|f_\varepsilon (c)|=\,O(|c|^{-N}) & \text{ for } \ \  c\to
-\infty\,.
\end{array}
\end{equation}

This implies that the Fourier transform
\begin{equation}\label{fouriernorm}
F(f_\varepsilon )(u)=\, \int\,f_\varepsilon (c)\,e^{-icu}\,dc
\end{equation}
is holomorphic in the strip $I_\beta$, bounded and smooth on the
boundary. In fact, one has $F(f_\varepsilon )\in \cS(I_\beta)$, as
follows from Lemma \ref{covarpi} below, based on the covariance of
the map $\Tr\circ \pi$.

We let $\Hol(I_\beta)$ be the algebra of multipliers of
$\cS(I_\beta)$, \ie of holomorphic functions $h$ in the strip
$I_\beta$ such that $h\,\cS(I_\beta)\subset \cS(I_\beta)$. For such
a function $h$ its restriction to the boundary component $\R$ is a
function of tempered growth and there is a unique distribution $\hat
h$ such that
$$
\int \,\hat h(\lambda)\,\lambda^{it}\,d^*\lambda=\,h(t) .
$$

One then has in full generality the following result.

\begin{Lem} \label{covarpi} 1) For $h\in \Hol(I_\beta)$ the
operator $\theta(\hat h)=\,\int \,\hat h(\lambda)\,\theta_\lambda
\,d^*\lambda$ acts on $\hat \cA_\beta$. This turns $\hat \cA_\beta$
into a module over the algebra $\Hol(I_\beta)$.

2) The map $F\circ \Tr\circ \pi$ is an $\Hol(I_\beta)$-module map.
\end{Lem}

\proof 1) One has
$$
\theta(\hat h)(\int \,x(t)\,U_t\,dt)=\,\int \int\,\hat
h(\lambda)\,\lambda^{it} \,x(t)\,U_t\,dt\,d^*\lambda=\,\int
\,h(t)\,x(t)\,U_t\,dt
$$
and $h\,x\in \cS(I_\beta,\cA)$.

2) The equivariance of $\pi$, \ie the equality (\cf \eqref{equi})
$$
\pi_{\varepsilon,H}\circ \theta_\lambda=\,\pi_{(\varepsilon,H+\log
\lambda)} \qqq (\varepsilon,H)\in \tilde \Omega_\beta,
$$
gives
$$
(\Tr\circ \pi\circ \theta_\lambda)(x)_{(\varepsilon,H)}=\,(\Tr\circ
\pi)(x)_{(\varepsilon,H+\log \lambda)} .
$$

\medskip
Thus, with $f_\varepsilon$ as in \eqref{feps}, one gets
$$
(\Tr\circ \pi\circ \theta(\hat h))(x)_{(\varepsilon,H)}=\,\int
\,\hat h(\lambda)\,f_\varepsilon(\log \lambda)\,d^*\lambda .
$$
Replacing $H$ by $H+c$ one gets that the function associated to
$\theta(\hat h)(x)$ by \eqref{feps} is given by
$$
g_\varepsilon(c)=\,\int \,\hat h(\lambda)\,f_\varepsilon(\log
\lambda + c)\,d^*\lambda .
$$
One then gets
$$
F(g_\varepsilon)(u)=\,\int\,g_\varepsilon (c)\,e^{-icu}\,dc=\,\int
\int \,\hat h(\lambda)\,f_\varepsilon(\log \lambda +
c)\,e^{-icu}\,dc\,d^*\lambda
$$
$$
=\,\int \,\hat
h(\lambda)\,\lambda^{iu}\,d^*\lambda\,\int\,f_\varepsilon
(a)\,e^{-iau}\,da=\, h(u)\,F(f_\varepsilon)(u)\,.
$$
\endproof

This makes it possible to improve the estimate \eqref{decaycond} for
the BC-system immediately. Indeed, the function $h(s)=s$ belongs to
$\Hol(I_\beta)$ and this shows that \eqref{decaycond} holds for all
derivatives of $f_\varepsilon$. In other words it shows the
following.

\begin{Cor} \label{bcdecay1}
In the BC case, for any $\beta >1$ and any element $x\in
\hat\cA_\beta\cap \Ker \,\tau$,  one has
$$
F\circ \Tr\circ \pi(x)\in C(\Omega_\beta,\cS(I_\beta))\,.
$$
\end{Cor}

Thus, by Lemma \ref{covarpi} we obtain an $\Hol(I_\beta)$-submodule
of $C(\Omega_\beta,\cS(I_\beta))$ and as we will see shortly that
this contains all the information on the zeros of $L$-functions with
Gr\"ossencharakters over the global field $\Q$.

When expressed in the variable $\lambda=e^c$, the regularity of the
function $\Tr\circ \pi(x)$ can be written equivalently in the form
$\Tr\circ \pi(x)\in C(\Omega_\beta,\cS_\beta(\R_+^*))$ where, with
$\mu(\lambda)=\lambda, \,\forall \lambda \in \R_+^*$, we let
\begin{equation}\label{reglambda}
\cS_\beta(\R_+^*)=\,\cap_{[0,\beta]}\,\mu^{-\alpha}\,\cS(\R_+^*)
\end{equation}
be the intersection of the shifted Schwartz spaces of $\R_+^*\sim
\R$.

We let $\Sc^\natural(\tilde\Omega_\beta)$ be the cyclic submodule of
$C(\tilde \Omega_\beta)^\natural$ whose elements are functions with
restriction to the main diagonal that belongs to
$C(\Omega_\beta,\cS_\beta(\R_+^*))$.  We then obtain the following
result.

\begin{Pro}\label{beta0prop}
Assuming  \eqref{bcbc}, the map $\Tr\circ \pi$ defines a cyclic
morphism
\begin{equation}\label{Abeta0}
\delta\;:\; \hat \cA_{\beta,0}^\natural \to
\,\Sc^\natural(\tilde\Omega_\beta)\,.
\end{equation}
\end{Pro}

The subscript $0$ in $\hat \cA_{\beta,0}^\natural$
of \eqref{Abeta0} means that we are considering the
cyclic module defined as in \eqref{cyc} below.
The proof of Proposition \ref{beta0prop}
follows from Proposition \ref{partial}.

\medskip

\begin{Def}
We define the distilled $\Lambda$-module $D(\cA,\varphi)$ as the
cokernel of the cyclic morphism $\delta$.
\end{Def}

\medskip

\begin{Rem}\label{dixmierdouady}{\rm
In the general case \ie without the lifting hypothesis of Corollary
\ref{pi}, one is dealing (instead of the representations
$\tilde\pi_{\varepsilon,H}$ in a fixed Hilbert space) with the
continuous field of elementary $C^*$-algebras
$C_{\varepsilon,H}=\pi_{\varepsilon,H}(\cA)$ on $\tilde
\Omega_\beta$. Even though  the algebras $C_{\varepsilon,H}$ are
isomorphic to the algebra of compact operators in $\cH$ (or in a
finite dimensional Hilbert space) the Dixmier-Douady invariant (\cf
\cite{dixmier} Theorem 10.8.4) gives an obstruction to the Morita
equivalence between the algebra $C(\tilde \Omega_\beta,C)$ of
continuous sections of the field and the algebra of continuous
functions on $\tilde \Omega_\beta$. When this obstruction is
non-trivial it is no longer possible to use the cyclic morphism
coming from the Morita equivalence and one has to work with
$C(\tilde \Omega_\beta,C)$ rather than with the algebra of
continuous functions on $\tilde \Omega_\beta$. The cooling morphism
is now the canonical map
\begin{equation}\label{coolingmorphism1}
\begin{array}{c}
\pi : \hat \cA_\beta \to C(\tilde \Omega_\beta, C) \\[3mm]
\pi(x)(\varepsilon,H)=\tilde\pi_{\varepsilon,H}(x), \ \ \ \ \forall
(\varepsilon,H)\in \tilde \Omega_\beta .
\end{array}
\end{equation}
which defines an $\R_+^*$-equivariant algebra homomorphism. The
cokernel of this morphism  as well as its cyclic cohomology continue
to make sense but are more complicated to compute. Moreover for a
general pair $(\cA,\varphi)$ the decay condition \eqref{bcbc} will
not be fulfilled. One then needs to proceed with care and analyze
the behavior of these functions case by case. One can nevertheless
give a rough general definition, just by replacing $\cS\subset
C(\tilde \Omega_\beta)$ above by the subalgebra generated by
$\Tr\circ \pi(\Ker \,\tau)$.}
\end{Rem}

\bigskip
\subsection{Dual action on cyclic homology}

By construction, the cyclic morphism $\delta$ is equivariant for the
dual action $\theta_\lambda$ of $\R^*_+$ (\cf  \eqref{equi}). This
makes it possible to consider the corresponding representation of
$\R^*_+$ in the cyclic homology group, \ie
\begin{equation}\label{cohomrep}
\theta(\lambda) \in \Aut(HC_0(D(\cA,\varphi))) \qqq \lambda \in
\R^*_+ .
\end{equation}

The spectrum of this representation defines a very subtle invariant
of the original noncommutative space $(\cA,\varphi)$.

Notice that, since the $\Lambda$-module
$\Sc^\natural(\tilde\Omega_\beta)$ is commutative, so is the
cokernel of $\delta$ (it is a quotient of the above). Thus, the
cyclic homology group $HC_0(D(\cA,\varphi))$ is simply given by the
cokernel of $\delta$ in degree $0$.

It might seem unnecessary to describe this as $HC_0$ instead of
simply looking at the cokernel of $\delta$ in degree $0$, but the
main reason is the Morita invariance of cyclic homology which wipes
out the distinction between algebras such as $B$ and $B\otimes
\cL^1$.

We now compute the invariant described above in the case of the
BC-system described in Example \ref{BCex} of Section \ref{Sartin}.
Recall that one has a canonical action of the Galois group
$G=\Gal(\bar \Q/\Q)$ on the system, which factorizes through the
action on the cyclotomic extension $\Q^{ab}$, yielding an action of
$G^{ab}=\,\hat\Z^*$ as symmetries of the BC-system. Thus, one
obtains a decomposition as a direct sum using idempotents associated
to Gr\"ossencharakters.

\begin{Pro}  Let $\chi$ be a character of the compact group $\hat \Z^*$.
Then the expression
$$
p_\chi=\,\int_{\hat \Z^*}\,g\,\chi(g)\,dg
$$
determines an idempotent $p_\chi$ in $\End_\Lambda D(\cA,\varphi)$.
\end{Pro}

Notice that the idempotent $p_\chi$ already makes sense in the
category $\cE\cV_{K,\C}^0$ introduced above in Section \ref{Sartin}.

The idempotents $p_\chi$ add up to the identity, namely
\begin{equation}
\sum_{\chi}\,p_\chi=\, {\rm Id},
\end{equation}
and we get a corresponding direct sum decomposition of the
representation \eqref{cohomrep}.

\medskip

\begin{The}\label{specreal}
The representation of $\R^*_+$ in
\begin{equation}
\Mc=\,HC_0(p_\chi\,D(\cA,\varphi))
\end{equation}
gives the spectral realization of the zeros of the $L$-function
$L_\chi$. More precisely, let $z\in I_\beta$ be viewed as a
character of $\Hol(I_\beta)$ and $\C_z$ be the corresponding one
dimensional module. One has
\begin{equation}
\Mc\otimes_{\Hol(I_\beta)}\C_z \neq \{0\} \Longleftrightarrow
L_\chi(-i\, \,z)=\,0.
\end{equation}
\end{The}

\medskip

One can give a geometric description  of the BC-system  in terms  of
$\Q$-lattices in $\R$ (\cf \cite{CM}). The $\Q$-lattices in $\R$ are
labeled by pairs$(\rho,\lambda)\in \hat \Z \times \R_+^*$ and the
commensurability equivalence relation is given by the orbits of the
action of the semigroup $\N=\Z_{>0}$ given by
$n\,(\rho,\lambda)=\,(n\,\rho,n\,\lambda)$. The $\Q$-lattice
associated to the pair $(\rho,\lambda)$ is $(\lambda^{-1}\Z,
\lambda^{-1}\rho)$ where $\rho$ is viewed in $\Hom(\Q/\Z,\Q/\Z)$.
The BC-system itself corresponds to  the noncommutative space of
commensurability classes of $\Q$-lattices up to scale, which
eliminates the $\lambda$ (\cf \cite{CM}). The $\Q$-lattices up to
scale are labelled by $\rho \in \hat \Z$ and the commensurability
equivalence relation is given by the orbits of the multiplicative
action of the semigroup $\N$. The algebra $\cA$ of the BC-system is
the crossed product $C^\infty(\hat \Z)\rtimes \N$ where
$C^\infty(\hat \Z)$ is the algebra of locally constant functions on
$\hat \Z $. By construction every $f\in C^\infty(\hat \Z)$ has a
level \ie can be written as
\begin{equation}\label{level}
f=\,h\circ p_N\,,\quad p_N\;:\;\hat \Z\to \Z/N\Z .
\end{equation}
General elements of $\cA$ are finite linear combinations of
monomials, as in \eqref{monom}. Equivalently one can describe $\cA$
as a subalgebra of  the convolution algebra of the \'etale groupoid
$G_{bc}$ associated to the partial action of $\Q_+^*$ on $\hat \Z $.
One has $G_{bc}=\,\{(k,\rho)\in \Q_+^*\times \hat \Z\,|\,k\,\rho \in
\hat \Z\}$ and $(k_1,\rho_1)\circ (k_2,\rho_2)=(k_1 k_2,\rho_2)$ if
$\rho_1=\,k_2\,\rho_2$. The product in the groupoid algebra is given
by the associative convolution product
\begin{equation}\label{convol1}
 f_1 * f_2 \, (k,x)= \sum_{s\in \Q_+^*} f_1(k\,s^{-1}, s\,x)
f_2(s,x),
\end{equation}
and the adjoint is given by $f^*(k,x)=\overline{f(k^{-1},k\,x)}$.
The elements of the algebra $\cA$ are functions $f(k,\rho)$ with
finite support in $k \in \Q_+^*$ and finite level (\cf
\eqref{level}) in the variable $\rho$. The time evolution is given
by
\begin{equation}\label{sigma}
\sigma_t(g)(k,\rho)=\,k^{it} \,g(k,\rho)\qqq g \in \cA\,.
\end{equation}

\medskip

The dual system of the BC-system corresponds to the noncommutative
space of commensurability classes of $\Q$-lattices in $\R$ (\cf
\cite{CM}). Thus, the dual algebra $\hat \cA$ for the BC-system is a
subalgebra of the convolution algebra of the \'etale groupoid $\hat
G_{bc}$ associated to the partial action of $\Q_+^*$ on $\hat \Z
\times \R_+^*$. The product in the groupoid algebra is given by the
same formula as \eqref{convol1} with $x=(\rho,\lambda)\in \hat \Z
\times \R_+^*$.

Let $x=\,\int \,x(t)\,U_t\,dt \in  \hat \cA$. Then the corresponding
function $f=\,\iota(x)$ on $\hat G_{bc}$ is given explicitly by
\begin{equation}\label{iso1}
f(k,\rho,\lambda)=\,\int \,x(t)(k,\rho)\,\lambda^{it}\,dt .
\end{equation}
One checks that $\iota$ is an algebra homomorphism using
\eqref{sigma}. Using the definition \eqref{schwar} of the Schwartz
space $\cS(\R,\cA)$, one gets that elements of $\iota\hat \cA$ are
functions $f(k,\rho,\lambda)$ with finite support in $k$, finite
level in $\rho$, and rapid decay (Schwartz functions) on $\R_+^*\sim
\R$. Similarly, elements of $\iota\hat \cA_\beta$ are functions
$f(k,\rho,\lambda)$ with finite support in $k$, finite level in
$\rho$, and such that the finitely many scalar functions
$$
f_{k,\rho}(\lambda)=\,f(k,\rho,\lambda) \qqq \lambda \in \R_+^*
$$
are in the Fourier transform of the space $\cS(I_\beta)$,
\begin{equation}\label{characabeta}
f \in \iota\hat \cA_\beta\,\Longleftrightarrow \tilde
F(f_{k,\rho})\in \cS(I_\beta) \qqq (k,\rho).
\end{equation}
Here the Fourier transform $\tilde F$ in the variable $\lambda \in
\R_+^*$ corresponds to \eqref{fouriernorm}  in the variable
$c=\log(\lambda)$, \ie
\begin{equation}\label{fouriernorm1}
\tilde F(g )(u)=\, \int\,g(\lambda)\,\lambda^{-iu}\,d^*\lambda .
\end{equation}

\medskip
Up to Morita equivalence $\hat G_{bc}$ is the same as the groupoid
of the partial action of $\Q^*$ on $\hat \Z \times \R^*$. It is not
quite the same as the adele class space \ie the action of $\Q^*$ on
$\A_\Q$, not because of the strict inclusion $\hat \Z \subset \A_f$
of $\hat \Z $ in finite adeles, which is a Morita equivalence, but
because we are missing the compact piece $\hat \Z \times
\{0\}\subset \A_\Q$. Modulo this nuance, one can compare the
restriction map $\delta$ with the map $E$ of \cite{AC} given (up to
normalization by $|j|^{1/2}$) by
\begin{equation}
\label{QbigE} E ( f) (j) =  \sum_{q\in \Q^*}  f (q \,j)\,,\quad j
\in C_\Q\,,\quad f\in \cS(\A_\Q)\, ,
\end{equation}
where $\A_\Q$ are the adeles of $\Q$ and $C_\Q$ the idele class
group.

\medskip

\begin{Lem} \label{Emap}
1) For $\beta >1$ one has a canonical isomorphism $\tilde
\Omega_\beta \sim \hat \Z^*\times \R_+^* \sim C_\Q$ of $\tilde
\Omega_\beta $ with the space of invertible $\Q$-lattices.

2) Let $x\in \hat \cA$ and $f=\,\iota(x)$ be the corresponding
function on $\hat G_{bc}\subset \Q_+^*\times \hat \Z \times \R_+^*$.
Then
$$
\delta(x)(j)=\,\sum_{n\in \N^*}  f (1, n \,u,n\,\lambda)\qqq j
=(u,\lambda)\in \hat\Z^*\times \R_+^*=C_\Q\, .
$$
\end{Lem}

\begin{proof} 1) Let us give explicitly the covariant irreducible
representation of $\cA$ associated to a pair $(u,\lambda)\in
\hat\Z^*\times \R_+^*$. The algebra acts on the fixed Hilbert space
$\cH=\,\ell^2(\N^*)$ with
$$
(\pi_{(u,\lambda)}(x)\,\xi)(n)=\,\sum_{m\in
\N^*}\,x(n\,m^{-1},\,m\,u)\,\xi(m) \qqq x\in \cA,
$$
while $\pi_{(u,\lambda)}(H)$ is the diagonal operator $D_\lambda$ of
multiplication by $\log n + \log \lambda$.

The classification of KMS$_\beta$ states of the BC system
(\cite{CM}) shows that the map $(u,\lambda)\to \pi_{(u,\lambda)}$ is
an isomorphism of $\,\hat\Z^*\times \R_+^*$ with $\tilde
\Omega_\beta$ for $\beta>1$.

2) By construction, the function $\delta(x)(j)$ is given, for $j
=(u,\lambda)\in \hat\Z^*\times \R_+^*$, by $ \delta(x)(j)=\,\Trace
(\pi_{(u,\lambda)}(x)) $. One has
$$
\pi_{(u,\lambda)}(x)=\,\int \,\pi_{(u,\lambda)}(x(t))\,e^{it
D_\lambda} \,dt
$$
and the trace is the sum of the diagonal entries of this matrix.
Namely, it is given by
$$
\sum_{n\in \N^*} \,\int \,x(t)(1,nu)\,e^{it (\log n + \log \lambda)}
\,dt =\, \sum_{n\in \N^*} \,f(1,n\,u,n\,\lambda) ,
$$
where we used \eqref{iso1}.
\end{proof}

{\bf Proof of Lemma \ref{bcdecay}}. Using Lemma \ref{Emap} we need
to show that for any $u\in \hat\Z^*$ and any function
$f(k,\rho,\lambda)$ in $\iota\hat \cA_\beta$, one has when
$\lambda\to 0$,
\begin{equation}\label{estimate}
\sum_{n\in \N^*}  f (1, n
\,u,n\,\lambda)=\,\lambda^{-1}\,\int\,f(1,\rho,v)\,d\rho\, dv
+\,O(|\log(\lambda)|^{-N})
\end{equation}
for all $N>0$. One checks indeed that the dual trace (\cf \cite{Tak}) is
given explicitly by the additive Haar measure
$$
\tau(f)=\,\int\,f(1,\rho,v)\,d\rho\, dv
$$
so that \eqref{estimate} implies Lemma \ref{bcdecay}. We only
consider the function of two variables $f(\rho,v)=\,f(1,\rho,v)$ and
we let $\tilde f$ be its unique extension to adeles by $0$ by
setting $\tilde f$ equal to zero outside of $\hat \Z\times
\R^+\subset \A_\Q$.
 Thus, the left hand
side of \eqref{estimate} is given by
\begin{equation}\label{estimate1}
\sum_{n\in \N^*}  f (1, n \,u,n\,\lambda)=\,\sum_{q\in \Q^*} \tilde
f (q \,j)\,,\quad j=(u,\lambda) .
\end{equation}
If the extended function $\tilde f$ were in the Bruhat-Schwartz space
$\cS(\A_\Q)$, then one could use the Poisson summation formula
$$
\vert j \vert \sum_{q \in \Q} h(q\,j) = \sum_{q \in \Q} \wh h
(q\,j^{-1}) ,
$$
with $\wh h$ the Fourier transform of $h$, in order to get a better
estimate with a remainder in $O(|\lambda|^{N})$. However, in general one
has $\tilde f\notin \cS(\A_\Q)$,  because of the singularity at
$\lambda=0$. At least some functions $f$ with non-zero integral do
extend to an element $\tilde f\in \cS(\A_\Q)$ and this allows us to
restrict to functions of two variables $f(\rho,v)$ such that
\begin{equation}\label{intzero}
 \int\,f(\rho,v)\,d\rho\,
dv=0\,.
\end{equation}
We want to show that, if $f(\rho,v)=g(p_n(\rho),v)$ with $g$ a
map from $\Z/n\Z$ to $\cS(\R_+^*)$, and if \eqref{intzero} holds,
then one has
\begin{equation}\label{est1}
E(\tilde f)(j)=\sum_{q\in \Q^*} \tilde f (q
\,j)=\,O(|\log|j|\;|^{-N})
\end{equation}
for $|j|\to 0$ and for all $N>0$.
Furthermore, we can assume that there exists a
character $\chi_n$ of the multiplicative group $(\Z/n\Z)^*$ such
that
$$
g(m \,a, v)=\,\chi_n(m)\,g(a,v)\qqq a\in \Z/n\Z\,,\quad m\in
(\Z/n\Z)^*\,.
$$
With $\chi=\chi_n\circ p_n$, one then has
\begin{equation}\label{eigen}
\tilde f(m \,a, v)=\,\chi(m)\,\tilde f(a,v)\qqq a\in\hat \Z\,,\quad
m\in \hat\Z^*\,.
\end{equation}
In particular, one gets
\begin{equation}\label{eigen1}
E(\tilde f)(m\,j)=\,\chi(m)\,E(\tilde f)(j) \ \ \ \forall m\in
\hat\Z^*,\ \forall j\in C_\Q\,.
\end{equation}
This shows that, in order to prove \eqref{est1}, it is enough to prove
the estimate
\begin{equation}\label{est0}
|h(\lambda)|=\,O(|\log|\lambda|\;|^{-N}), \ \ \text{ for } \ \lambda
\to 0 ,
\end{equation}
where
$$
h(\lambda)=\,\int_{\hat\Z^*}\,E(\tilde f)(m\,j)\,\chi_0(m)\,d^*m
\qqq j\ \ \text{ with } \ |j|=\lambda,
$$
and with $\chi_0=\bar \chi$,  the conjugate of $\chi$. Using
\eqref{fouriernorm1}, the Fourier
transform of $h$ (viewed as a function on $\R_+^*$) is given by
\begin{equation}\label{fouriertrans}
\tilde F h(t)=\,
\int\,h(\lambda)\,\lambda^{-it}\,d^*\lambda=\,\int_{\A_\Q^*}\,\tilde
f(j)\,\chi_0(j)\,|j|^s\,d^*j \,,\ \ \text{ with } \ s=-it\,,
\end{equation}
where $d^*j$ is the multiplicative Haar measure on the ideles and
$\chi_0$ is extended as a unitary character of the ideles. More
precisely, we identify  $\A_\Q^*=\hat\Z^*\times \R_+^*\times \Q^*$
where $\Q^*$ is embedded as principal ideles. We extend $\chi_0$ by
$1$ on $\R_+^*\times \Q^*\subset \A_\Q^*$. Thus $\chi_0$ becomes
a unitary Gr\"ossencharakter.
By Weil's results on homogeneous
distributions (\cf \cite{AC} Appendix I, Lemma 1), one has, with
$\tilde f$ as above, the equality
\begin{equation}\label{weildis}
\int_{\A_\Q^*}\,\tilde f(j)\,\chi_0(j)\,|j|^s\,d^*j =\,L(\chi_0 ,
s)\,\langle D'(s),\tilde f\rangle \,,\ \ \ \text{ for } \ \Re(s)>1 ,
\end{equation}
where $D'=\otimes \,D'_v$ is an infinite tensor product over the
places of $\Q$ of homogeneous distributions on the local fields
$\Q_v$ and $L(\chi_0 , s)$ is the $L$-function with
Gr\"ossencharakter $\chi_0$. This equality is a priori only valid
for $\tilde f \in \cS(\A_\Q)$, but it remains valid for $\tilde f$ by
absolute convergence of both sides. In fact, we can assume that
$f(\rho,v)=f_0(\rho) f_\infty(v)$, where $f_\infty\in \cS(\R_+^*)$.
One then gets $\tilde f=\tilde f_0\otimes \tilde f_\infty$ and
$$
\langle D'(s),\tilde f\rangle=\, \langle D'_0(s),\tilde
f_0\rangle\,\langle D'_\infty(s),\tilde f_\infty\rangle ,
$$
where by construction (\cf \cite{AC}) one has
\begin{equation}\label{dprimeinf}
\langle D'_\infty(s),\tilde f_\infty \rangle=\,
\int_{\R_+^*}\,f_\infty(\lambda)\,\lambda^s d^*\lambda \,,\quad
d^*\lambda=\,d\lambda/\lambda .
\end{equation}
Using \eqref{characabeta}, one gets $\tilde
F(f_\infty)\in\cS(I_\beta)$. Thus,  the function $t \mapsto \langle
D'_\infty(-i t),\tilde f_\infty \rangle$ is in $\cS(I_\beta)$. The
term $\langle D'_0(s),\tilde f_0\rangle$ is simple to compute
explicitly and the function $t \mapsto \langle D'_0(-i t),\tilde
f_0\rangle$ belongs to $\Hol(I_\beta)$.

Indeed, it is enough to do the computation when $f_0$ is a finite
tensor product $f_0=\otimes_{v\in S}  f_p $ over a finite set $S$ of
primes, while $f_p$ is a locally constant function on $\Q_p$
vanishing outside $\Z_p$ and such that $f_p(ma)=\chi(m)\,f(a)$, for
all $a\in \Z_p$, $m\in \Z_p^*$, using \eqref{eigen} for the
restriction of $\chi$ to $\Z_p^*$. One has
$$
\langle D'_0(s),\tilde f_0\rangle=\,\prod_S\, \langle D'_p(s),\tilde
f_p\rangle .
$$
One can assume that $S$ contains the finite set $P$ of places at
which $\chi$ is ramified, \ie where the restriction of $\chi$ to $\Z_p^*$
is non-trivial. For $p\in P$ one gets that $f_p$ vanishes in a
neighborhood of $0\in \Z_p$, while
$$
\langle D'_p(s),\tilde
f_p\rangle=\,\int_{\Q_p^*}\,f_p(x)\,\chi_0(x)\,|x|^s\,d^*x ,
$$
which, as a function of $s$, is a finite sum of exponentials $p^{ks}$.

For $p\notin P$, one has the normalization condition
$$
\langle D'_p(s), 1_{\Z_p}\rangle=\,1
$$
and the function $f_p$ is a finite linear combination of the
functions $x\to 1_{\Z_p}(p^k\,x)$. It follows from homogeneity of
the distribution $D'_p(s)$ that $\langle D'_p(s),\tilde f_p\rangle$,
as a function of $s$, is a finite sum of exponentials $p^{ks}$.

Similarly, if $\chi_0$ is non-trivial, the function $L(\chi_0 , -i t)$
belongs to $\Hol(I_\beta)$ and $\zeta(-i t)=L(1 , -i t)$ is the sum
of $1/(-i t-1)$ with an element of $\Hol(I_\beta)$. The presence of
the pole at $t=i$ for $\chi_0=1$ is taken care of by the condition
\eqref{intzero}. All of this shows that
\begin{equation}\label{hath}
\tilde F h \in \cS(I_\beta)\,, \quad \tilde F h (z)=0 \ \
\ \forall z\in I_\beta\,, \ \ \text{ with } \ L(\chi_0 , -i z)=0 .
\end{equation}
This shows in particular that $h\in \cS(\R_+^*)$, which gives the
required estimate \eqref{est0} on $|h(\lambda)|$, for $\lambda \to
0$. This ends the proof of Lemma \ref{bcdecay}.

\medskip

{\bf Proof of Theorem \ref{specreal}}.

It follows from \eqref{hath} that all elements $h$ of the range of
$\delta$ in degree zero have Fourier transform  $\tilde F h
\in\cS(I_\beta)$ that vanishes at any $z$ such that $L(\chi_0 , -i
z)=0$. Moreover, explicit choices of the test function  show that, if we
let $L_\Q$ be the complete $L$-function  (with its Archimedean Euler
factors), then, for $\chi_0$ non-trivial,
$L_\Q(\chi_0 , -i z)$ is the Fourier transform of an
element of the range of $\delta$. When
$\chi_0$ is trivial, the range of $\delta$ contains
$\xi(z-\frac{i}{2})$, where $\xi$ is the Riemann $\xi$ function
\cite{[R]}
$$
\xi (t):=-\frac{1}{2}\left(\frac{1}{4} +  t^2\right)\Gamma\left(\frac{1}{4} +
\frac{i t}{2}\right)\pi^{-\frac{1}{4} - \frac{i t}{2}} \zeta\left(\frac{1}{2} +
i t\right) .
$$
Thus, we see that, for $h$ in the range of $\delta$ in degree zero,
the set $Z_{\chi_0}$ of common zeros of the $\tilde F h$ is given by
$$Z_{\chi_0} =\,\{z\in I_\beta\;|\;L_\Q(\chi_0
, -i z)=0\}\,.$$

Using Lemma \ref{covarpi} one gets the required spectral
realization.

Indeed, this lemma shows that the $\Hol(I_\beta)$-module $\cM$ is
the quotient of $\cS(I_\beta)$ by the $\Hol(I_\beta)$-submodule
$\cN=\,\{\tilde F h\,|\,h\in \delta(\hat \cA_\beta)\}$. Let $z\in
I_\beta$ with $L_\Q(\chi_0 , -i z)=0$. Then the map $\cS(I_\beta) \ni
f\mapsto f(z)$ vanishes on $\cN$ and yields a non-trivial map
of $\Mc\otimes_{\Hol(I_\beta)}\C_z$ to $\C$. Let $z\in I_\beta$ with
$L_\Q(\chi_0 , -i z)\neq 0$. Any element of
$\Mc\otimes_{\Hol(I_\beta)}\C_z$ is of the form $f\otimes 1$, with
$f\in \cS(I_\beta)$. There exists an element $\xi= g\otimes 1\in
\cN$ such that $g(z)\neq 0$. Thus, we can assume, since we work
modulo $\cN$, that $f(z)=0$. It remains to show that $f\otimes 1$ is in
$\cS(I_\beta)\otimes_{\Hol(I_\beta)}\C_z$. To see this, one writes
$f(t)=\,(t-z)\,f_1(t)$ and checks that $f_1\in \cS(I_\beta)$, while
$$
f\otimes 1=\,(t-z)\,f_1\otimes 1\sim f_1\otimes (t-z)=\,0\,.
$$

\medskip

\subsection{Type III factors and unramified extensions}\label{IJ}

\medskip

In this section we briefly recall  the classification of type III
factors (\cf \cite{Co-th}) and we explain the analogy between the reduction
from type III to type II and automorphisms, first done in
\cite{Co-th}, and the use by Weil of the unramified extensions
$K\subset \,K\otimes_{\F_q}\,\F_{q^n}$ of global fields of positive
characteristic.

Central simple algebras over a (local) field $\L$ naturally arise  as
the commutant of isotypic representations of semi-simple algebraic
structures over $\L$. When one works over $\L=\C$ it is natural to
consider unitary representations in Hilbert space and to restrict to
the algebras $M$ (called von Neumann algebras) which appear as
commutants of unitary representations. The central von Neumann
algebras are called {\em factors}.

The module of
a factor $M$ was first defined in \cite{Co-th} as a closed subgroup of
$\R_+^*$ by the equality
\begin{equation}\label{modfac1}
S(M) = \bigcap_{\varphi} \ {\rm Spec} (\Delta_{\varphi}) \subset
\R_+ ,
\end{equation}
where $\varphi$ varies among (faithful, normal) states on $M$, while
the operator $\Delta_{\varphi}$ is the  modular operator
\eqref{modop}.

The central result of \cite{Co-th}, which made it possible to perform
the classification of factors, is that the 1-parameter family
$$
\sigma_t^{\varphi}\in \Out(M)=\Aut(M)/\Inn(M)
$$
is {\em independent} of the choice of the state $\varphi$.
Here $\Inn(M)$ is
the group of inner automorphisms of $M$. This gives a canonical one
parameter group of automorphism classes
\begin{equation}
\label{timefac} \delta \;:\; \R \to \Out(M)\,.
\end{equation}

This implies, in particular, that the crossed product $\hat
M=\,M\rtimes_{\sigma_t^{\varphi}}\R$ is {\em independent} of the
choice of $\varphi$ and one can show (\cf \cite{ct}) that the dual
action $\theta_\lambda$ is also independent of the choice of
$\varphi$.

The classification of {\it approximately finite dimensional} factors
involves the following steps.
\begin{itemize}
\item  The definition of the invariant ${\rm Mod} (M)$ for arbitrary
factors (central von Neumann algebras).
\item  The equivalence of all possible notions of approximate finite
dimensionality.
\item  The proof that Mod is a complete invariant\footnote{
we exclude the trivial case $M = M_n (\C)$ of matrix algebras, also
${\rm Mod} (M)$ does not distinguish between type II$_1$ and
II$_\infty$} and that all virtual subgroups are obtained.
\end{itemize}

In the general case,
\begin{equation}\label{modfac}
{\rm Mod} (M) \mathop{\subset}_{\sim} \ \R_+^* \, .
\end{equation}
was defined in \cite{ct} (see also \cite{Co-th})
and is a virtual subgroup of $\R_+^*$ in the sense of G.~Mackey,
i.e. an ergodic action of $\R_+^*$. All ergodic flows appear and
$M_1$ is isomorphic to $M_2$ iff ${\rm Mod} (M_1) \cong {\rm Mod}
(M_2)$.

The period group $T(M)=\,\Ker\,\delta$ of a factor was defined in
\cite{Co-th} and it was shown that, for $T\in T(M)$, there exists a
faithful normal state on $M$ whose modular automorphism is periodic with
$\sigma_T^\varphi =1$. The relevance in the context described above
should now be clear since for such a {\em periodic state}
the discussion above regarding
the distilled module simplifies notably. In fact, the multiplicative
group $\R_+^*$ gets replaced everywhere by its discrete subgroup
$$
\lambda^\Z\subset \R_+^*\,,\quad \lambda=\,e^{-2\pi/T} .
$$

This situation happens in the context of the next section, for
function fields, where one gets $\lambda=\,\frac{1}{q}$ in terms of
the cardinality of the field of constants $\F_q$. In  that situation
the obtained factor is of type III$_\frac{1}{q}$.

It is an open
question to find relevant number theoretic examples dealing with
factors of type III$_0$.

We can now state a
proposition showing that taking the cross product by the modular
automorphism is an operation entirely analogous to the unramified
extension
$$
\K\subset \,\K\otimes_{\F_q}\,\F_{q^n}
$$
of a global field $\K$ of positive characteristic with field of
constant $\F_q$.

\medskip

\begin{Pro} Let $M$ be a factor of type III$_1$.

 1) For $\lambda\in ]0,1[$ and $T=\,-\frac{2\pi}{\ln \lambda}$, the
algebra
$$ M_T=\,M\rtimes_{ \sigma_T}\,\Z $$
is a factor of type III$_\lambda$ with
$$
{\rm Mod} (M_T)=\,\lambda^\Z \subset \R_+^* .
$$

2) One has a natural inclusion $M\subset M_T $ and  the dual action
$$ \theta\;:\;\R_+^*/\,\lambda^\Z\to {\rm Aut}\,M_T $$
of $U(1)=\,\R_+^*/\,\lambda^\Z$ admits $M$ as fixed points,
$$
 M=\,\{x\in
M_T\;;\;\theta_\mu(x)=\,x \, ,\quad \forall \mu\in
\R_+^*/\,\lambda^\Z\} .
$$

\end{Pro}

\medskip

In our case, the value of $\lambda$ is given, in the analogy with the
function field case, by $\lambda=\,\frac{1}{p}$. Thus, we get
$$
T=\,\frac{2\pi}{\ln p}\,,\quad {\rm Mod} (M_T)=\,p^\Z \subset \R_+^*
$$
 We can
summarize some aspects of the analogy between global fields and
factors in the set-up decsribed above in the following table.

\bigskip

\begin{center}
\begin{tabular}{|c|c|}
\hline &  \\
 \bf{Global field $\K$} &\ \bf{Factor $M$} \\
&\\ \hline &  \\
${\rm Mod}\,\K\subset \R_+^*$  &\ ${\rm Mod}\,M\subset \R_+^*$  \\
&\\
\hline &  \\
 $\K\to \,\K\otimes_{\F_q}\,\F_{q^n}
$ &\ $M\to \,M\rtimes_{ \sigma_T}\,\Z$ \\
&\\ \hline &  \\
 $\K\to \,\K\otimes_{\F_q}\,\bar \F_{q}
$ &\ $M\to \,M\rtimes_\sigma \R$ \\
& \\ \hline
\end{tabular}
\end{center}

\bigskip
\bigskip

Notice also the following. For a global field $\K$ of positive
characteristic, the union of the fields $\K\otimes_{\F_q}\,\F_{q^n}$,
\ie the field $\K\otimes_{\F_q}\,\bar \F_{q}$, is the field of
rational functions from the associated curve $C$ to $\bar \F_{q}$.
In a similar manner, the crossed product $ M \rtimes_\sigma\,\R $
corresponds, at the geometric level, to the space of commensurability
classes of one dimensional $\Q$-lattices (\cf \cite{CM}) closely
related  to the noncommutative adele class space $X_\Q$ to which we
now turn our attention.

It is of course crucial to follow the above analogy in the specific
example of the BC-system. We shall show in \S \ref{classical} that
there is a natural way to associate to a noncommutative space $X$
not only its set of ``classical points" but in a more subtle manner
the set of all its classical points that can be defined over an
unramified extension, through the following basic steps:
$$
X \stackrel{Dual\, System} \longrightarrow \hat X
\stackrel{Periodic\, Orbits}\longrightarrow \cup\, {\hat X}_v
\stackrel{Classical\, Points}\longrightarrow \cup\,\Xi_v=\Xi
$$
When applied to the BC-system this procedure yields a candidate for
an analogue in characteristic zero of the geometric points $C(\bar
\F_q)$ of the curve over $\bar \F_q$ in the function field case.

\bigskip
\section{Geometry of the adele class space}\label{Aclass}

We apply the cohomological construction of the previous section
in the general
context of global fields and we obtain a cohomological interpretation
of the spectral realization of zeros of Hecke $L$-functions of
\cite{AC}. This also provides a geometric framework in which one can
begin to transpose the main steps of Weil's proof of RH in the case
of positive characteristic to the case of number fields.

We give in
this section a brief summary of these results, which will be dealt
with in details in our forthcoming paper \cite{CoCM}.

\subsection{The spectral realization on ${\mathcal H}^1_{\K,\C}$}

We let  $\Ac=\, \Sc(\A_\K) \rtimes \K^*$ be the noncommutative
algebra of coordinates on the adele class space $X_\K$. It is the
crossed product of the Bruhat-Schwartz algebra $\Sc(\A_\K)$ by the
action of $\K^*$. We let  $\cG_\K=\K^*\ltimes \A_\K$ be the locally
compact etale groupoid associated to the action of $\K^*$ on $\A_\K$
and write $\Ac=\, \Sc(\cG_\K)$. We then consider the cyclic module
\begin{equation}\label{cyc}
\Ac^\natural_0=\, \cap \;{\rm Ker}\;\varepsilon_j^\natural
\end{equation}
where the two cyclic morphisms
$$
\varepsilon_j^\natural \;:\;\Ac^\natural \to \C, \ \ \ \text{ for } \
j=0,1,
$$
are given by the traces on $\Ac$ associated, respectively,
to the evaluation at
$0\in \A_\K$ for $\varepsilon_0$ and to the integral on $\A_\K$ for
$\varepsilon_1$. The integral is taken with respect to
the additive Haar measure, which is $\K^*$-invariant and hence gives
a trace on $\Ac=\, \Sc(\A_\K) \rtimes \K^*$.  In
other words, one has
\begin{equation}\label{eps}
\varepsilon_0(\sum\,f_g\,U_g)=\,f_1(0)\,,\ \ \ \ \
\varepsilon_1(\sum\,f_g\,U_g)=\,\int\,f_1(a)\,da ,
\end{equation}
and in higher degree
$$
\varepsilon_j^\natural(x^0\otimes \cdots \otimes x^n)=\,
\varepsilon_j (x^0\, \cdots \, x^n) .
$$

\medskip

In this context, as we saw in the proof of Theorem \ref{specreal},
the cooling morphism of Section \ref{Frobenius} is simply given by
the restriction map to the idele class group $C_\K\subset X_\K$. In
other words, it is given by the restriction map
\begin{equation}\label{rho0restr}
\rho(\sum\,f_g\,U_g)=\,\sum\,f_g|_{\A_\K^*}\;\,U_g \,,
\end{equation}
where $f|_{\A_\K^*}\in C(\A_\K^*)$ is the restriction of $f\in
\Sc(\A_\K)$ to the ideles $\A_\K^*\subset \A_\K$. The action of
$\K^*$ on $\A_\K^*$ is free and proper and the corresponding
groupoid $\K^*\ltimes \A_\K^*$ is Morita equivalent to $C_\K$. We
use the exact sequence of locally compact groups
\begin{equation}\label{akbykstar}
1\to \K^*\to \A_\K^*\stackrel{p}\to C_\K\to 1
\end{equation}
to parameterize the orbits of $\K^*$ as the fibers $p^{-1}(x)$ for
$x\in C_\K$. By construction the Hilbert spaces
\begin{equation}\label{hilbfiber}
\cH_x=\ell^2(p^{-1}(x)) \qqq x\in C_K
\end{equation}
form a continuous field of Hilbert spaces over $C_\K$. We let
$\cL^1(\cH_x)$ be the Banach algebra of trace class operators in
$\cH_x$, these form a continuous field over $C_\K$.

\medskip

\begin{Pro}\label{restrAKCK1}
The restriction map $\rho$ of \eqref{rho0restr} extends to an
algebra homomorphism
\begin{equation}\label{rhoAKCKrestrtr}
\rho: \cS(\cG_\K) \to C(C_\K,\cL^1(\cH_x))\,.
\end{equation}
\end{Pro}

\proof Each $p^{-1}(x)$ is globally invariant under the action of
$\K^*$ so the crossed product rules in $C_\rho(\A^*_\K)\rtimes \K^*$
are just   multiplication of operators in $\cH_x$. To show that the
obtained operators are in $\cL^1$ we just need to consider monomials
$f_k\,U_k$. In that case the only non-zero matrix elements
correspond to $k=x y^{-1}$. It is enough to show that, for any $f\in
\cS(\A_\K)$, the function $k\mapsto f(k\,b)$ is summable. This
follows from the discreteness of $b\,\K\subset \A_\K$ and the
construction of the Bruhat--Schwartz space $\cS(\A_\K)$, \cf
\cite{Bruhat}, \cite{AC}. In fact the associated operator is of
finite rank when $f$ has compact support. In general what happens is
that the sum will look like the sum over $\Z$ of the values $f(nb)$
of a Schwartz function $f$ on $\R$.
\endproof

We let, in the number field case,
\begin{equation}\label{SCK}
{\bf S\,}(C_\K)=\,\cap_{\beta \in \R}\,\mu^\beta \cS(C_\K),
\end{equation}
where $\mu \in C(C_\K)$ is the module morphism from $C_\K$ to
$\R_+^*$. In the function field case one can simply use for ${\bf
S\,}(C_\K)$ the Schwartz functions with compact support.

\medskip

\begin{Def}\label{strongS}
We define ${\bf S\,}^\natural(C_\K, \cL^1(\cH_x))$ to be the cyclic
submodule of the cyclic module $C(C_\K, \cL^1(\cH_x))^\natural$,
whose elements are continuous functions such that the trace of the
restriction to the main diagonal belongs to ${\bf S\,}(C_\K)$.
\end{Def}

Note that for $T\in C(C_\K, \cL^1(\cH_x))^\natural$ of degree $n$,
$T(x_0,\ldots,x_n)$ is an operator in $\cH_{x_0}\otimes \ldots
\otimes \cH_{x_n}$. On the diagonal, $x_j=x$ for all $j$, the trace
map corresponding to $\Tr^\natural$ is given by
\begin{equation}\label{Trnatural}
\Tr^\natural( T_0\otimes T_1\otimes\ldots\otimes T_n)=\Tr( T_0\,
T_1\,\ldots\, T_n)\,.
\end{equation}
This makes sense since on the diagonal all the Hilbert spaces
$\cH_{x_j}$ are the same.

\begin{Def}\label{motH1AKCK}
We define ${\mathcal H}^1_{\K,\C}$ as the cokernel of the cyclic
morphism $$ \rho^\natural \ :\ \cS(\cG_\K)^\natural_0 \to  {\bf
S\,}^\natural(C_\K, \cL^1(\cH_x))
$$
\end{Def}

We can use the result of \cite{CoExt}, describing the cyclic
(co)homology in terms of derived functors in the category of cylic
modules, to write the cyclic homology as
\begin{equation}\label{HCnTorC}
HC_n(\cA) = \Tor_n (\C^\natural,\cA^\natural) .
\end{equation}
Thus, we obtain a cohomological realization of the  cyclic module
$\cH^1(\A_\K/\K^*,C_\K)$ by setting
\begin{equation}\label{TorH1mot}
H^1(\A_\K/\K^*,C_\K):= \Tor(\C^\natural,\cH^1(\A_\K/\K^*,C_\K)).
\end{equation}
We think of this as an $H^1$ because of its role as a relative term
in a cohomology exact sequence of the pair $(\A_\K/\K^*,C_\K)$.

\medskip

We now show that   $H^1(\A_\K/\K^*,C_\K)$ carries an action of
$C_\K$, which we can view as the abelianization $W_\K^{ab}\sim C_\K$
of the Weil group. This action is induced by the multiplicative
action of $C_\K$ on $\A_\K/\K^*$ and on itself. This generalizes to
global fields the action of $C_\Q=\hat\Z^*\times \R^*_+$ on
$HC_0(D(\cA,\varphi))$ for the Bost--Connes endomotive.

\begin{Pro}\label{CKactionH1}
The cyclic modules $\cS(\cG_\K)^\natural_0$ and ${\bf
S\,}^\natural(C_\K, \cL^1(\cH_x))$ are endowed with an action of
$\A^*_\K$ and the morphism $\rho^\natural$ is $\A^*_\K$-equivariant.
This induces an action of $C_\K$ on $H^1(\A_\K/\K^*,C_\K)$.
\end{Pro}

\proof For $\gamma\in \A^*_\K$ one defines an action by
automorphisms of the algebra $\cA=\cS(\cG_\K)$ by setting
\begin{equation}\label{actCKAKbis}
\vartheta_a(\gamma)(f)(x):= f(\gamma^{-1} x), \ \ \text{ for } f\in
\cS(\A_\K),
\end{equation}
\begin{equation}\label{actCKAK}
\vartheta_a(\gamma)(\sum_{k\in \K^*} f_k\,U_k):=\sum_{k\in \K^*}
\vartheta_a(\gamma)(f_k)\,U_k\,.
\end{equation}
This action is inner for $\gamma\in \K^*$ and induces an outer
action
\begin{equation}\label{actCKAKout}
C_\K\to \Out(\cS(\cG_\K))\,.
\end{equation}
 Similarly, the continuous field $\cH_x=\ell^2(p^{-1}(x))$  over
 $C_\K$
 is $\A^*_\K$-equivariant for the action of $\A^*_\K$ on $C_\K$ by
 translations, and the equality
 \begin{equation}\label{actAKltwo}
(V(\gamma)\xi)(y):=\xi(\gamma^{-1}\,y) \qqq y\in p^{-1}(\gamma
x)\,,\ \xi \in \ell^2(p^{-1}(x))\,,
\end{equation}
defines an isomorphism $\cH_x\stackrel{V(\gamma)}\longrightarrow
\cH_{\gamma x}$. One obtains then an action of $\A^*_\K$ on $C(C_\K,
\cL^1(\cH_x))$ by setting
\begin{equation}\label{actCKCK}
\vartheta_m(\gamma)(f)(x):=V(\gamma)\,f(\gamma^{-1}\,x)\,V(\gamma^{-1}),
\ \ \ \forall f\in C(C_\K, \cL^1(\cH_x))\,.
\end{equation}
The morphism $\rho$ is $\A^*_\K$-equivariant, so that one obtains an
induced action on the cokernel $\cH^1(\A_\K/\K^*,C_\K)$. This action
is inner for $\gamma\in \K^*$ and thus induces an action of $C_\K$
on $H^1(\A_\K/\K^*,C_\K)$.
\endproof

We  denote by
\begin{equation}\label{psigactH1}
C_\K \ni \gamma\mapsto {\underline\vartheta_m}(\gamma)
\end{equation}
the induced action on   $H^1(\A_\K/\K^*,C_\K)$.

\medskip

\begin{The}
The representation of $C_\K$ on
$$
H^1_{\K,\C_\chi}=HC_0({\mathcal H}^1_{\K,\C_\chi})
$$
gives the spectral realization of the zeros of the $L$-function with
Gr\"ossencharakter $\chi$.
\end{The}

This is a variant of Theorem 1 of \cite{AC}. We use the same
notation as in \cite{AC} and, in particular, we take a non canonical
isomorphism
$$
C_\K \sim \,C_{\K,1}\,\times N ,
$$
where $N\subset \R_+^*$ is the range of the module. For number
fields one has $N=\R_+^*$, while for fields of nonzero
characteristic $N \sim \Z$ is the subgroup $q^{\Z} \subset \R_+^*$,
where $q=p^{\ell}$ is the cardinality of the field of constants.

Given a character $\chi$ of $C_{\K,1}$, we let $\wt{\chi}$ be the
unique extension of $\chi$ to $C_\K$ which is equal to $1$ on $N$.

The map $E$ used in \cite{AC} involved a slightly different
normalization since the result of the summation on $\K^*$ was then
multiplied by $\vert x\vert^{1/2}$. Thus in essence the
representation $W$ given in \cite{AC} is related to
\eqref{psigactH1} above by
$$
{\underline\vartheta_m}(\gamma)=\,\vert
\gamma\vert^{1/2}\,W(\gamma).
$$
Instead of working at the Hilbert space level we deal with
function spaces dictated by the algebras at hand. This implies in
particular that we no longer have the restriction of unitarity
imposed in \cite{AC} and we get the spectral side of the trace formula
in the following form.

\begin{The} \label{spectral}  For any  function $f
\in {\bf S\,}(C_\K)$, the operator $${\underline\vartheta_m}(f) =
\int {\underline\vartheta_m}(\gamma) \, f(\gamma) \, d^* \, \gamma
$$ acting on $H^1_{\K,\C_\chi}$ is of trace class, and its trace is
given by
$$
\Trace \; \left({\underline\vartheta_m}(f) |_{ H^1} \right)=
\sum_{{L\left(\wt{\chi} , \rho \right) =0 \,,\;\rho \in {\C / N^{
\bot}}}} \wh f (\wt{\chi} ,\rho) ,
$$
where the Fourier transform $\wh f$ of $f$ is defined by
$$
\wh f (\wt{\chi} ,\rho) = \int_{C_\K} f(u) \, \wt{\chi} (u) \, \vert
u \vert^\rho \, d^* \, u \, . $$
\end{The}

\medskip
This is a variant of Corollary 2 of \cite{AC}. At the technical
level the work of Ralf Meyer \cite{Meyer} is quite relevant and
Theorem \ref{spectral} follows directly from \cite{Meyer}.

\medskip

\subsection{The Trace Formula}

A main result in \cite{AC} is a trace formula, called the $S$-local trace
formula, for the action of the class group on the simplified
noncommutative space obtained by considering only finitely many
places $v\in S$ of the global field $\K$. The remarkable feature of
that formula is that it produces exactly the complicated principal
values which enter in the Riemann-Weil explicit formulas of number
theory. It is then shown in \cite{AC} that, in the Hilbert space
context, the global trace formula is in fact equivalent to the
Riemann hypothesis for all $L$-functions with Gr\"ossencharakter.

It was clear from the start that relaxing from the Hilbert space
framework to the softer one of nuclear spaces would eliminate the
difficulty coming from the potential non-critical zeros, so that the
trace formula would reduce to the Riemann-Weil explicit formula.
However, it was not obvious how to obtain a direct geometric proof of
this formula from the $S$-local trace formula of \cite{AC}. This was
done in \cite{Meyer}, showing that the noncommutative geometry
framework makes it possible to give a geometric interpretation of
the Riemann-Weil explicit formula. While the spectral side of the
trace formula was given in Theorem \ref{spectral}, the geometric
side is given as follows.

\medskip

\begin{The} \label{geom} Let $h
\in {\bf S\,} (C_\K)$. Then the following holds:
\begin{equation}\label{trace3}
\Trace \; \left({\underline\vartheta_m}(h)|_{ H^1} \right)=\,\wh
h(0) + \wh h(1) -\,\Delta\bullet \Delta\;h(1) -\,\sum_{v}
\int_{(\K^*_v,e_{\K_v})}' \, \frac{h(u^{-1})}{
 |1-u|}\, d^* u \, .
\end{equation}
\end{The}

\medskip

We used the following notation. For a local field $\L$ one
chooses a preferred additive character $e_{\L}$ with the following
properties.

\begin{itemize}

\item $e_\R(x)=\,e^{-2\pi\,ix}\qqq x\in \R$

\item $e_\C(z)=\,e^{-2\pi\,i(z+\bar z)}\qqq z\in \C$

\item for $\L$ non-archimedian with maximal compact subring $\Oc$
the character $e_{\L}$ fulfills the condition ${\rm Ker}\,e_{\L}=\,\Oc$.

\end{itemize}

One lets $d^* \lambda $ be the multiplicative Haar measure
normalized by
$$
\int_{1 \leq \vert \lambda \vert \leq \Lambda} \, d^*\lambda \sim
\log \Lambda \qquad \hbox{when} \ \Lambda \to \infty \,
$$
For $\L=\K_v$, the definition of the finite value
$$ \int_{(\K^*_v,\beta)}' \,
\frac{h(u^{-1})}{
 |1-u|}\, d^* u
$$
relative to a given additive
character $\beta$ of $\L$ is obtained as follows \cite{AC}. Let
$\varrho_\beta$ be the unique distribution extending  $ d^*u$ at
$u=0$ whose Fourier transform relative to $\beta$ \ie
$$
\int \,\varrho(x)\,\beta(xy)dx=\,\hat\varrho(y).
$$
vanishes at 1, $\wh \varrho (1) = 0$. Notice that this does not
depend upon the normalization of the additive Haar measure $dx$ of
$\K_v$. One then has by definition
$$
\int'_{(\L,\beta)} \frac{h(u^{-1}) }{\vert 1-u \vert} \, d^* u
=\,\left\langle \varrho_\beta
 ,  \, g\right\rangle
\,, \quad g(\lambda) = h ((\lambda + 1)^{-1}) \, \vert \lambda + 1
\vert^{-1} \, .
$$
The slight shift of notation with \cite{AC} comes from the choice of
the additive characters $e_{\K_v}$ instead of characters
$\alpha_v(x)=\,e_{\K_v}(a_v\,x)$ such that $\K\subset \A_\K$ be
self-dual relative to $\prod\,\alpha_v$. The idele $a=(a_v)$ is
called a {\em differental idele}. This introduces the term
$\Delta\bullet \Delta$ given by
\begin{equation}\label{euler1}
\Delta\bullet \Delta=\log |a|
\end{equation}
This coincides, up to the overall factor $\log q$, with the Euler
characteristic in the function field case and with $-\log |D|\leq 0$
in the case of number fields, $D$ being the discriminant.

\subsection{The ``periodic classical points" of $X_\K $}\label{classical}

The origin (\cf \cite{AC}) of the terms in the geometric side of the
trace formula (Theorem \ref{geom}) comes from the Lefschetz formula
by Atiyah-Bott \cite{atbo} and its adaptation by Guillemin-Sternberg
(\cf \cite{gui}) to the distribution theoretic trace for flows on
manifolds, which is a variation on the theme of \cite{atbo}. For the
action of $C_\K$ on the adele class space $X_\K$ the relevant
periodic points are
\begin{equation} \label{globalper}
P=\{ (x,u)\in  X_\K \times C_\K \,|\,  u\,x=x\}
\end{equation}
and one has (\cf \cite{AC})
\begin{Pro}\label{preorbitckak} Let $(x,u)\in P$, with $u\neq 1$.
There exists a place $v\in \Sigma_\K$ such that
\begin{equation} \label{globalper1}
x\in X_{\K,v}=\{ x\in  X_\K \,|\,  x_v=0\}
\end{equation}
The isotropy subgroup of any $x\in X_{\K,v}$ contains the cocompact
subgroup
\begin{equation} \label{globalper2}
\K^*_v\subset C_\K\,,\  \  \ \K^*_v=\{ (k_w)\,|\, k_w=1\ \forall
w\neq v\}
\end{equation}
\end{Pro}
The spaces $X_{\K,v}$ are noncommutative spaces but they admit
 classical points  distilled (in the
sense of Section \ref{Frobenius}). One then obtains the following
description: For each place $v\in\Sigma_\K$, one lets $[v]$ be the
adele
\begin{equation}\label{base}
[v]_w=1\qqq w\neq v\,,\quad [v]_v=0\,.
\end{equation}

\smallskip

\begin{Def} \label{classsub} Let $\K$ be a global field and
$X_\K=\A_\K/ \K^*$ be the adele class space of $\K$. The set of {\em
periodic classical} points of $X_\K$ is
\begin{equation}\label{defcurve}
\Xi_\K=\,\cup \; C_\K\,[v]\,,\quad v\in \Sigma_\K .
\end{equation}
\end{Def}

\smallskip

The reason why the space $\Xi_{\K,\,v}=C_\K\,[v]$ inside the
periodic orbit $X_{\K,v}$ can be regarded as its set of classical
points  will be discussed more in detail in \cite{CoCM}. In essence,
one shows first that each $X_{\K,v}$ is a noncommutative space of
type III with a canonical time evolution $\sigma_t^v$. While points
$x\in X_{\K,v}$ yield irreducible representations $\pi_x$ of the
corresponding crossed product algebra, one shows that this
representation has positive energy (\ie the generator of
$\sigma_t^v$  has positive spectrum) if and only if $x\in
\Xi_{\K,\,v}\subset X_{\K,v}$.

\smallskip

While the structure of the space $\Xi_\K$ only comes from its being
embedded in the ambient noncommutative space $X_\K \smallsetminus
C_\K$, it is still worthwhile, in order to have a mental picture by
which to translate to standard geometric language, to describe
$\Xi_\K$ explicitly both in the function field case and in the case
$\K=\Q$.

The reduction to eigenspaces of
Gr\"ossencharakters corresponds to induced bundles over the quotient
$X_\K/C_{\K,1}$ so that we may as  well concentrate on the quotient
space $\Xi_\K/C_{\K,1}$ and see how it looks like under the
action of the quotient group $C_\K/C_{\K,1}$.

\smallskip

In the function field case, one has a {\em non-canonical}
isomorphism of the following form.

\begin{Pro} \label{Zspace}
Let $\K$ be the function  field of an algebraic curve $C$ over
$\F_q$. Then the  action of the Frobenius on $Y=C(\bar \F_q)$ is
isomorphic to the action of $q^\Z$ on the quotient
$$
\Xi_\K/C_{\K,1}.
$$
\end{Pro}

\begin{center}
\begin{figure}
\includegraphics[scale=1]{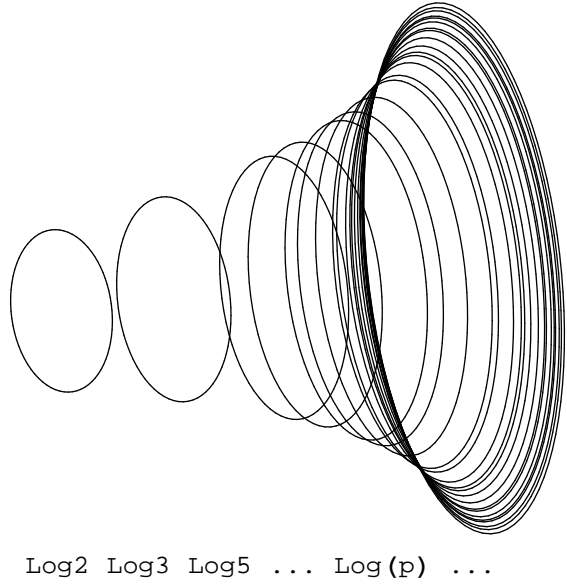}
\caption{The classical points of the adeles class space
\label{Figpoints}}
\end{figure}
\end{center}

In the case $\K=\Q$ the space $\Xi_\Q/C_{\Q,1}$ appears as the
union of periodic orbits of period $\log p$ under the action of
$C_\Q/C_{\Q,1}\sim \R$ (\cf Figure \ref{Figpoints}).

\subsection{Weil positivity}

It is a well known result of A. Weil (\cf \cite{weilpos}, \cite{EB})
that RH is equivalent to the positivity of the distribution entering
in the explicit formulae. This can be stated as follows.

\begin{The}  \label{pos} The following two conditions are
equivalent.
\begin{itemize}
\item All $L$-functions with Gr\"ossencharakter on $\K$ satisfy the Riemann
Hypothesis.
\item $\Trace\,{\underline\vartheta_m}(f\,\star\,
f^\sharp)|_{ H^1} \, \geq \,0$, for all $f\in  {\bf S\,} (C_\K)$.
\end{itemize}
\end{The}

Here we used the notation
\begin{equation}\label{conv}
f=\,f_1 \star f_2 \,,\ \ \ \text{ with } \ (f_1 \star f_2)(g)=\,\int
\,f_1(k)\, f_2(k^{-1}\,g) \,d^*g
\end{equation}
for the convolution of functions, using the multiplicative Haar
measure $d^*g$, and for the adjoint
\begin{equation}\label{adj1}
f\to f^{\sharp}\,,\quad f^{\sharp}(g)=\,\vert\, g\vert^{-1}\,\bar
f(g^{-1}) .
\end{equation}

The role of the specific correspondences used in Weil's proof of RH
in positive characteristic is played by the test functions $ f\in
{\bf S\,}(C_\K) $. More precisely the scaling map which replaces
$f(x)$ by $f(g^{-1} x)$ has a graph, namely the set of pairs $(x,
g^{-1} x) \in X_\K \times X_\K$, which we view as a correspondence
$Z_g$. Then, given a test function $f$ on the ideles classes, one
assigns to $f$ the linear combination
\begin{equation}\label{corresp}
Z(f)=\,\int f(g) Z_g d^*g
\end{equation}
of the above graphs, viewed as a ``divisor" on $X_\K \times X_\K$.

The analogs of the degrees $d(Z)$ and codegrees
$d^{\,'}(Z)=\,d(Z')$ of correspondences in the context of Weil's
proof are given, for the degree, by
\begin{equation}\label{degree}
d(Z(h))=\,\wh h(1)=\,\int\,h(u)\, \vert u \vert  \, d^* \, u ,
\end{equation}
so that the degree $d(Z_g)$ of the correspondence $Z_g$ is equal to
$|g|$. Similarly, for the codegree one has
\begin{equation}\label{codegree}
d^{\,'}(Z(h))=\,d(Z(\bar h^\sharp))=\,\int\,h(u)\,  d^* \, u =\,\wh
h(0) ,
\end{equation}
so that the codegree $d^{\,'}(Z_g)$ of the correspondence $Z_g$ is
equal to $1$.

\medskip

The role of principal divisors is to be found, but already one can
see that there is an interesting subspace $\cV$ of the linear space
of correspondences described above. It is given by the range of the
map $E$ or, more precisely, since there is a small shift in the
normalization, by the subspace
\begin{equation}\label{V}
\cV\subset  {\bf S\,}(C_\K) \,,\quad \cV=\{\,g(x)=\,
\sum\,\xi(k\,x)\,|\,\xi\in \Sc(\A_\K)_0\} ,
\end{equation}
where the subspace $\Sc(\A_\K)_0\subset \Sc(\A_\K)$ is defined by
the two boundary conditions
$$
\xi(0)=0 \,,\quad \int \,\xi(x)\,dx\,=\,0 .
$$

\medskip

\begin{Lem} \label{vanish}
For any $f\in \cV\subset {\bf S\,} (C_\K)$, one has
$$
{\underline\vartheta_m}(f)|_{ H^1}=\,0 .
$$
\end{Lem}

This shows that the Weil pairing of Theorem \ref{pos} admits a huge
radical given by all functions which extend to adeles. Thus, one is
led anyway to divide by this radical and hence to work with the
cohomology $H^1_{\K,\C}$ defined above. We will show in a
forthcoming paper \cite{CoCM} that several of the steps of Weil's
proof can be transposed in the framework described above. This
constitutes a clear motivation to develop noncommutative geometry
much further.

\bigskip
\section{Higher dimensional virtual correspondences}\label{Shigher}

In section \ref{Sartin} we have compared correspondences for
motives, given by algebraic cycles in the product $X\times Y$, and
correspondences for noncommutative spaces, given by bimodules
(elements in $KK$-theory), in the very special zero dimensional case
of Artin motives.

It is interesting to also consider higher dimensional cases. For
example, it appears from the recent work \cite{CMR}, \cite{CMR2},
\cite{HaPau} that the categories of motives of abelian varieties
should be enriched by adding degenerations of abelian varieties to
higher dimensional analogs of noncommutative tori, so that the
resulting moduli spaces are the noncommutative spaces considered in
\cite{CM}, \cite{CMR2} and \cite{HaPau}, having the classical
Shimura varieties as the set of classical points. Another compelling
reason for considering a unified setting for motives and
noncommutative spaces in higher dimensions is given by the results
on the Lefschetz formula for archimedean local factors of
$L$-functions of motives discussed in Section \ref{SLfactors} below.

Thus, it is useful to discuss the compatibility of morphisms given
by correspondences in the higher dimensional setting. To this
purpose, we compare correspondences defined by algebraic cycles with
the topological correspondences of \cite{Baum}, \cite{CoSka}, by
reformulating the case of algebraic cycles as a particular case of
the $KK$ correspondences of \cite{CoSka}.

\subsection{Geometric correspondences and $KK$-theory}

In the topological case, one considers smooth manifolds $X$ and $Y$
(in fact in \cite{CoSka} one only assumes $Y$ a smooth manifold
while $X$ is only a locally compact parameter space). A topological
(geometric) correspondence is given by data $(Z,E,f_X,g_Y)$ where
$Z$ is a smooth manifold with continuous maps $f_X:Z\to X$ and
$g_Y:Z\to Y$. One assumes $f_X$ continuous and proper and $g_Y$
continuous and $K$-oriented (orientation in $K$-homology).

The remaining piece of data $E$ is a complex vector bundle over $Z$.
Notice that in this setting one does not require that $Z\subset
X\times Y$, while, on the other hand, one has the extra piece of
data given by the vector bundle $E$. To the data $(Z,E,f_X,g_Y)$ one
associates an element $k(Z,E,f_X,g_Y)$ in $KK(X,Y)$ as in
\cite{CoSka}, defined as
\begin{equation}\label{kZEfg}
k(Z,E,f_X,g_Y) = (f_X)_* ( (E)\otimes (g_Y)! ),
\end{equation}
where $(E)$ is the class of $E$ in $KK(Z,Z)$. The element $(g_Y)!$
in $KK$-theory satisfies the following property.

Given $X_1$ and $X_2$ smooth manifolds and a continuous oriented map
$f: X_1 \to X_2$, the element $f! \in KK(X_1,X_2)$ gives the
Grothendieck Riemann--Roch formula
\begin{equation}\label{GRR}
{\rm ch}(F \otimes f!) =f_! ({\rm Td}(f) \cup {\rm ch}(F)),
\end{equation}
for all $F\in K^*(X_1)$, with ${\rm Td}(f)$ the Todd genus
\begin{equation}\label{Toddf}
{\rm Td}(f) = {\rm Td}(TX_1)/{\rm Td}(f^* TX_2).
\end{equation}

The composition of two correspondences $(Z_1,E_1,f_X,g_Y)$ and
$(Z_2,E_2,f_Y,g_W)$ is given in this setting by taking the fibered
product $Z=Z_1\times_Y Z_2$ and the bundle $E=\pi_1^*E_1 \times
\pi_2^* E_2$, with $\pi_i:Z\to Z_i$ the projections. One needs to
assume a transversality conditions on the maps $g_Y$ and $f_Y$ in
order to ensure that the fibered product is a smooth manifold.
Theorem 3.2 of \cite{CoSka} shows that the Kasparov product in
$KK$-theory is given by this product of correspondences, namely
\begin{equation}\label{prodGeomKK}
k(Z_1,E_1,f_X,g_Y)\circ k(Z_2,E_2,f_Y,g_W) =k(Z,E,f_X,g_W) \ \ \ \in
KK(X,W).
\end{equation}

\subsection{Cycles and $K$-theory}

In the algebro-geometric setting, the formulation of correspondences
that comes closest to the setting described above is obtained by
considering maps induced by algebraic cycles on $K$-theory (\cf
\cite{Man}).

Let us assume that $X$ and $Y$ are given smooth projective algebraic
varieties. Also we consider an algebraic cycle $Z=\sum_i n_i Z_i$ in
$X\times Y$,  
namely an element of the free abelian group generated by closed irreducible
algebraic subvarieties $Z_i$ of $X\times Y$. We may assume, for the purpose
of this discussion that $Z=Z_i$ is an irreducible subvariety of
$X\times Y$. We denote by
$p_X$ and $p_Y$ the projections of $X\times Y$ onto $X$ and $Y$,
respectively, and we assume that they are proper. We denote by
$f_X=p_X|_Z$ and $g_Y=p_Y|_Z$ the restrictions. Notice that, in the
case of complex algebraic varieties, the complex structure
in particular determines a $K$-orientation. 

To an irreducible subvariety $T\stackrel{i}{\hookrightarrow} Y$
we associate the coherent
$\cO_Y$-module $i_* \cO_T$. For simplicity of notation we write
it as $\cO_T$. We use
a similar notation for the coherent sheaf $\cO_Z$, for
$Z\hookrightarrow X\times Y$ as
above. The pullback
\begin{equation}\label{XYmod}
p_Y^* \cO_T =p_Y^{-1} \cO_T \otimes_{p_Y^{-1}\cO_Y} \cO_{X\times Y}.
\end{equation}
is naturally a $\cO_{X\times Y}$-module.

We consider the map on sheaves that corresponds to the cap product on
cocycles. This is given by
\begin{equation}\label{cyclemap2}
Z: \cO_T \mapsto (p_X)_* \left( p_Y^* \cO_T \otimes_{\cO_{X\times
Y}} \cO_Z  \right),
\end{equation}
where the result is a coherent sheaf, since $p_X$ is proper. Using
\eqref{XYmod}, we can write equivalently
\begin{equation}\label{cyclemap3}
Z: \cO_T \mapsto (p_X)_* \left( p_Y^{-1} \cO_T \otimes_{p_Y^{-1}
\cO_X} \cO_Z  \right).
\end{equation}

Recall that $f_!$ is right adjoint to $f^*$, \ie $f^*f_!=id$,
and that it satisfies the Grothendieck Riemann--Roch formula
\begin{equation}\label{GRR2}
{\rm ch}(f_!(F)) =f_! ({\rm Td}(f) \cup {\rm ch}(F)).
\end{equation}
Using this fact, we can equally compute the intersection product of
\eqref{cyclemap2} by first computing
\begin{equation}\label{shriekZ}
\cO_T \otimes_{\cO_Y} (p_Y)_! \cO_Z
\end{equation}
and then applying $p_Y^*$. Using \eqref{GRR2} and \eqref{GRR} we
know that we can replace \eqref{shriekZ} by $\cO_T \otimes_{\cO_X}
(\cO_Z \otimes (p_Y)!)$ with the same effect on $K$-theory.

Thus, to a correspondence in the sense of \eqref{cyclemap2} given by
an algebraic cycle $Z\subset X\times Y$ we associate the geometric
correspondence $(Z,E,f_X,g_Y)$ with $f_X=p_X|_Z$ and $g_Y=p_Y|Z$ and
with the bundle $E=\cO_Z$.

Now we consider the composition of correspondences. Suppose given
smooth projective varieties $X$, $Y$, and $W$, and (virtual)
correspondences $U=\sum a_i Z_i$ and $V=\sum c_j Z_j '$, with
$Z_i\subset X\times Y$ and $Z_j '\subset Y\times W$ closed reduced
irreducible subschemes. The composition of correspondences is then
given in terms of the intersection products
\begin{equation}\label{intprod}
U \circ V= (\pi_{13})_* ((\pi_{12})^* U \bullet (\pi_{23})^* V),
\end{equation}
with the projection maps $\pi_{12}:X\times Y\times W \to X\times Y$,
$\pi_{23}: X\times Y\times W \to Y\times W$, and $\pi_{13}: X\times
Y\times W \to X\times W$. Under an assumption analogous to the
transversality requirement for the topological case, we obtain the
following result.

\begin{Pro}\label{corrcorr2}
Suppose given smooth projective varieties $X$, $Y$, and $W$ and
correspondences $U$ given by a single $Z_1\subset X\times Y$ and $V$
given by a single $Z_2\subset Y\times W$. Assume that $(\pi_{12})^*
Z_1$ and $(\pi_{23})^* Z_2$ are in general position in $X\times
Y\times W$. Then assigning to a cycle $Z$ the topological
correspondence $\cF(Z)=(Z,E,f_X,g_Y)$ satisfies
\begin{equation}\label{prodprod}
\cF(Z_1\circ Z_2)=\cF(Z_1)\circ \cF(Z_2),
\end{equation}
where $Z_1\circ Z_2$ is the product of algebraic cycles and
$\cF(Z_1)\circ \cF(Z_2)$ is the Kasparov product of the topological
correspondences.
\end{Pro}

\proof Recall the following facts. Assume that $T_1,T_1\subset X$
are closed subschemes of $X$, with $\cO_{T_i}$ the corresponding
structure sheaves. If $T_1$ and $T_2$ are in general position, then
their product in $K_*(X)$ agrees with the intersection product
(\cite{Man}, Theorem 2.3)
\begin{equation}\label{prodKint}
[\cO_{T_1}]\cdot [\cO_{T_2}]= [\cO_{T_1\bullet T_2}].
\end{equation}

Using \eqref{prodKint} we write the composition \eqref{intprod} in
the form
\begin{equation}\label{intprodK}
Z_1\circ Z_2=(\pi_{13})_* ([\pi_{12}^* \cO_{Z_1}] \cdot [\pi_{23}^*
\cO_{Z_2}]).
\end{equation}
Notice that, for $Z_1\subset X\times Y$ and $Z_2\subset Y\times W$
the intersection in $X\times Y\times W$ is exactly the fiber product
considered in \cite{CoSka}, $\pi_{12}^{-1} Z_1 \cap \pi_{23}^{-1}
Z_2= Z=Z_1\times_Y Z_2$, with $\cO_Z=\pi_{12}^* \cO_{Z_1}\otimes
\pi_{23}^*\cO_{Z_2}$. Moreover, we have $\pi_{13}=(f_X,g_W):Z\to
X\times W$, where both $f_X$ and $g_W$ are proper, so that we can
identify $(f_X)_* ( \cO_Z \otimes (g_W)!)$ with $(\pi_{13})_*
(\pi_{12}^* \cO_{Z_1}\otimes \pi_{23}^* \cO_{Z_2})$ as desired. In
fact, Theorem 3.2 of \cite{CoSka} shows that $(f_X)_* ( \cO_Z
\otimes (g_W)!)$ represents the Kasparov product.
\endproof

Notice that, while in the topological (smooth) setting
transversality can always be achieved by a small deformation (\cf \S
III, \cite{CoSka}), in the algebro-geometric case one needs to
modify the above construction in the case of the product of cycles
that are not in general position. In this case the formula
\eqref{prodKint} is modified by Tor corrections and one obtains
(\cite{Man}, Theorem 2.7)
\begin{equation}\label{prodKint2}
[\cO_{T_1}]\cdot [\cO_{T_2}]= \sum_{i=0}^n (-1)^i \left[{\rm
Tor}_i^{\cO_X} (\cO_{T_1},\cO_{T_2})\right].
\end{equation}

\subsection{Algebraic versus topological $K$-theory}

We have considered so far only correspondences in $KK$-theory and
the cohomological realization given by an absolute
cohomology, defined for noncommutative spaces by the cyclic category
and cyclic cohomology. However, we found in the work \cite{CoCM}
that the ``primary'' invariants obtained by the Chern character map
from $K$-theory to cyclic cohomology is trivial, for the simple
reason that we are looking at an action of a continuous group (the
scaling action). In such cases, one wants to be able to consider a
refined setting in which secondary invariants appear and a more
refined version of $K$-theory is considered, which combines
algebraic and topological $K$-theory.

It is known from the results of \cite{CoKa} that, by viewing both
algebraic and topological $K$-theory in terms of homotopy groups of
corresponding classifying spaces, one obtains a fibration and a
corresponding long exact sequence. This sequence is related to the
long exact sequence for Hochschild and cyclic cohomology of
\cite{CoIHES} through the Chern character and a regulator map, so
that one has a commutative diagram with exact rows

\begin{center}
\diagram K^{top}_{n+1}(A) \rto\dto^{ch_{n+1}} & K^{rel}_n (A)
\rto\dto^{ch_n^{rel}}  & K_n^{alg} (A) \rto\dto^{D_n/n!} & K^{top}_n
(A) \rto\dto^{ch_n} & K^{rel}_{n-1}
(A) \dto^{ch_{n-1}^{rel}} \\
HC_{n+1}(A) \rto^S & HC_{n-1}(A)  \rto^{B/n} & HH_n(A) \rto^I &
HC_n(A) \rto^S & HC_{n-2}(A).
\enddiagram
\end{center}

\section{A Lefschetz formula for archimedean local
factors}\label{SLfactors}

As we mentioned already at the beginning of Section \ref{Sartin} we
can consider the $L$-functions $L(H^m(X),z)$ of a smooth projective
algebraic variety over a number field $\K$. This is written as an
Euler product over finite and archimedean places $v$ of $\K$,
\begin{equation}\label{Lfunction}
L(H^m(X),z) =\prod_v L_v (H^m (X),z).
\end{equation}
The local $L$-factors of \eqref{Lfunction} at finite primes encode
the action of the geometric Frobenius on the inertia invariants
$H^m(\bar X, \Q_\ell)^{I_v}$ of the \'etale cohomology in the form
(\cf \cite{Se})
\begin{equation}\label{L-factor}
L_v (H^m (X),z) = \det\left( 1-Fr_v^* N(v)^{-z} | H^m(\bar X,
\Q_\ell)^{I_v} \right)^{-1},
\end{equation}
with $N$ the norm map.

\smallskip

Serre showed in \cite{Se} that the local factors one needs to
consider at the archimedean primes, dictated by the expected form of
the functional equation for the $L$-function $L(H^m(X),z)$, depend
upon the Hodge structure
\begin{equation}\label{Hodge-str}
 H^m(X_v(\C))= \oplus_{p+q=m} H^{p,q}(X_v(\C)),
\end{equation}
where $v:\K \hookrightarrow \C$ is the archimedean place and
$X_v(\C)$ is the corresponding complex algebraic variety.

\smallskip

The archimedean local factors have an explicit formula (\cf
\cite{Se}) given in terms of Gamma functions, with shifts in the
argument and powers that depend on the Hodge numbers. Namely, one
has a product of Gamma functions according to the Hodge numbers
$h^{p,q}$ of the form
\begin{equation}\label{archL-factor}
L(H^*,z)= \left\{ \begin{array}{l} \prod_{p,q}
\Gamma_\C(z-\text{min}(p,q))^{h^{p,q}} \\[2mm]
\prod_{p<q}\Gamma_\C(z-p)^{h^{p,q}}\prod_p
\Gamma_\R(z-p)^{h^{p,+}}\Gamma_\R(z-p+1)^{h^{p,-}}
\end{array}\right.
\end{equation}
where the two cases correspond, respectively, to the complex and the
real embeddings, $\K_v=\C$ and $\K_v=\R$. Here $h^{p,\pm}$ is the
dimension of the $\pm(-1)^p$-eigenspace of the involution on
$H^{p,p}$ induced by the real structure and
\begin{equation}\label{GammaCR}
\Gamma_\C(z) := (2\pi)^{-z}\Gamma(z), \ \ \ \Gamma_\R(z)
:=2^{-1/2}\pi^{-z/2}\Gamma(z/2).
\end{equation}

\smallskip

It is clear that it is desirable to have a reformulation of the
formulae for the local factors in such a way that the archimedean
and the non-archimedean cases are treated as much as possible on
equal footing. To this purpose, Deninger in \cite{Den1}, \cite{Den2}
expressed both \eqref{L-factor} and \eqref{archL-factor} in the form
of zeta-regularized infinite determinants. Here we propose a
different approach. Namely, we reinterpret the archimedean local
factors as a Lefschetz Trace formula over a suitable geometric
space. The main idea is that we expect a global (or at least
semi-local) Lefschetz trace formula to exist for the $L$-functions
$L(H^m(X),z)$, over a noncommutative space which will be a
noncommutative generalization of the pure motive $h^m(X)$ and should
be obtained as an extension of $h^m(X)$ by a suitable modification
of the ad\`eles class space (by certain division algebras instead of
the local fields $\K_v$ at the places $v$ of bad reduction and at
the real archimedean place). Here we give some evidence for this
idea, by giving the local formula for a single archimedean place
$v$. One sees already in the case of a single real place that the
underlying geometric space for the trace formula is obtained by
passing to the division algebra of quaternions.

\smallskip

As we recalled already in Section \ref{Aclass} above, it was shown in
\cite{AC} that the noncommutative space of adele classes over a
global field $\K$ provides both a spectral realization of zeros of
$L$-functions and an interpretation of the explicit formulas of
Riemann-Weil as a Lefschetz formula. The corresponding trace formula
was proved in \cite{AC} in the semilocal case (finitely many
places). It is natural then to ask whether a similar approach can be
applied to the $L$-functions of motives, also in view of the fact
that the results of \cite{ConsMar}, \cite{ConsMar2} showed that
noncommutative geometry can be employed to describe properties of
the fibers at archimedean places (and at places of totally
degenerate reduction) of arithmetic varieties.

\smallskip

In this section we consider the archimedean factors of the
Hasse-Weil $L$-function attached to a non-singular projective
algebraic variety $X$ defined over a number field $\K$. We obtain
the real part of the logarithmic derivative of the archimedean
factors of the $L$-function $L(H^m, z)$ on the critical line as a
trace formula for the action of a suitable Weil group on a complex
manifold attached to the archimedean place.

In \S \ref{SWeilform} we make some preliminary calculations. First
we show (Lemma \ref{princval}) that for $\K_v=\C$ or $\R$, one
obtains the logarithmic derivative of the imaginary part of the
Gamma functions $\Gamma_{\K_v}$ \eqref{GammaCR} on the critical line
as the principal value on $\K_v^*$ of a distribution $$ \frac{\vert
u \vert^{ \frac{1}{2}+ i \,s}}{\vert 1-u \vert}. $$ We then show
(Lemma \ref{minpq}) that the shift by $\min\{ p,q \}$ in the
argument of $\Gamma_\C$ appears when one considers the principal
value on $\C^*$ of
$$ \frac{u^{-p}{\overline{u}}^{-q}\vert u \vert_\C^{ z}}
{\vert 1-u \vert_\C}. $$ This is sufficient to obtain in Theorem
\ref{complexpl} of \S \ref{SLefC} the logarithmic derivative of the
archimedean local factor as a trace formula for the action of $\C^*$
on a space with base $\C$ and fiber $H^m(X,\C)$ with the
representation
\begin{equation}\label{Gmhscomp}
 \pi(H^m,u)\,\xi = u^{-p} \overline{u}^{-q} \xi
 \end{equation}
of $\C^*$ on the cohomology. The real case is more delicate and it
is treated in \S \ref{SLefR}. In this case, to obtain a similar
result, one needs to consider a space that has base the quaternions
$\H$ and fiber the cohomology. We then obtain the logarithmic
derivative of the local factor as a trace formula for an action on
this space of the Weil group, with representation
\begin{equation}\label{Gmhsreal}
 \pi(H^m,w j)\,\xi =\,i^{p+q}\,w^{-p} \overline{w}^{-q}
F_\infty(\xi)
\end{equation}
 on the cohomology, with $F_\infty$ the linear
involution induced by complex conjugation. The properties of the
$L$-functions of algebraic varieties defined over a number field
$\K$ are mostly conjectural (\cite{Se}). In particular one expects
from the functional equation (\cite{Se}, C$_9$) that the zeros of
$L(H^m(X),z)$ are located on the critical line
$\Re(z)=\frac{1+m}{2}$. By analogy with the case of the Riemann zeta
function, one expects that the number of non-trivial zeros
\begin{equation}
N_s(E) =\# \{ \rho \, | \,  L(H^m(X),\rho)=0, \ \text{ and } \
-E\leq \Im(\rho)\leq E \}\,, \label{countingLzerossym}
\end{equation}
can be decomposed as the sum of an average part and an oscillatory
part
\begin{equation}
N_s(E) =\langle N_s(E) \rangle + N_{\rm osc} (E)\,,
\label{countingLzerosdec}
\end{equation}
where the average part is given in terms of the Archimedean local
factors described above
 \begin{equation}\label{countLzerosfinal}
\langle N_s(E) \rangle =\,\sum_{v|\infty} \, \frac{1}{\pi}
\int_{-E}^E\frac{d}{ds}\Im \log L_v(H^m(X),\frac{1+m}{2}+is) ds\, .
\end{equation}
We shall show in Theorems \ref{complexpl} and \ref{real} that the
key ingredient of this formula can be expressed in the form of a
Lefschetz contribution, using the above representations
\eqref{Gmhscomp} and \eqref{Gmhsreal} of the local Weil group $W_v$
\begin{equation}\label{Lefschetzloc}
\frac{1}{\pi}\,\frac{d}{ds}\Im \log L_v(H^m(X),\frac{1+m}{2}+is)=
-\frac{1}{2\pi}\, \int'_{W_v}\frac{\Trace(\pi_v(H^m,u))\vert u
\vert_{W_v}^{ z}}{ \vert 1-u \vert_{\H_v}} \, d^* u
\end{equation}
where $z=\frac{1+m}{2}+is$. The  local Weil group $W_v$ is  a
modulated group. It embeds in the algebra $\H_v$ which is the
algebra $\H$ of quaternions for $v$ a real place and the algebra
$\C$ for a complex place. We shall describe below, in more detail,
in \S \ref{SLefC} and \ref{SLefR} the natural representation
$\pi_v(H^m,u)$  of $W_v$.

In \S \ref{SQuest} we describe the problem of a semilocal trace
formula involving several archimedean places and explain how the
conjectured trace formula combines with \eqref{Lefschetzloc} to
suggest that the (non-trivial) zeros of the $L$-function appear as
an absorption spectrum as in \cite{AC}. In \S \ref{Scurve} we
address the problem of considering simultaneously archimedean and
nonarchimedean places and in particular the fact that in the usual
approaches one is forced to make a choice of embeddings of $\Q_\ell$
in $\C$. We suggest that, at least in the case of a curve $X$, one
may obtain a different approach by considering, as the fiber of the
space one would like to construct for a semi-local Lefschetz trace
formula, the adele class space of the function field of the curve
over the residue field $k_v$ at the place $v$.

\subsection{Weil form of logarithmic derivatives of local
factors}\label{SWeilform}

We begin by checking the formulae (\cite{weil}, \cite{AC}) that will
relate the Fourier transform of the local Lefschetz contribution
(viewed as a distribution on the multiplicative group) with the
derivative of the imaginary part of the logarithm of the archimedean
local factor $\Gamma_{\K_v}$, for $\K_v=\R$ ($v$ a real archimedean
place) or $\K_v=\C$ ($v$ a complex archimedean place).

We use the notations of \cite{AC} for the principal values. We
recall the following basic facts.

\begin{Lem}\label{princval}
For $\K_v=\R$ or $\K_v=\C$ and for $s$ real, one has
\begin{equation}\label{prval}
\int'_{\K_v^*}  \frac{\vert u \vert^{ \frac{1}{2}+ i \,s}}{\vert 1-u
\vert} d^* u =
 -2\, \frac{d}{ds}\Im \log \Gamma_{\K_v}\left( \frac{1}{2}+ i \,s\right) \, .
\end{equation}
\end{Lem}

\proof First notice that, for $s$ real, one has
\begin{equation}\label{imlogGamma}
2i\, \Im \log \Gamma_{\K_v}\left( \frac{1}{2}+ i \,s\right) =\log
\Gamma_{\K_v}\left( \frac{1}{2}+ i \,s\right) -\log
\Gamma_{\K_v}\left( \frac{1}{2}- i \,s\right),
\end{equation}
since $\Gamma_{\K_v}$ is a ``real" function, \ie it fulfills
$f(\overline z)=\overline{ f(z)}$.

Thus, one can rewrite the equality above in the form
\begin{equation}\label{prval2}
\int'_{\K_v^*}  \frac{\vert u \vert^{ \frac{1}{2}+ i \,s}}{\vert 1-u
\vert} d^* u = - \left(\frac{\Gamma'_{\K_v}}{\Gamma_{\K_v}} \left(
\frac{1}{2} + i \,s \right) + \frac{\Gamma'_{\K_v}}{\Gamma_{\K_v}}
\left( \frac{1}{2} - i \,s \right)\right) \, .
\end{equation}
The $i $ in $i\,s$ brings up the real part of
$\frac{\Gamma'_{\K_v}}{\Gamma_{\K_v}}$ which is again a ``real"
function.

Notice also that \eqref{prval2} now holds without the restriction
$s\in \R$, since both sides are analytic functions of $s\in \C$.

Let us first take $\K_v=\R$. In this case we have
\begin{equation}\label{GammaR}
\Gamma_\R (x) = 2^{-1/2} \pi^{-x/2}\Gamma(x/2),
\end{equation}
hence the equality takes the explicit form in terms of the usual
$\Gamma$-function
\begin{equation}\label{Gammacl}
\int'_{\R^*} \vert u \vert^{is}  \frac{\vert u \vert^{1/2}}{\vert
1-u \vert} d^* u =
  \log \pi - \frac{1}{2}\left(\frac{\Gamma'}{\Gamma} \left( \frac{1}{4} + i
\frac{s}{2} \right) + \frac{\Gamma'}{\Gamma} \left( \frac{1}{4} - i
\frac{s}{2} \right)\right) \, .
\end{equation}
which can be deduced from \cite{AC}, Appendix 2. Let us then take
$\K_v=\C$. In this case we have
\begin{equation}\label{GammaC}
\Gamma_\C (x) = (2\,\pi)^{-x} \Gamma(x),
\end{equation}
and again, in terms of the usual $\Gamma$-function, the equality
takes the form
\begin{equation}\label{Gammacl2}
\int'_{\C^*} \vert
 u \vert^{is}  \frac{\vert u \vert^{1/2}}{\vert 1-u \vert} d^* u =
  2\,\log 2\pi - \left(\frac{\Gamma'}{\Gamma} \left( \frac{1}{2} + i \,s
\right) + \frac{\Gamma'}{\Gamma} \left( \frac{1}{2} - i \,s
\right)\right) \, .
\end{equation}

This can be extracted from \cite{AC} (Appendix 2) and will be proved
as a special case of Lemma \ref{minpq} below.
\endproof

The next step consists of obtaining a similar formula for
\begin{equation}\label{derlogGamma}
\frac{d}{ds}\Im \log \Gamma_{\K_v}\left( \frac{1}{2}+ \frac{\vert n
\vert}{2}+i \,s\right) \, ,
\end{equation}
where $n \in \Z$ is an integer.

Notice that the function $\Gamma_{\K_v}( \frac{\vert n \vert}{2}+
z)$ is still real.

We take $\K_v =\C$ and we let
\begin{equation}\label{f0}
f_0 (\nu) = \inf (\nu^{1/2} , \nu^{-1/2}), \ \ \ \ \forall \nu \in
\R_+^* \,.
\end{equation}
We then obtain the following result.

\begin{Lem} \label{absn}
For $n\in\Z$ and $\rho\in \R_+^*$ with $\rho\neq 1$, one has
\begin{equation}\label{rho2eq}
\frac{1}{ 2\pi} \int_0^{2\pi} \frac{e^{i n\theta}}{ \vert 1 -
e^{i\theta} \rho\vert^2} \, d\theta = \frac{f_0 (\nu)^{\vert n
\vert}} {\vert 1 - \nu \vert} \,,\ \ \ \ \nu = \rho^2\,.
\end{equation}
\end{Lem}

\proof Let us first consider the case $\rho <1$. Then we have
\begin{equation}\label{rho2eq2}
\frac{1-\rho^2}{ \vert 1 - e^{i\theta} \rho\vert^2} =\frac{1}{1-\rho
e^{i\theta}}+\,\frac{1}{1-\rho e^{-i\theta}} -1  \, ,
\end{equation}
whose Fourier coefficients are  $\rho^{\vert n \vert}$. This gives
the equality in the case $\rho <1$, with $f_0 (\nu)=\rho$. One then
checks that both sides of \eqref{rho2eq} fulfill
$$
f(\rho^{-1})=\,\rho^2\,f(\rho)\,,
$$
which gives the desired equality for $\rho >1$.
\endproof

We then obtain the following result.

\begin{Lem} \label{minpq}
Let $\K_v=\C$ and suppose given $p,q \in \N$ with $m=p+q$. For
$z=\frac{1+m}{2}+\,i s$, with $s\in \R$, one obtains
\begin{equation}\label{intminpq}
\int'_{\C^*}  \frac{u^{-p}{\overline{u}}^{-q}\vert u \vert_\C^{ z}}
{\vert 1-u \vert_\C} d^* u =
 -2\, \frac{d}{ds}\Im \log \Gamma_\C\left(  z
-\min (p,q)\right) \, .
\end{equation}
\end{Lem}

\proof Let $n=p-q$. One has $\min (p,q)=\frac{m}{2}-\frac{\vert n
\vert}{2}$. Since $\vert u \vert_\C= u\overline{u}$, one has
$u^{-p}{\overline{u}}^{-q}= e^{-i n\theta}\,\vert u
\vert_\C^{\frac{-m}{2}}$, where $\theta$ is the argument of $u$. One
can then rewrite the desired equality in the form
\begin{equation}\label{intminpq2}
\int'_{\C^*}  \frac{ e^{-i n\theta}\,\vert u \vert_\C^{ \frac{1}{2}+
i \,s}} {\vert 1-u \vert_\C} d^* u =
 -2\, \frac{d}{ds}\Im\log\Gamma_\C \left( \frac{1}{2}+ i \,s
+\frac{\vert n \vert}{2}\right) \, .
\end{equation}

Recall now that $\int'$ coincides with the Weil principal value
$Pfw$ of the integral which is obtained as
\begin{equation}\label{Pfw}
Pfw \int_{\C^*} \varphi (u) \, d^* u = PF_0 \int_{\R_+^*} \psi (\nu)
\, d^* \nu \, ,
\end{equation}
where $\psi (\nu) = \int_{\vert u \vert_\C = \nu} \varphi (u) \,
d_{\nu} u$ is obtained by integrating $\varphi$ over the fibers,
while one has
\begin{equation}\label{PF0}
PF_0 \int \psi (\nu) d^* \nu = 2 \log (2\pi) \, c + \lim_{t \to
\infty} \left( \int (1-f_0^{2t} ) \, \psi (\nu) \, d^* \nu - 2c \log
t \right) \, .
\end{equation}
Here one assumes that $\psi - c \, \, f_1^{-1}$ is integrable on
$\R_+^*$, with $ f_1 = f_0^{-1} - f_0 $. In our case we have
$$\varphi (u)=\,\frac{ e^{-i n\theta}\,\vert u \vert_\C^{
\frac{1}{2}+ i \,s}} {\vert 1-u \vert_\C}.$$ By Lemma \ref{absn},
integration over the fibers gives
$$
\psi (\nu)=\,{\frac{\nu^{\frac{1}{2}+ i \,s}\,f_0 (\nu)^{\vert n
\vert}} {\vert 1 - \nu \vert}} \,.
$$
We have $c=1$ independently of $s$ and $n$.

One has the identity (\cite{AC}, Appendix 2, (46)),
$$
PF_0 \int f_0 \,f_1^{-1}   d^* \nu=  2 (\log 2\pi + \gamma) \, .
$$
This allows one to check the equality for $n=1$ and $s=0$. Indeed,
in that case one gets $\psi (\nu)=\,f_0 \,f_1^{-1}$, while
$$ -2\left(\frac{\Gamma'_\C}{\Gamma_\C}\left( \frac{1+n}{2} \right)
\right)= -2 (-\log(2\pi)-\gamma) $$

Having checked the result for some value while $c$ is independent of
both $n$ and $s$ one can now use any regularization to compare other
values of
\begin{equation}\label{finpart0}
\int_{\R_+^*} \,{\frac{\nu^{\frac{1}{2}+ i \,s}\,f_0 (\nu)^{\vert n
\vert}} {\vert 1 - \nu \vert}} \,
  d^* \nu \,.
\end{equation}
We write the integral as $\int_0^1+\int_1^\infty$ and use the
minimal substraction as regularization \ie the substraction of the
pole part in $\varepsilon$ after replacing the denominator $\vert 1
- \nu \vert$ by $\vert 1 - \nu \vert^{1-\varepsilon}$. The first
integral gives
$$
\int_0^1 \nu^{\left(\frac{1}{2} + i s+\frac{\vert n
\vert}{2}\right)} (1-\nu)^{-1 + \varepsilon} \frac{d\nu}{\nu} = B
\left( \frac{1}{2} + i s+\frac{\vert n \vert}{2}
 , \varepsilon \right) =
\frac{\Gamma \left( \frac{1}{2} + i s+\frac{\vert n \vert}{2}
\right) \Gamma (\varepsilon)}{\Gamma \left( \frac{1}{2} + i
s+\frac{\vert n \vert}{2} + \varepsilon \right)} \, .
$$

The residue at $\varepsilon = 0$ is equal to one and the finite part
gives
\begin{equation}\label{finpart}
- \frac{\Gamma'}{\Gamma} \left( \frac{1}{2} + i s+\frac{\vert n
\vert}{2} \right) - \gamma \, .
\end{equation}
The other integral $\int_1^\infty$ gives the complex conjugate,
\begin{equation}\label{finpartconj}
- \frac{\Gamma'}{\Gamma} \left( \frac{1}{2} - i s+\frac{\vert n
\vert}{2} \right) - \gamma \, .
\end{equation}

Thus, having checked the additive constant term, we get
\begin{equation}\label{intderimlog}
\int'_{\C^*}  \frac{ e^{-i n\theta}\,\vert u \vert_\C^{ \frac{1}{2}+
i \,s}} {\vert 1-u \vert_\C} d^* u =
 -2\, \frac{d}{ds}\Im\log\Gamma_\C \left( \frac{1}{2}+ i \,s
+\frac{\vert n \vert}{2}\right) \, ,
\end{equation}
and the required equality follows.
\endproof

\medskip
\subsection{Lefschetz formula for complex places}\label{SLefC}

We now look at a complex archimedean place with the local $L$-factor
$L_{\C}(H^m,z)$. We let $\pi(H^m,u)$  be the canonical
representation (\cite{De}) of $\C^* $ on $H^m=\oplus H^{p,q}$,
\begin{equation}\label{canreppq}
\pi(H^m,u)\,\xi = u^{-p} \overline{u}^{-q} \xi\,,\ \ \ \  \forall
\xi \in H^{p,q}\,, \ \ \forall u\in\C^*\,.
\end{equation}

We obtain the following formula for the archimedean local factor at
a complex place.

\begin{The} \label{complexpl}
Let $\K_v=\C$. For $m \in \N$, let $\pi(H^m,u)$ be the canonical
representation of $\C^*$ on $H^m$. Then, for $z=\frac{1+m}{2}+\,i s$
with $s\in \R$, we have
\begin{equation}\label{LfactLefC}
\int'_{\C^*} \frac{ \Trace (\pi(H^m,u))\vert u \vert_\C^{ z}} {\vert
1-u \vert_\C} d^* u =
 -2\, \frac{d}{ds}\Im\log L_{\C}(H^m,z) \, .
\end{equation}
\end{The}

\proof This follows directly from Lemma \ref{minpq} and the formula
\eqref{archL-factor} (\cf \cite{Se}) that expresses $L_{\C}(H^m,z)$
as a product of powers of the $\Gamma_\C (  z - \min(p,q)) $.

\medskip
\subsection{Lefschetz formula for real places}\label{SLefR}

We now look at a real archimedean place with the local $L$-factor
$L_\R (H^m,z)$. We let $W$ be the Weil group \cite{weil} which, in
this case, is the normalizer of $\C^*$ in $\H^*$ where $\H$ here
denotes the division algebra of quaternions $\H=\C + \C j$, with
$j^2=-1$ and $j\,w\,j^{-1}=\,\overline{w}$, $\forall w\in \C$.
Elements of $W$ are of the form
$$
u=\,w\;j^\epsilon\,,\quad w\in\C^*\,,\ \ \ \  \epsilon \in
\{0,1\}\,,
$$
where we use the notation $j^0=1$, $j^1=j$.

We let $\pi(H^m,u)$ be the canonical representation (\cite{De}) of
the Weil group in $H^m =\oplus H^{p,q}$. The subgroup $\C^*\subset
W$ acts as above, while elements in $\C^* j$ act by
\begin{equation}\label{actwj}
\pi(H^m,w j)\,\xi =\,i^{p+q}\,w^{-p} \overline{w}^{-q}
F_\infty(\xi)\,,\ \ \ \ \forall \xi \in H^{p,q}\, .
\end{equation}
Here $F_\infty$ is the linear involution associated to the geometric
action of complex conjugation (once translated in the cohomology
with complex coefficients) as in Serre (\cite{Se}).

One has to check that $\pi$ is a representation, and in particular
that $\pi(H^m, j)^2=\,\pi(H^m,-1)$. This follows for $\xi \in
H^{p,q}$ {}from $$\pi(H^m, j)^2
\xi=(-1)^{p+q}F_\infty^2(\xi)=(-1)^{p+q}\xi\,.$$

One checks in the same way that  $\pi(H^m, j) \,\pi(H^m,
w)=\,\pi(H^m, \overline{w}) \, \pi(H^m, j)$.

We now investigate the integral
\begin{equation}
\int'_{W} \frac{ \Trace (\pi(H^m,u))\vert u \vert_\H^{ z}} {\vert
1-u \vert_\H} d^* u\,,
\end{equation}
in which we follow the conventions of Weil \cite{weil}. Thus, for
$u\in W$, the module $\vert u \vert_\H$ is the same as the natural
module $\vert u \vert_W$ of the Weil group and $\vert 1-u \vert_\H$
is the reduced norm in quaternions.

Let us start by the case $m$ odd, since then the subtle term
involving $\Gamma_\R$ does not enter. Then the action of $W$ on
$H^{p,q}\oplus H^{q,p}$ will give the same trace as in the complex
case for elements $u=\,w\;j^0 =w\in \C^*$, but it will give a zero
trace to any element $u=\,w\;j^1=w\,j$. In fact, in this case the
action is given by an off diagonal matrix. This gives an overall
factor of $\frac{1}{2}$ in the above expression, since when one
integrates on the fibers of the module map $u\mapsto \vert u\vert$
from $W$ to $\R_+^*$ the only contribution will come from the
$u=\,w\;j^0=w$. This will give the same answer as for the complex
case, except for an overall factor $\frac{1}{2}$ due to the
normalization of the fiber measure as a probability measure.

Thus, things work when $m$ is odd, since then the local factor
$L_\R(H^m,z)$ is really the square root of what it would be when
viewed as complex.

When $m= 2 p$ is even and $h^{p+}= h^{p-}$, the same argument does
apply irrespectively of the detailed definition of the
representation of $W$. Indeed, elements of $W$ of the form
$u=\,w\;j$ have pairs of eigenvalues of opposite signs $\pm
w^{-p}\overline{w}^{-p}$. One thus gets the required result using
\eqref{archL-factor} and the duplication formula
\begin{equation}\label{dupl}
\Gamma_\R (  z) \,\Gamma_\R (  z+1)=\,\Gamma_\C (  z)  \, .
\end{equation}

When  $h^{p+}\neq h^{p-}$, we use the detailed definition of the
action of $j$. Notice that the computation will not be reducible to
the previous ones, since one uses the quaternions to evaluate
$$
\vert 1-u \vert_\H =\, 1 + \vert w \vert^2 \,,\quad u=w j\,,
$$

The action of $j$ is $(-1)^p $ times the geometric action of complex
conjugation (once translated in the cohomology with complex
coefficients) on the space $H^{p,p}$.

We compute in the following lemma the relevant integral.

\begin{Lem}\label{lemrelint}
For $z= \frac{1}{2}+i s$ with $s\in \R$, one has
\begin{equation}\label{trans}
\int_{\R_+^* } \frac{ u^{z}}{ 1+u } d^* u=\,- 2\frac{d}{ds}\Im
\log(\Gamma_\R (  z) /\Gamma_\R (  z+1)) \, .
\end{equation}
\end{Lem}

\proof First for $z\in \C$ with positive real part, one has
\begin{equation}\label{intm2p}
\int_0^1  \frac{ u^{z}}{ 1+u } d^* u = \frac{1}{2}
\left(\frac{\Gamma'}{\Gamma}\left(\frac{z+1}{2}\right)-
\frac{\Gamma'}{\Gamma}\left(\frac{z}{2}\right)\right)  \, .
\end{equation}

Let us prove \eqref{trans}. The factors in $\pi^{-z/2}$ in
\eqref{GammaR} for $\Gamma_\R(z) $ do not contribute to the right
hand side, which one can replace by
$$
-\frac{1}{2}\left(\frac{\Gamma'}{\Gamma} \left( \frac{1}{4} + i
\frac{s}{2} \right) + \frac{\Gamma'}{\Gamma} \left( \frac{1}{4} - i
\frac{s}{2} \right)\right)+\frac{1}{2}\left(\frac{\Gamma'}{\Gamma}
\left( \frac{3}{4} + i \frac{s}{2} \right) + \frac{\Gamma'}{\Gamma}
\left( \frac{3}{4} - i \frac{s}{2} \right)\right)\,,
$$

The equality follows from \eqref{intm2p}, using the symmetry of the
integral to integrate from $1$ to $\infty$. \endproof

We then obtain the following formula for the archimedean local
factor at a real place.

\begin{The} \label{real} Let $\K_v=\R$. For $m \in \N$, let
$\pi(H^m,u)$ be the representation of the Weil group $W $ on $H^m$
described above. Then, for $z=\frac{1+m}{2}+\,i s$ with $s\in \R$,
we have
\begin{equation}\label{LfactLefR}
\int'_{W} \frac{ \Trace (\pi(H^m,u))\vert u \vert_\H^{ z}} {\vert
1-u \vert_\H} d^* u =
 -2\, \frac{d}{ds}\Im \log L_\R(H^m,z) \, .
\end{equation}
\end{The}

\proof The discussion above shows that we are reduced to the case of
even $m=2p$ and the subspace $H^{p,p}$. Let then $k= h^{p+}-
h^{p-}$. When $k=0$ the discussion above already gives the result.
When $k\neq 0$ the left hand side gets an additional term from the
trace of elements of the form $u=wj\in W$. Their contribution is
given by
$$
\frac{k}{2}\,\int_{\R_+^* } \frac{ v^{z-p}}{ 1+v } d^* v ,
$$
where the factor $\frac{1}{2}$ comes from the normalization (as a
probability measure)  of the fiber measure of the module $W \mapsto
{\R_+^* }$. The term in $v^{-p}$ comes from the action of the scalar
part $w\in \C^*$ with $w \overline{w}= v$ in $u=wj$. The denominator
comes from the equality
$$
\vert 1-u \vert_\H =\, 1 + \vert w \vert^2 \,,\quad u=w j\,.
$$

The other side changes in the same way using (\ref{trans}). Indeed
one has
$$
h^{p+}\, \log \Gamma_\R (z-p)+\,h^{p-}\, \log \Gamma_\R (z-p+1)=
$$
$$
\, \frac{h^{p+}+ h^{p-}}{2}\,\log (\Gamma_\R (z-p)\Gamma_\R
(z-p+1))+ \frac{k}{2}\,\log (\Gamma_\R (z-p)/\Gamma_\R (z-p+1))
$$
$$
= \,\frac{h^{p+}+ h^{p-}}{2}\,\log \Gamma_\C (z-p)+
\frac{k}{2}\,\log (\Gamma_\R (z-p)/\Gamma_\R (z-p+1))\,.
$$
This completes the proof.

\endproof

\medskip

The space which gives the formula of theorem \ref{real} as a
Lefschetz local contribution is obtained by taking the quaternions
$\H$ as the base, but endowing them with an additional structure,
namely their complex structure when viewed as a right vector space
over $\C$. We will use  the Atiyah-Bott Lefschetz formula (\cf
\cite{atbo}) applied to the $\bar{\partial}$-complex, which
generates the crucial term of theorem \ref{real}, that is,
\begin{equation}
\frac{ 1} {\vert 1-u \vert_\H}\,,
\end{equation}
with the reduced norm $\vert u \vert_\H$ as above, while a more
naive approach without the use of an elliptic complex would involve
the square of the reduced norm. The  Atiyah-Bott Lefschetz formula
 involves a numerator $\chi(u)$ which  yields,
\begin{equation}
\frac{\chi(u)} {\vert 1-u \vert_\H^2}=\frac{ 1} {\vert 1-u
\vert_\H}\,.
\end{equation}
In fact, the use of the $\bar{\partial}$-complex brings in the
powers $\wedge^jT_\C (\H)$ of the complex tangent space $T_\C (\H)$
of the complex manifold $\H$ and the alternating sum
\begin{equation}
\chi(u)=\,\sum_j (-1)^j\,\chi_j(u)=\,\sum_j (-1)^j\,{\rm
trace}(\wedge^j (u)) .
\end{equation}
This gives the determinant of the quaternion $1-u\in \H$ viewed as a
two by two complex matrix (even if some complex conjugates are
involved) and
\begin{equation}
\chi(u)=\,\vert 1-u \vert_\H .
\end{equation}

\medskip
\subsection{The question of the spectral realization}\label{SQuest}

A more general problem will be to obtain for archimedian local
factors a trace formula analogous to the semilocal case of Theorem 4
of \cite{AC}. This means that one considers a number field $K$, with
$S$ be the set of all archimedean places, and $X$ a non-singular
projective variety over $\K$. One would like to obtain the real part
of the logarithmic derivative of the full archimedian factor
$L_S(H^m, z)=\prod_{v\in S} L_{\K_v}(H^m, z)$ on the critical line
as a trace formula for the action of a suitable Weil group on a
complex space.

For a single place, one can work with a vector bundle over a base
space that is $B=\C$ for $v$ complex and $B=\H$ for $v$ real, with
fiber the $\Z$-graded vector space $E^{(m)}$ given by the cohomology
$H^m(X_v,\C)$, with $X_v$ the compact K\"ahler manifold determined
by the embedding $v:\K \hookrightarrow \C$. One expects a trace
formula of the following sort, analogous to the semilocal case of
\cite{AC}, modelled on the trace formula of \cite{gui}.

\begin{Con}\label{conjtrace}
Let $h \in \cS (\R_+^*)$ with compact support be viewed as an
element of $\cS(W)$ by composition with the module. Then, for
$\Lambda \to \infty$, one has
\begin{equation}\label{traceform}
\Trace (R_{\Lambda} \, \vartheta_a(h)) = 2h (1) B_m  \log \Lambda +
\sum_{v \in S} \int'_{W_v} \frac{h(\vert u
\vert^{-1})\,\Trace(\pi_v(H^m,u)) }{ \vert 1-u \vert_{\H_v}} \, d^*
u + o(1)
\end{equation}
where $B_m $ is the $m$-th Betti number of $X$, and  $\int'$ is
uniquely determined by the pairing with the unique distribution on
$\K_v$ which agrees with $\frac{du }{ \vert 1-u \vert}$ for $u \not=
1$ and whose Fourier transform relative to $\alpha_v$ vanishes at
$1$.
\end{Con}

For example, the above formula can be proved in simplest case of a
single complex place. The bundle $E=\oplus_m E^{(m)}$ comes endowed
with a representation of $\C^*$, by
\begin{equation}\label{CrepE}
\lambda: (z,\omega) \mapsto (\lambda z, \lambda^{-p}
\bar\lambda^{-q} \omega),
\end{equation}
for $\omega \in H^{p,q}(X_v,\C)$. We then let $\cH=L^2(B,E^{(m)})$
be the Hilbert space of $L^2$-sections of $E^{(m)}$ (for its
hermitian metric of trivial bundle). The action of $W=\C^*$ on $\cH$
is given by
\begin{equation}\label{repE1}
(\vartheta_a(\lambda)\xi)(b)=\,\lambda^{p} \bar\lambda^{q}
 \xi(\lambda^{-1} b)\,,\quad \forall \xi \in L^2(B,H^{p,q})\,,
\end{equation}
where we identify sections of $E^{(m)}$ with $H^m(X_v,\C)$ valued
functions on $B$.

One can use a cutoff as in \cite{AC}, by taking orthogonal
projection $P_{\Lambda}$ onto the subspace
\begin{equation}\label{cut}
P_{\Lambda} = \{ \xi \in L^2(B,E^{(m)}) \, ; \ \xi (b) = 0\,, \;
\forall b\in B \, ,\; \vert b \vert_\C > \Lambda \} \, .
\end{equation}
Thus, $P_{\Lambda}$ is the multiplication operator by the function
$\rho_{\Lambda}$, where $\rho_{\Lambda} (b) = 1$ if $\vert b
\vert_\C \leq \Lambda$, and $\rho_{\Lambda} (b) = 0$ for $\vert b
\vert_\C > \Lambda$. This gives an infrared cutoff and to get an
ultraviolet cutoff we use $\widehat{P}_{\Lambda} = F P_{\Lambda}
F^{-1}$ where $F$ is the Fourier transform which depends upon the
choice of the basic character $\alpha_v$ for the place $v$. We let
$$
R_{\Lambda} = \widehat{P}_{\Lambda} \, P_{\Lambda} \, .
$$

\begin{Pro}\label{traceform1C}
For the set of places $S$ consisting of a single complex place the
trace formula \eqref{traceform} holds.
\end{Pro}

\proof Both sides of the formula are additive functions of the
representation of $\C^*$ in $H^m(X_v,\C)$. We can thus assume that
this representation corresponds to a one dimensional $H^{p,q}$. Let
then $h_1(\lambda)=\,\lambda^{p} \bar\lambda^{q} \, h(\vert \lambda
\vert)$, one has
$$
h_1(u^{-1}) =\,h(\vert u \vert^{-1})\,\Trace(\pi_v(H^m,u))\,, \ \ \
\ \forall u\in \C^*
$$
while by \eqref{repE1} we get that $\vartheta_a(h)$ is the same
operator as $ U(h_1)$ in the notations of \cite{AC}. Thus applying
Theorem 4 of \cite{AC} to $h_1$ gives the desired result.
\endproof

In discussing here the case of a single complex place $v$, we have
taken the trivial bundle over $\C$ with fiber the cohomology
$H^m(X_v,\C)$. Already in order to treat the case of several complex
places, one needs to use the following result (\cite{Se2}, proof of
Proposition 12): the integers $h^{p,q}(X_v)$ are independent of the
archimedean place $v\in S$. This suggests that it is in fact more
convenient to think of the fiber $E^{(m)}$ as the motive $h^m(X)$
rather than its realization. This leads naturally to the further
question of a semi-local trace formula where the finite set of
places $S$ involves both archimedean and non-archimedean places. We
discuss in section \ref{Scurve} below, in the case of curves, how
one can think of a replacement for the $\ell$-adic cohomology at the
non-archimedean places, using noncommutative geometry.

\smallskip
\begin{Rem} {\rm The representation $\vartheta_a$ of Problem \ref{conjtrace} is not unitary
but  the product
\begin{equation}\label{scaleVrepbis}
u\mapsto \vert u \vert_W^{-(1+m)/2}\vartheta_a(u)
\end{equation}
should be unitary. In particular
 the spectral projection $N_E$, for the scaling action,
 associated to
the interval $[-E,E]$ in the dual group $\R$ of $\R_+^*$  is then
given by
\begin{equation}\label{NEproj3bis}
N_E=\vartheta_a(h^{(m)}_E),  \ \ \text{ with } \  h^{(m)}_E (u) =
\vert u \vert_W^{-(1+m)/2} \frac{1}{2\pi} \int_{-E}^E \vert u
\vert_W^{is} ds \, .
\end{equation}
Applying the conjectured formula \eqref{traceform} to the function
$h^{(m)}_E$ the left hand side gives the counting of quantum states,
$\Trace (R_{\Lambda} \, N_E) $ as in \cite{AC}. The first term in
the right hand side of \eqref{traceform} gives the contribution of
the regular representation of the scaling group. Finally, using
\eqref{Lefschetzloc} and \eqref{countLzerosfinal}, the last terms of
\eqref{traceform} combine to give
\begin{small}
$$
  \sum_{v|\infty
} \int'_{W_v} \frac{h^{(m)}_E( u^{-1})\,\Trace(\pi_v(H^m,u)) }{
\vert 1-u \vert_{\H_v}} \, d^* u  =\sum_{v|\infty
}\frac{1}{2\pi}\int_{-E}^E\,
\int'_{W_v}\frac{\Trace(\pi_v(H^m,u))\vert u
\vert_{W_v}^{\frac{1+m}{2}+is}}{ \vert 1-u \vert_{\H_v}} \, d^* u\,
ds$$
$$
  -\,\sum_{v|\infty
}\frac{1}{\pi}\,\int_{-E}^E\,\frac{d}{ds}\Im \log
L_v(H^m(X),\frac{1+m}{2}+is) ds =-\langle N_s(E)\rangle
$$
\end{small}
which shows that one should expect that the zeros of the
$L$-function appear as an absorption spectrum as in \cite{AC}.
Unlike in the case of the Riemann zeta function this remains
conjectural in the above case of $L$-functions of arithmetic
varieties.}
\end{Rem}

\medskip
\subsection{Local factors for curves}\label{Scurve}

In the previous part of this section, we showed how to
write the archimedean local factors of the $L$-function of an algebraic
variety $X$ in the form of a Lefschetz trace formula. Eventually, one
would like to obtain a semi-local formula, like the one conjectured
above, not just for the archimeden places, but for the full $L$-function.

We are only making some rather speculative remarks at this point,
and we simply want to illustrate the type of construction one expects should
give the geometric side of such a Lefschetz trace formula. We just discuss the
case where $X$ is a curve for simplicity. In this case we can concentrate on
the $L$-function for $H^1(X)$.

Let $X$ be a curve over a number field $\K$.
The local Euler factor at a place
$v$ has the following description (\cf \cite{Se}):
\begin{equation}\label{L-factorH1}
 L_v (H^1(X),z)= \det\left( 1-Fr_v^* N(v)^{-z} | H^1(\bar X,
\Q_\ell)^{I_v} \right)^{-1}.
\end{equation}
Here $Fr_v^*$ is the geometric Frobenius acting on $\ell$-adic
cohomology of $\bar X=X \otimes \Sp(\bar \K)$, with $\bar \K$ an
algebraic closure and $\ell$ a prime with $(\ell,q)=1$, where $q$ is
the cardinality of the residue field $k_v$ at $v$. We denote by $N$
the norm map. The determinant is evaluated on the inertia invariants
$H^1(\bar X, \Q_\ell)^{I_v}$ at $v$ (this is all of $H^1(\bar X,
\Q_\ell)$ when $v$ is a place of good reduction). The $L$-function
has then the Euler product description
$$ L(H^1(X),z) =\prod_v L_v (H^1(X),z), $$
where for $v$ a non-archimedean place the local factor is of the
form \eqref{L-factorH1} and if $v$ is a complex or real archimedean
place then the local factor is given by the corresponding
$\Gamma$-factor as discussed in the previous part of this section.

Usually, in using the expression \eqref{L-factorH1}  for the local
factors, one makes use of a choice of an embedding of $\Q_\ell$ in
$\C$. In our setting, if one wants to obtain a semi-local trace
formula for the $L$-function $L(H^1(X),z)$, one needs a geometric
construction which does not depend on such choices and treats the
archimedean and non-archimedean places on equal footing.  One
expects that a geometric space on which the geometric side of the
desired Lefschetz trace formula will concentrate should be obtained
as a fibration, where the base space should be a noncommutative
space obtained from the adele class space of the number field $\K$,
modified by considering, at the places of bad reduction and at the
real archimedean places, suitable division algebras over the local
field. The fiber should also be a noncommutative space in which the
special fiber embeds (at least at the places of good reduction).
This will have the effect of replacing the use of the $\ell$-adic
cohomology and the need for a choice of an embedding of $\Q_\ell$ in
$\C$.

In the case of a curve $X$, one can obtain this by embedding the
special fiber $X_v$ at a place $v$, which is a curve over the
residue field $k_v$ of cardinality $q$, in the adele class space of
the function field of $X_v$. Indeed by Theorem \ref{spectral}
applied to the global  field $\K$ of functions on $X_v$, one obtains
the local factor $L_v (H^1(X),z)$  directly  over $\C$. Thus, at
least in the case of curves,  the adele class spaces of the global
fields of functions on the curves $X_v$ should be essential
ingredients in the construction. This would then make it possible to
work everywhere with a cohomology defined over $\C$, in the form of
cyclic (co)homology.


\bigskip



\begin{thebibliography}{99}

\bibitem{Andre} Y.~Andr\'e, {\em Une introduction aux motifs},
Panoramas et Synth\`eses, Vol.17, Soci\'et\'e math\'ematique de
France, 2005.

\bibitem{atbo} M.F. Atiyah and R. Bott, {\it A Lefschetz fixed
point formula for elliptic complexes: I}, Annals of Math, Vol.86
(1967) 374--407.

\bibitem{BaJu} S.~Baaj, P.~Julg, {\em
Th\'eorie bivariante de Kasparov et op\'erateurs non born\'es dans
les $C\sp{*} $-modules hilbertiens}. C. R. Acad. Sci. Paris S‰r. I
Math. 296 (1983), no. 21, 875--878.

\bibitem{Baum} P.~Baum, {\em Cycles, cocycles and $K$-theory},
unpublished manuscript.

\bibitem{Bla} B.~Blackadar, {\em $K$-theory for operator algebras},
Second edition, MSRI Publications, 1998.

\bibitem{EB} E.~Bombieri, {\em Problems of the Millenium:
 The Riemann Hypothesis}, Clay Mathematical Institute (2000).

\bibitem{BC} J.B.~Bost, A.~Connes, {\em Hecke algebras, type III
factors and phase transitions with spontaneous symmetry breaking in
number theory}.  Selecta Math. (N.S.)  1  (1995),  no. 3, 411--457.

\bibitem{BR}  O.~Bratteli, D.W.~Robinson.
{\em Operator Algebras and Quantum Statistical Mechanics} I and II,
Springer Verlag, 1979 and 1981.

\bibitem{Bruhat} F.~Bruhat, {\em Distributions sur un groupe 
localement compact et
applications \`a l'\'etude des repr\'esentations des groupes
$p$-adiques}.   Bull. Soc. Math. France,   89  (1961), 43--75.

\bibitem{Co-th} A.~Connes, {\em Une classification des facteurs de
type ${\rm III}$}.  Ann. Sci. \'Ecole Norm. Sup. (4)  6  (1973),
133--252.

\bibitem{CoExt} A.~Connes, {\em  Cohomologie cyclique et foncteurs
${\rm Ext}\sp n$}.  C. R. Acad. Sci. Paris S\'er. I Math. 296
(1983), no. 23, 953--958.

\bibitem{CoIHES} A.~Connes, {\em Noncommutative differential
geometry}.  Inst. Hautes \'Etudes Sci. Publ. Math.  No. 62 (1985),
257--360.

\bibitem{Co94} A.~Connes, {\em Noncommutative geometry}. Academic
Press, 1994.

\bibitem{AC} A.~Connes, {\em Trace formula in noncommutative
geometry and zeros of the Riemann zeta function}, Selecta Math.
(N.S.) 5 (1999), no. 5, 29--106.

\bibitem{CoCM} A.~Connes, C.~Consani, M.~Marcolli, {\em The Weil proof
and the geometry of the adeles class space}, preprint.

\bibitem{CoKa} A.~Connes, M.~Karoubi, {\em Caract\`ere multiplicatif
d'un module de Fredholm}. $K$-Theory  2  (1988),  no. 3, 431--463.

\bibitem{CM} A.~Connes, M.~Marcolli, {\em From physics to number
theory via noncommutative geometry. Part I. Quantum statistical
mechanics of $\Q$-lattices}, in ``Frontiers in Number Theory, Physics,
and Geometry, I" pp.269--350, Springer Verlag, 2006.  

\bibitem{CM2} A.~Connes, M.~Marcolli, {\em From physics to number
theory via noncommutative geometry. Part II. Renormalization, the
Riemann-Hilbert correspondence, and motivic Galois theory},
in ``Frontiers in Number Theory, Physics, and Geometry, II"
pp.617--713, Springer Verlag, 2006. 

\bibitem{cmln} A.~Connes, M.~Marcolli, {\em Renormalization and
motivic Galois theory}, IMRN (2004) N.76, 4073--4092.

\bibitem{CManom} A.~Connes, M.~Marcolli, {\em Dimensional
regularization, anomalies, and noncommutative geometry}, in
preparation.

\bibitem{CMR} A.~Connes, M.~Marcolli, N.~Ramachandran, {\em KMS states
and complex multiplication}, Selecta Math. (N.S.)  11  (2005),
no. 3-4, 325--347. 

\bibitem{CMR2} A.~Connes, M.~Marcolli, N.~Ramachandran, {\em KMS states
and complex multiplication. Part II}, in ``Operator Algebras: The Abel
Symposium 2004'', pp.15--59, Abel Symp., 1, Springer Verlag, 2006. 

\bibitem{CMo} A.~Connes, H.~Moscovici,
 {\em  Hopf Algebras, Cyclic Cohomology and the Transverse
Index Theorem}, Commun. Math. Phys., 198, (1998), 199--246.

\bibitem{CoSka} A.~Connes, G.~Skandalis, {\em The longitudinal index
theorem for foliations}, Publ. RIMS Kyoto Univ. 20 (1984)
1139--1183.

\bibitem{ct} A.~Connes, M.~Takesaki, {\it The flow of weights on
factors of type III}. Tohoku Math. J., 29, (1977) 473--575.

\bibitem{ConsMar} C.~Consani, M.~Marcolli, {\em Noncommutative
geometry, dynamics, and $\infty$-adic Arakelov geometry}.  Selecta
Math. (N.S.)  10  (2004),  no. 2, 167--251.

\bibitem{ConsMar2} C.~Consani, M.~Marcolli, {\em New perspectives
in Arakelov geometry}.  Number theory,  81--102, CRM Proc.
Lecture Notes, 36, Amer. Math. Soc., Providence, RI, 2004.

\bibitem{De} P.~Deligne, {\em Valeurs de fonctions $L$ et p\'eriodes
d'int\'egrales}, Proc. Symp. Pure Math. Vol.33 (1979) part II,
313--346.

\bibitem{Del3} P.~Deligne, {\em \`A quoi servent les motifs?} in
``Motives'' (Seattle, WA, 1991), 143--161, Proc. Sympos .
  Pure Math., 55, Part 1, Amer. Math. Soc., Providence, RI, 1994.

\bibitem{Den1} C.~Deninger, {\it On the $\Gamma$-factors attached to
motives}, Invent. Math. 104 (1991) 245--261.

\bibitem{Den2} C.~Deninger, {\it Motivic $L$-functions and regularized
determinants}, in ``Motives'', Proceedings of Symposia in Pure
Mathematics, Vol. 55 (1994) Part I, 707--743.

\bibitem{dixmier} J.~Dixmier, {\em Les $C^*$-alg\`ebres et leurs
repr\'esentations}. Gauthier-Villars, Paris, 1964.

\bibitem{gui} V. Guillemin, {\it Lectures on spectral theory of
elliptic operators}, Duke Math. J., Vol.44  (1977) N.3, 485--517.

\bibitem{HaPau} E.~Ha, F.~Paugam, {\em A Bost--Connes--Marcolli system
for Shimura varieties}, IMRP 5 (2005) 237--286.

\bibitem{[Hilskand} M. Hilsum and G. Skandalis, {\em Morphismes
$K$-orient\'es d'espaces de feuilles et fonctorialit\'e en th\'eorie
de Kasparov}.   Ann. Sci. \'Ecole Norm. Sup. (4)
 20 (1987), 325-390.

\bibitem{Jacob} B.~Jacob, {\em Bost--Connes type systems for function
fields}, J. Noncommutative Geometry, Vol.1 (2007) N.2, 141--211. 

\bibitem{Jannsen} U.~Jannsen, {\em Motives, numerical equivalence,
and semi-simplicity}.
Invent. Math. 107 (1992), no. 3, 447--452.

\bibitem{Kal} D.~Kaledin, {\em Non-commutative Cartier
operator and Hodge-to-de Rham degeneration}, preprint
math.AG/0511665.

\bibitem{Kas} G.G.~Kasparov, {\em The operator $K$-functor and
extensions of $C\sp{*} $-algebras}. Izv. Akad. Nauk SSSR Ser. Mat.
44 (1980), no. 3, 571--636, 719.

\bibitem{Maxim} M.~Kontsevich, {\em Noncommutative motives}, talk
at Princeton, October 2005.

\bibitem{lr1} M. Laca and I. Raeburn, {\em Semigroup crossed
products  and the Toeplitz algebras of nonabelian groups}, J. Funct.
Anal., {\bf 139} (1996), 415--440.

\bibitem{l1} M. Laca, {\em From endomorphisms to automorphisms and
back: dilations and full  corners},  J. London Math. Soc. (2)  {\bf
61}  (2000),  no. 3, 893--904.

\bibitem{Lo} J.L.~Loday, {\em Cyclic homology}. Grundlehren der
Mathematischen Wissenschaften, 301. Springer-Verlag, Berlin, 1998.
xx+513 pp.

\bibitem{Man} Yu.I.~Manin, {\em Lectures on the $K$-functor in
algebraic geometry}. Uspehi Mat. Nauk  24  1969 no. 5 (149), 3--86.

\bibitem{Man-mot} Yu.I.~Manin, {\em Correspondences, motifs, and
monoidal transformations}, Mat. Sb. (N.S.) 77 (1968) 119, 475--507.

\bibitem{Meyer} R.~Meyer, {\em On a representation of the idele class
group related to primes and zeros of $L$-functions}, preprint arXiv
math.NT/0311468.

\bibitem{MeyerNest} R.~Meyer, R.~Nest,
{\em The Baum-Connes conjecture via localization of categories}.
Lett. Math. Phys. 69 (2004), 237--263.

\bibitem{[R]} B. Riemann, {\it Mathematical Works}, Dover,
New York, 1953.

\bibitem{Se2} J.~P.~Serre, {\it G\'eom\'etrie alg\'ebrique
et g\'eom\'etrie analytique}. Annales Inst. Fourier, Grenoble 6,
(1956), 1-42.

\bibitem{Se} J.~P.~Serre, {\it Facteurs locaux des fonctions z\^eta
des vari\'et\'es alg\'ebriques (d\'efinitions et conjectures)}.
S\'em. Delange-Pisot-Poitou, exp.~19, 1969/70.

\bibitem{tt}  M.~Takesaki, {\it Tomita's theory of modular
Hilbert algebras and its applications}. Lecture Notes in Math., 28,
Springer, 1970.

\bibitem{Tak} M. Takesaki, {\em  Duality for crossed products and the structure
of von Neumann algebras of type III},  Acta Math. (131) (1973),
249-310.

\bibitem{Ta} J.~Tate, {\em Number theoretic background}, in ``Automorphic
forms, representations, and $L$-functions'', Proc. Symp. Pure Math.
Vol.33, Part 2, 1979, pp.3--26.

\bibitem{We} A.~Weil, {\em Sur la th\'eorie du corps de classes},
J. Math. Soc. Japan, 3 (1951) 1--35.

\bibitem{weilpos} A.~Weil, {\it Sur les formules explicites de la
th\'eorie des nombres premiers}, Oeuvres complètes, Vol. 2, 48--62.

\bibitem{weil} A. Weil, {\it Sur les formules explicites de la
th\'eorie des nombres}, Izv. Mat. Nauk., (Ser. Mat.) Vol.36 (1972)
3--18.

\end{thebibliography}
\end{document}